\pdfoutput=1
\PassOptionsToPackage{pdfusetitle}{hyperref}
\PassOptionsToPackage{capitalize}{cleveref}
\PassOptionsToPackage{style=authoryear-comp}{biblatex}
\documentclass[a4paper]{article}

\usepackage{xr}
\usepackage{subfiles}
\usepackage{biblatex}
\DeclareNameAlias{author}{given-family}
\newcommand{\myCite}{\parencite}
\newcommand{\myTextCite}{\textcite}
\usepackage{enumitem}
\usepackage{amsthm}
\usepackage{my-preamble}
\usepackage{xassoccnt}

\makeatletter
\newcommand{\counter@within}{section}
\DeclareCoupledCountersGroup{thm@group}
\newcommand{\myNewTheorem}[2]{%
  \newtheorem{#1}{#2}[\counter@within]
  \AddCoupledCounters[name=thm@group]{#1}
  \Crefname{#1}{#2}{#2s}
}
\makeatother
\theoremstyle{definition}
\myNewTheorem{myDefinition}{Definition}
\myNewTheorem{myConstruction}{Construction}
\myNewTheorem{myNotation}{Notation}
\theoremstyle{remark}
\myNewTheorem{myExample}{Example}
\myNewTheorem{myRemark}{Remark}
\theoremstyle{plain}
\myNewTheorem{myProposition}{Proposition}
\myNewTheorem{myCorollary}{Corollary}
\myNewTheorem{myLemma}{Lemma}
\myNewTheorem{myTheorem}{Theorem}

\title{\myTitle}
\author{\myAuthor}

\begin{document}

\maketitle

\begin{abstract}

We propose elementary definitions of opetopes and opetopic sets. We
directly define opetopic sets by a simple structure and several
axioms. Opetopes are then opetopic sets satisfying one more axiom. We
show that our definition is equivalent to the polynomial monad
definition given by Kock, Joyal, Batanin, and Mascari. We also show
that our category of opetopes is equivalent to the one given by Ho
Thanh.

\end{abstract}

\section{Introduction}
\label{sec:introduction}

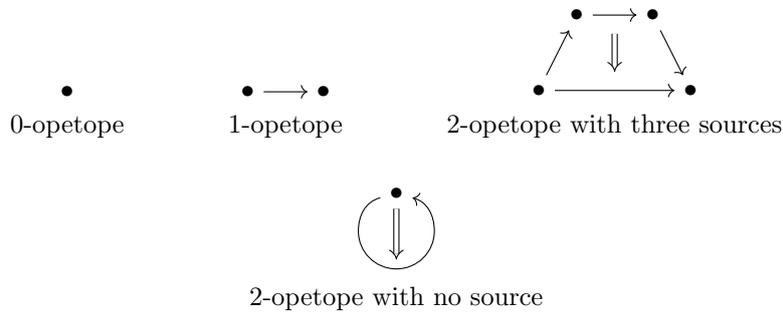
\begin{figure}
  \begin{mathpar}
    \begin{tikzpicture}
      \node (x) at (0, 0) {\(\bullet\)};
      \node[below] (title) at (current bounding box.south) {\(0\)-opetope};
    \end{tikzpicture}
    \and
    \begin{tikzpicture}
      \node (x) at (0, 0) {\(\bullet\)};
      \node (y) at (1, 0) {\(\bullet\)};
      \draw[->] (x) to (y);
      \node[below] (title) at (current bounding box.south) {\(1\)-opetope};
    \end{tikzpicture}
    \and
    \begin{tikzpicture}
      \node (x0) at (0, 0) {\(\bullet\)};
      \node (x1) at (0.5, 1) {\(\bullet\)};
      \node (x2) at (1.5, 1) {\(\bullet\)};
      \node (x3) at (2, 0) {\(\bullet\)};
      \draw[->] (x0) to (x1);
      \draw[->] (x1) to node[below,name=y0] {} (x2);
      \draw[->] (x2) to (x3);
      \draw[->] (x0) to node[above,name=y1] {} (x3);
      \draw[double equal sign distance,-Implies] (y0) to (y1);
      \node[below] (title) at (current bounding box.south)
      {\(2\)-opetope with three sources};
    \end{tikzpicture}
    \and
    \begin{tikzpicture}
      \node (x0) at (0.5, 1) {\(\bullet\)};
      \draw[->] (x0)
      to [out=195,in=90] (0, 0.5)
      to [out=-90,in=180] (0.5, 0)
      to [out=0,in=-90] (1, 0.5)
      to [out=90,in=-15] (x0);
      \node (x1) at (0.5, 0) {};
      \draw[double equal sign distance,-Implies] (x0) to (x1);
      \node[below] (title) at (current bounding box.south)
      {\(2\)-opetope with no source};
    \end{tikzpicture}
  \end{mathpar}
  \caption{Examples of opetopes of dimension \(0\), \(1\), and \(2\)}
  \label{fig-opetopes}
\end{figure}

\begin{figure}
  \begin{mathpar}
    \begin{tikzpicture}
      \node (x0) at (0, 0) {\(\bullet\)};
      \node (x1) at (-0.5, 1) {\(\bullet\)};
      \node (x2) at (0.5, 1) {\(\bullet\)};
      \node (x3) at (1.5, 2.5) {\(\bullet\)};
      \node (x4) at (2.5, 1) {\(\bullet\)};
      \node (x5) at (3, 0) {\(\bullet\)};
      \draw[->] (x0) to (x1);
      \draw[->] (x0) to node[above,name=y0] {} (x2);
      \draw[->] (x0) to node[above,name=y1] {} (x5);
      \draw[->] (x1) to (x2);
      \draw[->,out=90,in=180] (x2) to (x3);
      \draw[->] (x2) to node[above,name=y2] {} node[below,name=y3] {} (x4);
      \draw[->] (x3)
      to [out=210,in=180] (1.5, 1.75)
      to node[above,name=y4] {} node[below,name=y5] {} (1.5, 1.75)
      to [out=0,in=330] (x3);
      \draw[->,out=0,in=90] (x3) to (x4);
      \draw[->] (x4) to node[above,name=y6] {} (x5);
      \draw[->,out=30,in=30] (x4) to node[below,name=y7] {} (x5);
      \draw[double equal sign distance,-Implies] (x1) to (y0);
      \draw[double equal sign distance,-Implies] (x3) to (y4);
      \draw[double equal sign distance,-Implies] (y5) to (y2);
      \draw[double equal sign distance,-Implies] (y7) to (y6);
      \draw[double equal sign distance,-Implies] (y3) to (y1);
    \end{tikzpicture}
    \and
    \begin{tikzpicture}
      \node (x0) at (0, 0) {\(\bullet\)};
      \node (x1) at (-0.5, 1) {\(\bullet\)};
      \node (x2) at (0.5, 1) {\(\bullet\)};
      \node (x3) at (1.5, 2.5) {\(\bullet\)};
      \node (x4) at (2.5, 1) {\(\bullet\)};
      \node (x5) at (3, 0) {\(\bullet\)};
      \draw[->] (x0) to (x1);
      \draw[->] (x0) to node[above,name=y1] {} (x5);
      \draw[->] (x1) to (x2);
      \draw[->,out=90,in=180] (x2) to (x3);
      \draw[->,out=0,in=90] (x3) to (x4);
      \draw[->,out=30,in=30] (x4) to node[below,name=y7] {} (x5);
      \draw[double equal sign distance,-Implies] (x3) to (y1);
    \end{tikzpicture}
  \end{mathpar}
  \caption{Example of opetope of dimension \(3\), the pasting diagram
    of sources on the left and the target on the right.}
  \label{fig-pasting-diagram}
\end{figure}
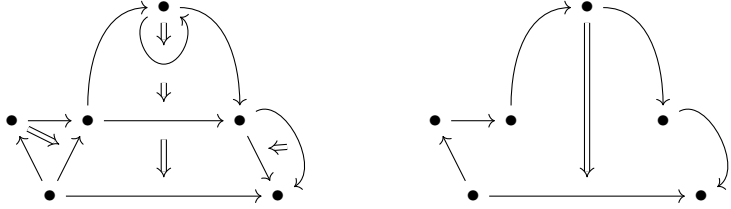

\emph{Opetopes} and \emph{optopic sets} were introduced by
\myTextCite{baez1998higher3} as a combinatorial approach to weak
\(\omega\)-categories. An opetope is a geometric shape of a
many-in-single-out operator in higher dimension. Examples of opetopes
of low dimensions are drawn in \cref{fig-opetopes}. There is only one
opetope of dimension \(0\), the point. There is only one opetope of
dimension \(1\), the arrow with one target and one source. There are
countably many opetopes of dimension \(2\). The sources of an opetope
of dimension \(2\) are opetopes of dimension \(1\) that form a
``pasting diagram''. It can be the case that an opetope of dimension
\(2\) has no source in which case it looks like an arrow filling a
loop. An opetope of dimension \(3\) is determined by its pasting
diagram of sources, and its target is the opetope of dimension \(2\)
that has the same ``boundary'' as the pasting diagram of sources; see
\cref{fig-pasting-diagram}. The opetopes form a category, and opetopic
sets are presheaves on the category of opetopes.

Several equivalent definitions of opetopes have been proposed:
\myTextCite{baez1998higher3} using operads;
\myTextCite{leinster2004higher} using cartesian monads;
\myTextCite{hermida2002weak3} (called multitopes there) using
multicategories. Comparison of these definitions is made by
\myTextCite{cheng2003opetopes,cheng2004weak,cheng2004multitopic}. More
recent accounts are given by \myTextCite{kock2010polynomial} using
polynomial monads and by \myTextCite{curien2022type} using type
theory. \myTextCite{ho-thanh2020equivalence-arxiv} gives an explicit
presentation of the category of opetopes by generators and relations
and shows that presheaves on it, that is, opetopic sets, are
equivalent to many-to-one polygraphs. Existing definitions of
opetopes, however, require some amount of prerequisites and are less
accessible than other notions of geometric shapes such as
simplices. What is missing is a concise definition of opetopes by a
structure and axioms just as simplices are inhabited finite total
orders.

In this paper, we propose elementary definitions of opetopes and
opetopic sets. We directly define opetopic sets by a simple structure
and several axioms. Our formal definition of opetopic sets given in
\cref{sec:defin-opet-sets} takes less than two pages, and the only
prerequisite is basic category theory. In \cref{sec:prop-opet-sets},
we show basic properties of opetopic sets and morphisms between
them. One goal is that every slice of the category of opetopic sets is
a presheaf category (\cref{prop-oset-slice-presheaf}).

Opetopes are defined in \cref{sec:opetopes} as opetopic sets
satisfying one more axiom. An interesting phenomenon is that the
category of opetopes is naturally turned into an opetopic
set. Moreover, the opetopic set of opetopes is the terminal object in
the category of opetopic sets. Consequently, the category of opetopic
sets is a presheaf category (\cref{prop-oset-presheaf}). In
particular, colimits of opetopic sets exist. We provide in
\cref{sec:categ-opet-sets} some tools to compute colimits of opetopic
sets.

In \cref{sec:bound-past-diagr}, we introduce \emph{boundaries} and
\emph{pasting diagrams}. We show that an opetope is completely
determined by its boundary (\cref{prop-boundary-equiv}) or the pasting
diagram formed by its sources (\cref{prop-opetope-equiv-pd}). We
provide in \cref{sec:subst-graft} two operators on pasting diagrams,
\emph{substitution} and \emph{grafting}. Since we already know the
presheaf category of opetopic sets, these operators are simply defined
by colimits.

Using the substitution and grafting operators, we show in
\cref{sec:baez-dolan-constr} that our definition of opetopes is
equivalent to the polynomial monad definition given by
\myTextCite{kock2010polynomial} (\cref{prop-ope-equiv-kock}). We also
see that the category of opetopes is presented by the generators and
relations described by \myTextCite{ho-thanh2020equivalence-arxiv}.

\paragraph{Foundations}

The results in the present paper are valid in any constructive
foundations of mathematics. For concreteness, we choose
\emph{Univalent Foundations} \myCite{hottbook} because it seems to be
a proper foundation especially for category theory in that isomorphic
objects are identical \myCite{ahrens2015univalent}. Only one univalent
universe \(\mathcal{U}\) is needed. An object is said to be \myDefine{small} if
it is equivalent to an object in \(\mathcal{U}\). We do not assume the law of
excluded middle, the axiom of choice, or the propositional resizing
axiom. We use notation \((x \myElemOf A) \myMorphism \myApp{B}{x}\)
for dependent function types (\(\Pi\)-types) and
\((x \myElemOf A) \myBinProd \myApp{B}{x}\) for dependent pair types
(\(\Sigma\)-types). A category in Univalent Foundations
\myCite{ahrens2015univalent} satisfies that the type of
identifications \(x \myId y\) between objects is equivalent to the
type of equivalences \(x \myEquiv y\), so equivalent objects satisfy
the same properties. \(\myApp{\myObj}{C}\) denotes the type of objects
in a category \(C\), and \(\myApp{\myArr_{C}}{x \myComma y}\) denotes
the set of arrows from \(x\) to \(y\) in \(C\).

All the types that appear in this paper are actually \(1\)-truncated,
so our language can be translated into set-theoretic foundations by
interpreting types as groupoids \myCite{hofmann1998groupoid}.

\section{Definition of opetopic sets}
\label{sec:defin-opet-sets}

We first introduce our formal definition of opetopic sets
(\cref{def-opetopic-set}) and then explain intuition.

\begin{myDefinition}
  We say a category \(C\) is \myDefine{gaunt} if its type of objects
  \(\myApp{\myObj}{C}\) is a set. The small gaunt categories and the
  functors between them form a category \(\myGauntCat\).
\end{myDefinition}

\begin{myRemark}
  In terms of set-theoretic foundations, a category is gaunt if its
  underlying groupoid is discrete. This is equivalent to that the
  identities are the only isomorphisms, which coincides with the
  definition given in \myCite[Definition 3.1]{barwick2021unicity}.
\end{myRemark}

\begin{myDefinition}
  An \myDefine{\(\myFinOrd\)-direct category} is a gaunt category
  \(A\) equipped with a conservative functor
  \(\myDegree_{A} \myElemOf A \myMorphism \myFinOrd\) called the
  \myDefine{degree functor}, where \(\myFinOrd\) is the poset of
  natural numbers. Let
  \(\myDirCat{\myFinOrd} \mySub \myGauntCat \mySlice \myFinOrd\)
  denote the full subcategory spanned by the \(\myFinOrd\)-direct
  categories.
\end{myDefinition}

\begin{myDefinition}
  Let \(k \myElemOf \myFinOrd\) and let \(A\) be an
  \(\myFinOrd\)-direct category. We write
  \(f \myElemOf x \myKMorphism{k} y\) to mean that
  \(f \myElemOf x \myMorphism y\) is an arrow in \(A\) between objects
  satisfying that
  \(\myApp{\myDegree}{x} + k \myId \myApp{\myDegree}{y}\). Such an
  arrow is called a \myDefine{\(k\)-step} arrow. Let
  \(\myApp{\myKArr{k}}{A}\) denote the set of \(k\)-step arrows in
  \(A\). We also define the \myDefine{\(k\)-step slice}
  \(A \myKSlice{k} x\) to be the full subcategory of \(A \mySlice x\)
  spanned by the \(k\)-step arrows into \(x\). Note that
  \(A \myKSlice{k} x\) is discrete.
\end{myDefinition}

\begin{myDefinition}
  A \myDefine{preopetopic set} is an \(\myFinOrd\)-direct category
  \(A\) equipped with a subset
  \(\myApp{\mySource}{A} \mySub \myApp{\myOneArr}{A}\) with complement
  \(\myApp{\myTarget}{A}\). An arrow in \(\myApp{\mySource}{A}\) is
  called a \myDefine{source arrow} and written as
  \(f \myElemOf x \mySourceMorphism y\). An arrow in
  \(\myApp{\myTarget}{A}\) is called \myDefine{target arrow} and
  written as \(f \myElemOf x \myTargetMorphism y\). A
  \myDefine{morphism of preopetopic sets} is a morphism of
  \(\myFinOrd\)-direct categories preserving source and target
  arrows. Let \(\myPreOSet\) denote the category of small preopetopic
  sets.
\end{myDefinition}

\begin{myDefinition}
  Let \(A\) be a preopetopic set and let
  \(f \myElemOf y \myOneMorphism x\) and
  \(g \myElemOf z \myOneMorphism y\) be \(1\)-step arrows in \(A\). We
  say \((f \myComma g)\) is \myDefine{homogeneous} if either both
  \(f\) and \(g\) are source arrows or both \(f\) and \(g\) are target
  arrows. We say \((f \myComma g)\) is \myDefine{heterogeneous} if
  either \(f\) is a source arrow and \(g\) is a target arrow or \(f\)
  is a target arrow and \(g\) is a source arrow. By a
  \myDefine{homogeneous/heterogeneous factorization} of a \(2\)-step
  arrow \(h\) we mean a factorization \(h \myId f \myComp g\) such
  that \((f \myComma g)\) is homogeneous/heterogeneous.
\end{myDefinition}

\begin{myDefinition}
  \label{def-opetopic-set}
  An \myDefine{opetopic set} is a preopetopic set \(A\) satisfying the
  following axioms.
  \begin{myWithCustomLabel}{axiom}
    \begin{enumerate}[label=\textbf{O\arabic*.},ref=O\arabic*]
    \item \label{ax-oset-first} \label{ax-oset-finite}
      \(A \myOneSlice x\) is finite for every \(x \myElemOf A\).
    \item \label{ax-oset-target} For every object \(x \myElemOf A\)
      of degree \(\myGe 1\), there exists a unique target arrow into
      \(x\).
    \item \label{ax-oset-source} For every object \(x \myElemOf A\)
      of degree \(1\), there exists a unique source arrow into \(x\).
    \item \label{ax-oset-homo} Every \(2\)-step arrow
      \(y \myTwoMorphism x\) in \(A\) has a unique homogeneous
      factorization.
    \item \label{ax-oset-hetero} Every \(2\)-step arrow
      \(y \myTwoMorphism x\) in \(A\) has a unique heterogeneous
      factorization.
    \item \label{ax-oset-path} For every object \(x \myElemOf A\)
      of degree \(\myGe 2\), there exists a \(2\)-step arrow
      \(r \myElemOf A \myTwoSlice x\) such that, for
      every \(2\)-step arrow \(f \myElemOf A \myTwoSlice x\),
      there exists a zigzag
      \begin{equation}
        \label{eq-oset-path}
        f \myId f_{0} \myXMorphism{s_{0}} g_{0} \myXMorphismAlt{t_{0}}
        f_{1} \myXMorphism{s_{1}} \myDots \myXMorphism{s_{m - 1}}
        g_{m - 1} \myXMorphismAlt{t_{m - 1}} f_{m} \myId r,
      \end{equation}
      where \(g_{i}\)'s are source arrows into \(x\), \(s_{i}\)'s are
      source arrows in \(A \mySlice x\), and \(t_{i}\)'s are target
      arrows in \(A \mySlice x\).
      \label{ax-oset-core-last}
    \item \label{ax-oset-kstep-unique} For every target arrow
      \(f \myElemOf y \myTargetMorphism x\) in \(A\) and object
      \(z \myElemOf A\) of degree \(\myLe \myApp{\myDegree}{y} - 2\),
      the postcomposition map
      \(f_{\myBang} \myElemOf \myApp{\myArr_{A}}{z \myComma y}
      \myMorphism \myApp{\myArr_{A}}{z \myComma x}\) is injective.
    \item \label{ax-oset-kstep-factor} For every \(k \myGe 3\),
      every \(k\)-step arrow \(y \myKMorphism{k} x\) in \(A\) factors as
      \(f \myComp g\) such that \(f\) is a \((k - 1)\)-step arrow and
      \(g\) is a \(1\)-step arrow.
      \label{ax-oset-last}
    \end{enumerate}
  \end{myWithCustomLabel}
  Let \(\myOSet \mySub \myPreOSet\) denote the full subcategory
  spanned by the opetopic sets.
\end{myDefinition}

Let \(A\) be an opetopic set. We think of objects in \(A\) as
\emph{cells}, and arrows in \(A\) determine the configuration of the
cells. A source arrow \(y \mySourceMorphism x\) exhibits \(y\) as a
source of \(x\), and a target arrow \(y \myTargetMorphism x\) exhibits
\(y\) as a target of \(x\). \Cref{ax-oset-finite} asserts that every
cell \(x\) has finitely many sources and targets. Recall that a set
\(X\) is finite if there (merely) exist a natural number
\(n \myElemOf \myNat\) and an equivalence
\(\mySetCompr{k \myElemOf \myNat}{k \myLt n} \myEquiv X\)
\myCite[Definition
16.3.1]{rijke2022introduction-arxiv}. \Cref{ax-oset-target} asserts
that every cell \(x\) of dimension \(\myGe 1\) has a unique target,
expressing the single-out nature of opetopes. An opetope may have many
sources with the exception of the opetope of dimension \(1\) which has
a unique source, so we introduce \cref{ax-oset-source}.

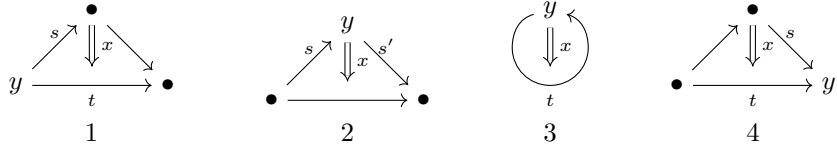
\begin{figure}
  \begin{mathpar}
    \begin{tikzpicture}
      \node (x0) at (0, 0) {\(y\)};
      \node (x1) at (1, 1) {\(\bullet\)};
      \node (x2) at (2, 0) {\(\bullet\)};
      \draw[->] (x0) to node[above] {\scriptsize \(s\)} (x1);
      \draw[->] (x0) to node[above,name=y0] {}
      node[below] {\scriptsize \(t\)} (x2);
      \draw[->] (x1) to (x2);
      \draw[double equal sign distance,-Implies] (x1) to node[right]
      {\scriptsize \(x\)} (y0);
      \node[below] (title) at (current bounding box.south)
      {\labelcref{item-ss-st}};
    \end{tikzpicture}
    \and
    \begin{tikzpicture}
      \node (x0) at (0, 0) {\(\bullet\)};
      \node (x1) at (1, 1) {\(y\)};
      \node (x2) at (2, 0) {\(\bullet\)};
      \draw[->] (x0) to node[above] {\scriptsize \(s\)} (x1);
      \draw[->] (x0) to node[above,name=y0] {} (x2);
      \draw[->] (x1) to node[above] {\scriptsize \(s'\)} (x2);
      \draw[double equal sign distance,-Implies] (x1) to node[right]
      {\scriptsize \(x\)} (y0);
      \node[below] (title) at (current bounding box.south)
      {\labelcref{item-ss-ts}};
    \end{tikzpicture}
    \and
    \begin{tikzpicture}
      \node (x0) at (0.5, 1) {\(y\)};
      \draw[->] (x0)
      to [out=195,in=90] (0, 0.5)
      to [out=-90,in=180] (0.5, 0)
      to node[above,name=x1] {} node[below] {\scriptsize \(t\)} (0.5, 0)
      to [out=0,in=-90] (1, 0.5)
      to [out=90,in=-15] (x0);
      \draw[double equal sign distance,-Implies] (x0) to node[right]
      {\scriptsize \(x\)} (x1);
      \node[below] (title) at (current bounding box.south)
      {\labelcref{item-tt-st}};
    \end{tikzpicture}
    \and
    \begin{tikzpicture}
      \node (x0) at (0, 0) {\(\bullet\)};
      \node (x1) at (1, 1) {\(\bullet\)};
      \node (x2) at (2, 0) {\(y\)};
      \draw[->] (x0) to (x1);
      \draw[->] (x0) to node[above,name=y0] {} node[below]
      {\scriptsize \(t\)} (x2);
      \draw[->] (x1) to node[above] {\scriptsize \(s\)} (x2);
      \draw[double equal sign distance,-Implies] (x1) to node[right]
      {\scriptsize \(x\)} (y0);
      \node[below] (title) at (current bounding box.south)
      {\labelcref{item-tt-ts}};
    \end{tikzpicture}
  \end{mathpar}
  \caption{Illustration of \cref{ax-oset-homo,ax-oset-hetero}. In Case
    \labelcref{item-ss-st}, \(y\) is a source of the source \(s\) of
    \(x\) and a source of the target \(t\) of \(x\). In Case
    \labelcref{item-ss-ts}, \(y\) is a source of the source \(s'\) of
    \(x\) and the target of the source \(s\) of \(x\). In Case
    \labelcref{item-tt-st}, \(y\) is the target of the target \(t\) of
    \(x\) and a source of the target \(t\) of \(x\). In Case
    \labelcref{item-tt-ts}, \(y\) is the target of the target \(t\) of
    \(x\) and the target of the source \(s\) of \(x\).}
  \label{fig-ax-oset-homo}
\end{figure}

\Cref{ax-oset-homo,ax-oset-hetero} assert that, for every \(2\)-step
arrow \(y \myTwoMorphism x\), exactly one of the following holds.
\begin{enumerate}
\item \label{item-ss-st} \(y\) is a source of a source of \(x\) and a
  source of the target of \(x\).
\item \label{item-ss-ts} \(y\) is a source of a source of \(x\) and
  the target of a source of \(x\).
\item \label{item-tt-st} \(y\) is the target of the target of \(x\)
  and a source of the target of \(x\).
\item \label{item-tt-ts} \(y\) is the target of the target of \(x\)
  and the target of a source of \(x\).
\end{enumerate}
\Cref{fig-ax-oset-homo} illustrates each of these situations. These
axioms are also related to the so-called diamond property in the
theory of abstract polytopes \myCite{mcmullen2002abstract}. Every
\(2\)-step arrow \(y \myTwoMorphism x\) factors in exactly two ways,
homogeneously and heterogeneously, forming a diamond shape.
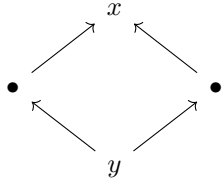
\begin{equation*}
  \begin{tikzcd}
    & x & \\
    \bullet
    \arrow[ur] & &
    \bullet
    \arrow[ul] \\
    & y
    \arrow[ul]
    \arrow[ur]
  \end{tikzcd}
\end{equation*}
Properties similar to \cref{ax-oset-homo,ax-oset-hetero} also appear
in \myCite[Definition 5]{hadzihasanovic2020combinatorial} and
\myCite[Theorem 1.36]{nguyen2018thesis}.

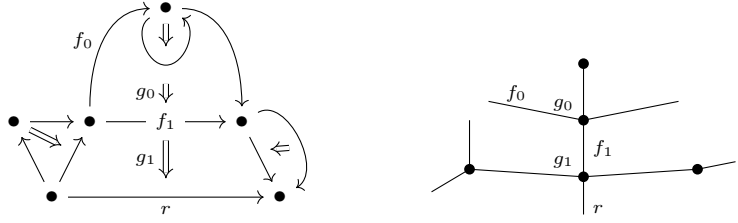
\begin{figure}
  \begin{mathpar}
      \begin{tikzpicture}
      \node (x0) at (0, 0) {\(\bullet\)};
      \node (x1) at (-0.5, 1) {\(\bullet\)};
      \node (x2) at (0.5, 1) {\(\bullet\)};
      \node (x3) at (1.5, 2.5) {\(\bullet\)};
      \node (x4) at (2.5, 1) {\(\bullet\)};
      \node (x5) at (3, 0) {\(\bullet\)};
      \draw[->] (x0) to (x1);
      \draw[->] (x0) to node[above,name=y0] {} (x2);
      \draw[->] (x0) to node[above,name=y1] {}
      node[below] {\scriptsize \(r\)} (x5);
      \draw[->] (x1) to (x2);
      \draw[->,out=90,in=180] (x2) to
      node[left] {\scriptsize \(f_{0}\)} (x3);
      \draw[->] (x2) to node[above,name=y2] {} node[below,name=y3] {}
      node[fill=white] {\scriptsize \(f_{1}\)} (x4);
      \draw[->] (x3)
      to [out=210,in=180] (1.5, 1.75)
      to node[above,name=y4] {} node[below,name=y5] {} (1.5, 1.75)
      to [out=0,in=330] (x3);
      \draw[->,out=0,in=90] (x3) to (x4);
      \draw[->] (x4) to node[above,name=y6] {} (x5);
      \draw[->,out=30,in=30] (x4) to node[below,name=y7] {} (x5);
      \draw[double equal sign distance,-Implies] (x1) to (y0);
      \draw[double equal sign distance,-Implies] (x3) to (y4);
      \draw[double equal sign distance,-Implies] (y5) to
      node[left] {\scriptsize \(g_{0}\)} (y2);
      \draw[double equal sign distance,-Implies] (y7) to (y6);
      \draw[double equal sign distance,-Implies] (y3) to
      node[left] {\scriptsize \(g_{1}\)} (y1);
    \end{tikzpicture}
    \and
    \begin{tikzpicture}
      \node (x0) at (1.5, 0) {};
      \node (x1) at (1.5, 0.5) {};
      \fill (x1) circle (0.2em) node[above left] {\scriptsize \(g_{1}\)};
      \node (x2) at (0, 0.6) {};
      \fill (x2) circle (0.2em);
      \node (x3) at (-0.5, 0.3) {};
      \node (x4) at (0, 1.25) {};
      \node (x5) at (1.5, 1.25) {};
      \fill (x5) circle (0.2em) node[above left] {\scriptsize \(g_{0}\)};
      \node (x6) at (0.25, 1.5) {};
      \node (x7) at (1.5, 2) {};
      \fill (x7) circle (0.2em);
      \node (x8) at (2.75, 1.5) {};
      \node (x9) at (3, 0.6) {};
      \fill (x9) circle (0.2em);
      \node (x10) at (3.5, 0.7) {};
      \draw (x0.center) to
      node[below right] {\scriptsize \(r\)} (x1.center);
      \draw (x1.center) to (x2.center);
      \draw (x2.center) to (x3.center);
      \draw (x2.center) to (x4.center);
      \draw (x1.center) to
      node[right] {\scriptsize \(f_{1}\)} (x5.center);
      \draw (x5.center) to
      node[above left] {\scriptsize \(f_{0}\)} (x6.center);
      \draw (x5.center) to (x7.center);
      \draw (x5.center) to (x8.center);
      \draw (x1.center) to (x9.center);
      \draw (x9.center) to (x10.center);
    \end{tikzpicture}
  \end{mathpar}
  \caption{Illustration of \cref{ax-oset-path}. The pasting diagram on
    the left has the tree structure on the right. Dots and lines in
    the tree correspond to \(2\)-dimensional cells and
    \(1\)-dimensional cells, respectively, in the pasting diagram.}
  \label{fig-pd-tree}
\end{figure}

\Cref{ax-oset-path}, combined with other axioms, expresses that the
sources of a cell form a ``tree''. For example, consider the pasting
diagram on the left of \cref{fig-pd-tree} which is the sources of a
\(3\)-dimensional cell. We see the zigzag
\begin{equation*}
  f_{0} \mySourceMorphism g_{0} \myTargetMorphismAlt f_{1}
  \mySourceMorphism g_{1} \myTargetMorphismAlt r,
\end{equation*}
and one can find a similar zigzag from any \(1\)-cell to
\(r\). Moreover, such a zigzag to \(r\) is unique; we will prove this
in \cref{prop-oset-tree} later. In this way we see the tree structure
of the pasting diagram displayed on the right of \cref{fig-pd-tree}.

\Crefrange{ax-oset-first}{ax-oset-core-last} are local conditions in
that they mention only \(1\)-step and \(2\)-step
arrows. \Cref{ax-oset-kstep-unique,ax-oset-kstep-factor} are global
conditions for tying cells in all dimensions together. There can be
several axiomatizations. We choose
\cref{ax-oset-kstep-unique,ax-oset-kstep-factor} as minimal
assumptions to prove basic properties of opetopic sets in
\cref{sec:prop-opet-sets}.

\section{Properties of opetopic sets}
\label{sec:prop-opet-sets}

We prove basic properties of opetopic sets and morphisms between them:
the underlying category of an opetopic set has a canonical
presentation (\cref{sec:generators-relations}); every slice of an
opetopic set is finite (\cref{sec:local-finiteness}); if two morphisms
of opetopic sets \(A \myMorphism A'\) agree at an object
\(x \myElemOf A\), then they agree on the slice \(A \mySlice x\)
(\cref{sec:local-uniq-morph}); any morphism of opetopic sets induces
an equivalence between slices (\cref{sec:local-equivalence}); every
slice \(\myOSet \mySlice A\) is a presheaf category
(\cref{sec:opetopic-sets-as}).

\subsection{Generators and relations}
\label{sec:generators-relations}

We show that the underlying category of an opetopic set is canonically
presented by generating \(1\)-step arrows and equations between
\(2\)-step arrows (\cref{prop-oset-canonical-presentation}).

\begin{myConstruction}
  \label{cst-oset-canonical-presentation}
  Let \(A\) be an opetopic set. The set \(\myApp{\myOSetGen}{A}\) of
  \myDefine{canonical generators for \(A\)} is the set of \(1\)-step
  arrows in \(A\). The set \(\myApp{\myOSetRel}{A}\) of
  \myDefine{canonical relations for \(A\)} is the set of equations of
  the form \(f_{1} \myComp g_{1} \myId f_{2} \myComp g_{2}\) that hold
  in \(A\) such that \((f_{1} \myComma g_{1})\) is heterogeneous and
  \((f_{2} \myComma g_{2})\) is homogeneous. Let
  \(\myApp{\myOSetAlt}{A}\) be the category with the same object as
  \(A\) and presented by generators \(\myApp{\myOSetGen}{A}\) and
  relations \(\myApp{\myOSetRel}{A}\). By definition, we have a
  canonical functor \(\myApp{\myOSetAlt}{A} \myMorphism A\), which is
  conservative because there is no way to construct a non-trivial
  equivalence in \(\myApp{\myOSetAlt}{A}\), by which we regard
  \(\myApp{\myOSetAlt}{A}\) an \(\myFinOrd\)-direct category.
\end{myConstruction}

We show that the canonical functor
\(\myApp{\myOSetAlt}{A} \myMorphism A\) is an equivalence
(\cref{prop-oset-canonical-presentation}). We prepare a lemma on a
tree structure in an opetopic set (\cref{prop-oset-tree}).

\begin{myDefinition}
  A \myDefine{graph} \(G\) consists of a set \(\myApp{\myVertex}{G}\)
  of \myDefine{vertices} and a set
  \(\myApp{\myEdge_{G}}{x \myComma y}\) of \myDefine{edges from \(x\)
    to \(y\)} for every
  \(x \myComma y \myElemOf \myApp{\myVertex}{G}\). We may write
  \(f \myElemOf x \myMorphism y\) instead of
  \(f \myElemOf \myApp{\myEdge_{G}}{x \myComma y}\). For vertices
  \(x \myComma y \myElemOf G\), a \myDefine{path from \(x\) to \(y\)}
  is a chain of edges
  \begin{equation*}
    x \myId z_{0} \myXMorphism{f_{0}} z_{1} \myXMorphism{f_{1}}
    \myDots \myXMorphism{f_{m - 1}} z_{m} \myId y.
  \end{equation*}
  We say a graph \(G\) is a \myDefine{tree} if there exists a vertex
  \(r \myElemOf G\) such that, for every vertex \(x \myElemOf G\),
  there exists a unique path from \(x\) to \(r\). In other words,
  \(r\) is the terminal object in the free category over \(G\), from
  which it follows that such a vertex \(r\) is unique. We refer to
  \(r\) as the \myDefine{root in \(G\)}.
\end{myDefinition}

When \(G\) is a tree, the following induction principle is valid. Let
\(P\) be a property on vertices in \(G\). Suppose:
\begin{itemize}
\item the root in \(G\) satisfies \(P\); and
\item for every edge \(f \myElemOf x \myMorphism y\) in \(G\), if
  \(y\) satisfies \(P\), then \(x\) satisfies \(P\).
\end{itemize}
Then every vertex \(x \myElemOf G\) satisfies \(P\). This induction
principle is justified by induction on the length of the unique path
from \(x\) to the root.

\begin{myConstruction}
  Let \(A\) be a preopetopic set and let \(x \myElemOf A\). We define
  the \myDefine{source slice}
  \(A \mySourceSlice x \mySub A \myOneSlice x\) to be the subset
  spanned by the source arrows into \(x\). We define the
  \myDefine{target slice}
  \(A \myTargetSlice x \mySub A \myOneSlice x\) to be the subset
  spanned by the target arrows into \(x\). By definition,
  \((A \mySourceSlice x) \myBinCoprod (A \myTargetSlice x) \myEquiv (A
  \myOneSlice x)\).
\end{myConstruction}

\begin{myConstruction}
  Let \(A\) be an opetopic set and let \(x \myElemOf A\) be an object
  of degree \(\myGe 1\). We refer to the unique target arrow into
  \(x\), which exists by \cref{ax-oset-target}, as
  \(\myApp{\myTargetOf_{A}}{x} \myElemOf \myApp{\myDomTargetOf_{A}}{x}
  \myTargetMorphism x\) or
  \(\myApp{\myTargetOf}{x} \myElemOf \myApp{\myDomTargetOf}{x}
  \myTargetMorphism x\) if \(A\) is clear from the context. For
  \(k \myLe \myApp{\myDegree}{x}\), we define
  \(\myApp{\myKTargetOf{k}}{x} \myElemOf \myApp{\myDomKTargetOf{k}}{x}
  \myKMorphism{k} x\) by
  \(\myApp{\myKTargetOf{0}}{x} \myDefEq \myIdMorphism_{x}\) and
  \(\myApp{\myKTargetOf{k + 1}}{x} \myDefEq \myApp{\myKTargetOf{k}}{x}
  \myComp \myApp{\myTargetOf}{\myApp{\myDomKTargetOf{k}}{x}}\).
\end{myConstruction}

\begin{myConstruction}
  Let \(A\) be an opetopic set and let \(x \myElemOf A\) be an object
  of degree \(\myGe 1\). We define a graph
  \(\myApp{\myOGraph}{A \myComma x}\) as follows. The set of vertices
  in \(\myApp{\myOGraph}{A \myComma x}\) is
  \((A \mySourceSlice x) \myBinCoprod (A \myTwoSlice x)\). There is no
  edge between vertices from \(A \mySourceSlice x\). There is no edge
  between vertices from \(A \myTwoSlice x\). An edge from
  \(f \myElemOf A \mySourceSlice x\) to
  \(g \myElemOf A \myTwoSlice x\) is a target arrow
  \(t \myElemOf g \myTargetMorphism f\) in \(A \mySlice x\). An edge
  from \(g \myElemOf A \myTwoSlice x\) to
  \(f \myElemOf A \mySourceSlice x\) is a source arrow
  \(s \myElemOf g \mySourceMorphism f\) in \(A \mySlice x\).
\end{myConstruction}

\begin{myLemma}
  \label{prop-oset-tree}
  Let \(A\) be an opetopic set and let \(x \myElemOf A\) be an object
  of degree \(\myGe 1\).
  \begin{enumerate}
  \item The graph \(\myApp{\myOGraph}{A \myComma x}\) is a tree.
  \item If \(\myApp{\myDegree}{x} \myId 1\), then
    \(\myApp{\myOGraph}{A \myComma x}\) is a singleton set with no
    edge.
  \item If \(\myApp{\myDegree}{x} \myGe 2\), then the root in
    \(\myApp{\myOGraph}{A \myComma x}\) is
    \(\myApp{\myKTargetOf{2}}{x} \myElemOf A \myKSlice{2} x\).
  \end{enumerate}
\end{myLemma}
\begin{proof}
  Suppose that \(\myApp{\myDegree}{x} \myId 1\). Then the set of
  vertices in \(\myApp{\myOGraph}{A \myComma x}\) is a singleton set
  since \(A \myTwoSlice x\) is empty and \(A \mySourceSlice x\) is a
  singleton by \cref{ax-oset-source}.
  \(\myApp{\myOGraph}{A \myComma x}\) has no edge by
  definition. \(\myApp{\myOGraph}{A \myComma x}\) is a tree whose root
  is the unique vertex. Suppose that \(\myApp{\myDegree}{x} \myGe
  2\). Take \(r\) as in \cref{ax-oset-path}. By the definition of
  \(\myApp{\myOGraph}{A \myComma x}\), \cref{eq-oset-path} is a path
  in \(\myApp{\myOGraph}{A \myComma x}\) from \(f\) to \(r\). Thus,
  there is a path \(p\) from \(\myApp{\myKTargetOf{2}}{x}\) to
  \(r\). By \cref{ax-oset-homo}, there is no edge from
  \(\myApp{\myKTargetOf{2}}{x}\). Hence, the length of \(p\) is \(0\),
  and thus \(\myApp{\myKTargetOf{2}}{x} \myId r\). For every
  \(f \myElemOf A \myTwoSlice x\), there is a path from \(f\) to
  \(\myApp{\myKTargetOf{2}}{x}\). Such a path is unique by
  \cref{ax-oset-homo,ax-oset-target} and by the fact that there is no
  edge from \(\myApp{\myKTargetOf{2}}{x}\). For every
  \(f \myElemOf A \mySourceSlice x\), we have the unique edge
  \(t \myElemOf f \myMorphism g\) by \cref{ax-oset-target} and then a
  unique path from \(g\) to \(\myApp{\myKTargetOf{2}}{x}\) as we have
  seen.
\end{proof}

We prove \(\myApp{\myOSetAlt}{A} \myEquiv A\) by a
\emph{normalization} procedure.

\begin{myDefinition}
  Let \(A\) be an opetopic set and let \(k \myGe 2\). We say a
  composable tuple of \(1\)-step arrows
  \((f_{1} \myComma \myDots \myComma f_{k})\) in \(A\) is \myDefine{in
    normal form} if \(f_{1} \myComma \myDots \myComma f_{k - 2}\) are
  target arrows and \((f_{k - 1} \myComma f_{k})\) is homogeneous. For
  \(x \myComma y \myElemOf A\) such that
  \(\myApp{\myDegree}{x} + k \myId \myApp{\myDegree}{y}\), let
  \(\myApp{\myKNF{k}_{A}}{x \myComma y}\) denote the set of composable
  tuples of \(1\)-step arrows
  \((f_{1} \myComma \myDots \myComma f_{k})\) in normal form such that
  \(f_{1} \myComp \myDots \myComp f_{k} \myElemOf x \myKMorphism{k}
  y\).
\end{myDefinition}

We first consider normalization of three \(1\)-step arrows.

\begin{myConstruction}
  Let \(A\) be an opetopic set and let
  \(f \myElemOf y \myThreeMorphism x\) be a \(3\)-step arrow. We
  define a graph \(\myApp{\myOGraphIII}{A \myComma f}\) as follows. A
  vertex in \(\myApp{\myOGraphIII}{A \myComma f}\) is a factorization
  \((p \myComma q \myComma r)\) of \(f\) into three \(1\)-step arrows
  \(f \myId p \myComp q \myComp r\). We add an edge from
  \((p_{1} \myComma q_{1} \myComma r_{1})\) to
  \((p_{2} \myComma q_{2} \myComma r_{2})\) when one of the following
  holds.
  \begin{myWithCustomLabel}{condition}
    \begin{enumerate}
    \item \label{item-sxx} \(p_{1} \myId p_{2}\) is a source arrow,
      \(q_{1} \myComp r_{1} \myId q_{2} \myComp r_{2}\),
      \((q_{1} \myComma r_{1})\) is heterogeneous, and
      \((q_{2} \myComma r_{2})\) is homogeneous.
    \item \label{item-xxs} \(r_{1} \myId r_{2}\) is a source arrow,
      \(p_{1} \myComp q_{1} \myId p_{2} \myComp q_{2}\),
      \((p_{1} \myComma q_{1})\) is homogeneous, and
      \((p_{2} \myComma q_{2})\) is heterogeneous.
    \item \label{item-xxt} \(r_{1} \myId r_{2}\) is a target arrow,
      \(p_{1} \myComp q_{1} \myId p_{2} \myComp q_{2}\),
      \((p_{1} \myComma q_{1})\) is heterogeneous, and
      \((p_{2} \myComma q_{2})\) is homogeneous.
    \end{enumerate}
  \end{myWithCustomLabel}
  These are all the edges in \(\myApp{\myOGraphIII}{A \myComma
    f}\). \Cref{tab-oset-graph3-st} lists all possible combinations of
  source and target arrows when there is an edge
  \((p_{1} \myComma q_{1} \myComma r_{1}) \myMorphism (p_{2} \myComma
  q_{2} \myComma r_{2})\) in \(\myApp{\myOGraphIII}{A \myComma f}\).

  \begin{table}
    \begin{center}
      \begin{tabular}{cll}
        Condition
        & \((p_{1} \myComma q_{1} \myComma r_{1})\)
        & \((p_{2} \myComma q_{2} \myComma r_{2})\)
        \\ \hline
        \labelcref{item-sxx}
        & \((\mySourceMark \myComma \mySourceMark \myComma
        \myTargetMark)\) or \((\mySourceMark \myComma \myTargetMark
        \myComma \mySourceMark)\)
        & \((\mySourceMark \myComma \mySourceMark \myComma
          \mySourceMark)\) or \((\mySourceMark \myComma \myTargetMark
          \myComma \myTargetMark)\)
        \\ \labelcref{item-xxs}
        & \((\mySourceMark \myComma \mySourceMark \myComma
        \mySourceMark)\) or \((\myTargetMark \myComma \myTargetMark
        \myComma \mySourceMark)\)
        & \((\mySourceMark \myComma \myTargetMark \myComma
          \mySourceMark)\) or \((\myTargetMark \myComma \mySourceMark
          \myComma \mySourceMark)\)
        \\ \labelcref{item-xxt}
        & \((\mySourceMark \myComma \myTargetMark \myComma
        \myTargetMark)\) or \((\myTargetMark \myComma \mySourceMark
        \myComma \myTargetMark)\)
        & \((\mySourceMark \myComma \mySourceMark \myComma
          \myTargetMark)\) or \((\myTargetMark \myComma \myTargetMark
          \myComma \myTargetMark)\)
      \end{tabular}
    \end{center}
    \caption{Possible combinations of source (\(\mySourceMark\)) and
      target (\(\myTargetMark\)) arrows when there is an edge
      \((p_{1} \myComma q_{1} \myComma r_{1}) \myMorphism (p_{2}
      \myComma q_{2} \myComma r_{2})\) in
      \(\myApp{\myOGraphIII}{A \myComma f}\).}
    \label{tab-oset-graph3-st}
  \end{table}
\end{myConstruction}

\begin{myLemma}
  \label{prop-oset-graph3-succ}
  Let \(A\) be an opetopic set, let \(f\) be a \(3\)-step arrow in
  \(A\), and let
  \(v \myDefEq (p \myComma q \myComma r) \myElemOf
  \myApp{\myOGraphIII}{A \myComma f}\) be a vertex. Then either there
  is no edge from \(v\) or there is a unique edge from
  \(v\). Moreover, the former holds if and only if
  \((p \myComma q \myComma r)\) is in normal form.
\end{myLemma}
\begin{proof}
  By \cref{ax-oset-hetero,ax-oset-homo}. See also
  \cref{tab-oset-graph3-st}.
\end{proof}

\begin{myLemma}
  \label{prop-oset-graph3-terminate}
  Let \(A\) be an opetopic set and let
  \(f \myElemOf y \myThreeMorphism x\) be a \(3\)-step arrow in
  \(A\). Then there exists a natural number \(n\) such that every path
  in \(\myApp{\myOGraphIII}{A \myComma f}\) is of length at most
  \(n\).
\end{myLemma}
\begin{proof}
  By definition, a path \(\pi\) in
  \(\myApp{\myOGraphIII}{A \myComma f}\) is of the form
  \begin{equation*}
    \myDots \myMorphism v_{i} \myMorphism v'_{i} \myMorphism v_{i + 1}
    \myMorphism \myDots,
  \end{equation*}
  where \cref{item-sxx} holds for the edge
  \(v_{i} \myMorphism v'_{i}\) and Condition \labelcref{item-xxs} or
  \labelcref{item-xxt} holds for the edge
  \(v'_{i} \myMorphism v_{i + 1}\). Let
  \(v_{i} \myDefEq (p_{i} \myComma q_{i} \myComma r_{i})\) and
  \(v'_{i} \myDefEq (p_{i} \myComma q'_{i} \myComma r_{i + 1})\), and
  let \(g_{i} \myDefEq p_{i} \myComp q'_{i}\) and
  \(h_{i} \myId q_{i} \myComp r_{i}\). Let \(z_{i}\) be the codomain
  of \(q_{i}\) (or the domain of \(p_{i}\)). Because a vertex
  \((p \myComma q \myComma r)\) with \(p\) a target arrow appears only
  at the start or the end of a path by \cref{tab-oset-graph3-st}, we
  may assume that all the \(p_{i}\)'s are source arrows. If
  \(r_{i + 1}\) is a source arrow, then \(q'_{i}\) is a source arrow,
  \(q_{i + 1}\) is a target arrow, and we have the edges
  \(p_{i} \myXMorphismAlt{q'_{i}} g_{i} \myXMorphismAlt{q_{i + 1}}
  p_{i + 1}\) in \(\myApp{\myOGraph}{A \myComma x}\). If \(r_{i + 1}\)
  is a target arrow, then \(q'_{i}\) is a target arrow, \(q_{i + 1}\)
  is a source arrow, and we have the edges
  \(p_{i} \myXMorphism{q'_{i}} g_{i} \myXMorphism{q_{i + 1}} p_{i +
    1}\) in \(\myApp{\myOGraph}{A \myComma x}\). In particular,
  \(p_{i}\)'s and \(g_{i}\)'s form a zigzag \(\zeta\) in
  \(\myApp{\myOGraph}{A \myComma x}\). Let
  \(p \myElemOf A \mySourceSlice x\) and let
  \(i_{1} \myLt i_{2} \myLt \myDots \) be the indexes \(i\) such that
  \(p_{i} \myId p\). Note that we can find these indexes because the
  identity type on \(A \mySourceSlice x\) is decidable by
  \cref{ax-oset-finite}.

  Suppose that \(r_{i_{j} + 1}\) is a source arrow and \(i_{j + 1}\)
  exists. Since \(\myApp{\myOGraph}{A \myComma x}\) is a tree
  (\cref{prop-oset-tree}), the zigzag \(\zeta\) between \(p_{i_{j}}\) and
  \(p_{i_{j + 1}}\) begins with
  \(p_{i_{j}} \myXMorphismAlt{q'_{i_{j}}} g_{i_{j}}
  \myXMorphismAlt{q_{i_{j} + 1}} p_{i_{j} + 1}\) and ends with
  \(p_{i_{j} + 1} \myXMorphism{q_{i_{j} + 1}} g_{i_{j}}
  \myXMorphism{q'_{i_{j}}} p_{i_{j}} \myId p_{i_{j + 1}}\). Hence,
  \(p_{i_{j + 1} - 1} \myId p_{i_{j} + 1}\),
  \(q'_{i_{j + 1} - 1} \myId q_{i_{j} + 1}\), and
  \(q_{i_{j + 1}} \myId q'_{i_{j}}\). Then \(q_{i_{j + 1}}\) is always
  a source arrow, and we have the edges
  \(q_{i_{j}} \myXMorphism{r_{i_{j}}} h_{i_{j}} \myXMorphism{r_{i_{j}
      + 1}} q'_{i_{j}} \myId q_{i_{j + 1}}\) in
  \(\myApp{\myOGraph}{A \myComma z_{i_{j}}}\) when \(q_{i_{j}}\) is
  also a source arrow.

  Suppose that \(r_{i_{j} + 1}\) is a target arrow and \(i_{j + 1}\)
  exists. Let \(p'\) be the last \(p_{i}\) between \(p_{i_{j}}\) and
  \(p_{i_{j + 1}}\) in the path \(\pi\) with minimum depth (i.e. the
  length of the path to the root), and let
  \(k_{1} \myComma k_{2} \myDots\) be the indices \(k\) between
  \(i_{j}\) and \(i_{j + 1}\) such that \(p_{k} \myId p'\). Then we
  have the edges
  \(p_{k_{l} - 1} \myXMorphism{q'_{k_{l} - 1}} g_{k_{l} - 1}
  \myXMorphism{q_{k_{l}}} p_{k_{l}} \myXMorphismAlt{q'_{k_{l}}}
  g_{k_{l}} \myXMorphismAlt{q_{k_{l} + 1}} p_{k_{l} + 1}\) in
  \(\myApp{\myOGraph}{A \myComma x}\). We see that \(q_{k_{l}}\) and
  \(q'_{k_{l}}\) are source arrows, and we have the edges
  \(q_{k_{l}} \myXMorphism{r_{k_{l}}} h_{k_{l}} \myXMorphism{r_{k_{l}
      + 1}} q'_{k_{l}}\) in
  \(\myApp{\myOGraph}{A \myComma z_{k_{l}}}\). From the observation in
  the previous paragraph, \(q'_{k_{l}} \myId q_{k_{l + 1}}\) is a
  source arrow. It then follows that \(p_{i_{j}}\) and
  \(p_{i_{j + 1}}\) lie in different branches of \(p'\) in
  \(\myApp{\myOGraph}{A \myComma x}\), which contradicts that
  \(p_{i_{j}} \myId p_{i_{j + 1}} \myId p\). Therefore, if
  \(r_{i_{j} + 1}\) is a target arrow, then \(p_{i_{j}}\) is the end
  of occurrence of \(p\) in \(\pi\).

  These observations give a bound
  \(2 \times \myApp{\myCardinality}{A \myOneSlice x'}\) for the number of
  occurrences of each \(p \myElemOf x' \mySourceMorphism x\) in a path
  in \(\myApp{\myOGraphIII}{A \myComma f}\), where
  \(\myApp{\myCardinality}{A \myOneSlice x'}\) is the cardinality of
  the finite set \(A \myOneSlice x'\) (\cref{ax-oset-finite}). Since
  \(A \mySourceSlice x\) is also finite by \cref{ax-oset-finite}, we
  obtain a bound
  \(\sum_{x' \myElemOf A \mySourceSlice x} 2 \times \myApp{\myCardinality}{A
    \myOneSlice x'}\) for the lengths of paths in
  \(\myApp{\myOGraphIII}{A \myComma f}\).
\end{proof}

\begin{myLemma}
  \label{prop-oset-graph3-forest}
  Let \(A\) be an opetopic set, let \(f\) be a \(3\)-step arrow in
  \(A\), and let \(v \myElemOf \myApp{\myOGraphIII}{A \myComma
    f}\). Then there exists a unique path in
  \(\myApp{\myOGraphIII}{A \myComma f}\) from \(v\) to a vertex in
  normal form.
\end{myLemma}
\begin{proof}
  By \cref{prop-oset-graph3-succ}, we can uniquely extend a path from
  \(v\) until it reaches a vertex in normal form. By
  \cref{prop-oset-graph3-terminate}, this procedure terminates at such
  a vertex.
\end{proof}

We now obtain the normalization procedure.

\begin{myLemma}
  \label{prop-oset-normal-form}
  Let \(A\) be an opetopic set and \(k \myGe 2\). Then every
  \(k\)-step arrow \(f\) in \(\myApp{\myOSetAlt}{A}\) factors into
  \(k\) \(1\)-step arrows \(g_{1} \myComp \myDots \myComp g_{k}\) in
  normal form.
\end{myLemma}
\begin{proof}
  We proceed by induction on \(k\). The case when \(k \myId 2\) is by
  \cref{ax-oset-homo}. Suppose that \(k \myGe 3\). By definition,
  \(f\) factors into \(1\)-step arrows
  \(f_{1} \myComp \myDots \myComp f_{k}\). By
  \cref{prop-oset-graph3-forest},
  \(f_{1} \myComp f_{2} \myComp f_{3} \myId g_{1} \myComp g_{2}
  \myComp g_{3}\) in \(\myApp{\myOSetAlt}{A}\) with
  \((g_{1} \myComma g_{2} \myComma g_{3})\) in normal form, since
  \(p_{1} \myComp q_{1} \myComp r_{1} \myId p_{2} \myComp q_{2}
  \myComp r_{2}\) in \(\myApp{\myOSetAlt}{A}\) whenever there is an
  edge
  \((p_{1} \myComma q_{1} \myComma r_{1}) \myMorphism (p_{2} \myComma
  q_{2} \myComma r_{2})\) in \(\myApp{\myOGraphIII}{A \myComma f}\) by
  definition. Then apply the induction hypothesis for
  \(g_{2} \myComp g_{3} \myComp f_{4} \myComp \myDots \myComp f_{k}\).
\end{proof}

\begin{myLemma}
  \label{prop-oset-normalization}
  Let \(A\) be an opetopic set, let \(k \myGe 2\), and let
  \(x \myComma y \myElemOf A\) be objects such that
  \(\myApp{\myDegree}{x} + k \myId \myApp{\myDegree}{y}\). Then
  \(\myApp{\myKNF{k}_{A}}{x \myComma y} \myEquiv
  \myApp{\myArr_{\myApp{\myOSetAlt}{A}}}{x \myComma y} \myEquiv
  \myApp{\myArr_{A}}{x \myComma y}\).
\end{myLemma}
\begin{proof}
  Let
  \(H \myElemOf \myApp{\myKNF{k}_{A}}{x \myComma y} \myMorphism
  \myApp{\myArr_{\myApp{\myOSetAlt}{A}}}{x \myComma y}\) and
  \(K \myElemOf \myApp{\myArr_{\myApp{\myOSetAlt}{A}}}{x \myComma y}
  \myMorphism \myApp{\myArr_{A}}{x \myComma y}\) denote the canonical
  maps. \(H\) is surjective by \cref{prop-oset-normal-form}. \(K\) is
  surjective by \cref{ax-oset-kstep-factor,ax-oset-homo}.
  \(K \myComp H\) is injective by
  \cref{ax-oset-target,ax-oset-kstep-unique,ax-oset-homo}. Thus,
  \(K \myComp H\) is an equivalence, and then \(H\) and \(K\) are also
  equivalences.
\end{proof}

\begin{myProposition}
  \label{prop-oset-canonical-presentation}
  Let \(A\) be an opetopic set. Then the canonical functor
  \(\myApp{\myOSetAlt}{A} \myMorphism A\) is an equivalence.
\end{myProposition}
\begin{proof}
  By construction, the canonical functor is an equivalence on objects
  and fully faithful on \(1\)-step arrows. It is fully faithful on
  \(k\)-step arrows for \(k \myGe 2\) by
  \cref{prop-oset-normalization}.
\end{proof}

\subsection{Local finiteness}
\label{sec:local-finiteness}

We show that every slice of an opetopic set is finite
(\cref{prop-oset-slice-finite}). We say a gaunt category \(C\) is
finite if \(\myApp{\myObj}{C}\) is finite and the set of arrows
\(\myApp{\myArr_{C}}{x \myComma y}\) is finite for all
\(x \myComma y \myElemOf C\).

\begin{myLemma}
  \label{prop-oset-kslice-normal-form}
  Let \(A\) be an opetopic set, let \(k \myGe 2\), and let
  \(x \myElemOf A\) be an object of degree \(\myGe k\). Then
  \(A \myKSlice{k} x \myEquiv ((y \myElemOf A \mySourceSlice
  \myApp{\myDomKTargetOf{k - 2}}{x}) \myBinProd (A \mySourceSlice y))
  \myBinCoprod ((y \myElemOf A \myTargetSlice \myApp{\myDomKTargetOf{k
      - 2}}{x}) \myBinProd (A \myTargetSlice y))\).
\end{myLemma}
\begin{proof}
  By \cref{prop-oset-normalization}.
\end{proof}

\begin{myLemma}
  \label{prop-oset-kslice-finite}
  Let \(A\) be an opetopic set, let \(x \myElemOf A\), and let
  \(k \myGe 0\). Then \(A \myKSlice{k} x\) is finite.
\end{myLemma}
\begin{proof}
  By case analysis on \(k\). The case when \(k \myId 0\) is because
  \(A \myKSlice{0} x \myEquiv \{\myIdMorphism_{x}\}\). The case when
  \(k \myId 1\) is \cref{ax-oset-finite}. The case when \(k \myId 2\)
  follows from \cref{prop-oset-kslice-normal-form,ax-oset-finite}.
\end{proof}

\begin{myLemma}
  \label{prop-category-arrow-slice-fiber}
  Let \(C\) be a category and let \(x \myComma y \myElemOf C\) be
  objects. Then \(\myApp{\myArr_{C}}{x \myComma y}\) is the fiber of
  \(\myApp{\myObj}{C \mySlice y} \myMorphism \myApp{\myObj}{C}\) over
  \(x \myElemOf \myApp{\myObj}{C}\).
\end{myLemma}
\begin{proof}
  By definition.
\end{proof}

\begin{myLemma}
  \label{prop-slice-category-arrow}
  Let \(C\) be a category, \(x \myElemOf C\), and
  \(y \myComma z \myElemOf C \mySlice x\). Then
  \(\myApp{\myArr_{C \mySlice x}}{y \myComma z}\) is the fiber of
  \(\myApp{\myObj}{C \mySlice z} \myMorphism \myApp{\myObj}{C \mySlice
    x}\) over \(y \myElemOf \myApp{\myObj}{C \mySlice x}\).
\end{myLemma}
\begin{proof}
  By \cref{prop-category-arrow-slice-fiber}.
\end{proof}

\begin{myLemma}
  \label{prop-gaunt-category-slice-finite}
  Let \(C\) be a gaunt category. Suppose that \(\myApp{\myObj}{C
    \mySlice x}\) is finite for every \(x \myElemOf C\). Then \(C
  \mySlice x\) is finite for every \(x \myElemOf C\).
\end{myLemma}
\begin{proof}
  By assumption, the object part of \(C \mySlice x\) is finite. The
  arrow part of it is also finite by \cref{prop-slice-category-arrow}.
\end{proof}

\begin{myProposition}
  \label{prop-oset-slice-finite}
  Let \(A\) be an opetopic set. Then \(A \mySlice x\) is finite for
  every \(x \myElemOf A\).
\end{myProposition}
\begin{proof}
  By \cref{prop-oset-slice-finite}, \(\myApp{\myObj}{A \mySlice x}\)
  is finite. Then apply \cref{prop-gaunt-category-slice-finite}.
\end{proof}

\subsection{Local uniqueness of morphisms}
\label{sec:local-uniq-morph}

The next goal is local uniqueness of morphisms of opetopic sets
(\cref{prop-mor-oset-slice-id}), which asserts that, if two morphisms
of opetopic sets \(A \myMorphism A'\) agree at an object
\(x \myElemOf A\), then they agree on the slice \(A \mySlice x\).

\begin{myConstruction}
  Let \(G\) be a graph and let \(x \myElemOf G\) be a vertex. We
  define
  \(G \mySlice x \myDefEq (y \myElemOf G) \myBinProd
  \myApp{\myEdge_{G}}{y \myComma x}\).
\end{myConstruction}

\begin{myLemma}
  \label{prop-oset-graph-slice-source}
  Let \(A\) be an opetopic set and let
  \(f \myElemOf y \mySourceMorphism x\) be a source arrow in
  \(A\). Then
  \(\myApp{\myOGraph}{A \myComma x} \mySlice f \myEquiv A
  \mySourceSlice y\).
\end{myLemma}
\begin{proof}
  By definition.
\end{proof}

\begin{myLemma}
  \label{prop-oset-graph-slice-2step}
  Let \(A\) be an opetopic set and let
  \(f \myElemOf y \myTwoMorphism x\) be a \(2\)-step arrow in
  \(A\). Then
  \(\myApp{\myOGraph}{A \myComma x} \mySlice f \myEquiv \myInitial\)
  or
  \(\myApp{\myOGraph}{A \myComma x} \mySlice f \myEquiv
  \myTerminal\). The former holds if and only if \(f\) factors as a
  source arrow followed by a target arrow. The latter holds if and
  only if \(f\) factors as a target arrow followed by a source arrow.
\end{myLemma}
\begin{proof}
  By definition and \cref{ax-oset-hetero}.
\end{proof}

\begin{myLemma}
  \label{prop-mor-oset-1slice-id}
  Let \(F_{1} \myComma F_{2} \myElemOf A \myMorphism A'\) be morphisms
  of opetopic sets, let \(x \myElemOf A\), and let \(x' \myElemOf A'\)
  such that \(\myApp{F_{1}}{x} \myId \myApp{F_{2}}{x} \myId x'\). Then
  the induced maps
  \(F_{1} \myOneSlice x \myComma F_{2} \myOneSlice x \myElemOf A
  \myOneSlice x \myMorphism A' \myOneSlice x'\) are identical.
\end{myLemma}
\begin{proof}
  We prove that the graph morphisms
  \(\myApp{\myOGraph}{F_{1} \myComma x} \myComma
  \myApp{\myOGraph}{F_{2} \myComma x} \myElemOf \myApp{\myOGraph}{A
    \myComma x} \myMorphism \myApp{\myOGraph}{A' \myComma x'}\) agree
  on vertices by induction on \(\myApp{\myDegree}{x}\). This implies
  that
  \(F_{1} \mySourceSlice x \myComma F_{2} \mySourceSlice x \myElemOf A
  \mySourceSlice x \myMorphism A' \mySourceSlice x'\) are
  identical. Then, since both \(F_{1}\) and \(F_{2}\) send
  \(\myApp{\myTargetOf}{x}\) to \(\myApp{\myTargetOf}{x'}\), we see
  that \(F_{1} \myOneSlice x \myId F_{2} \myOneSlice x\). The case
  when \(\myApp{\myDegree}{x} \myId 0\) is trivial since
  \(\myApp{\myOGraph}{A \myComma x}\) is empty. Suppose that
  \(\myApp{\myDegree}{x} \myGe 1\). We show that
  \(\myApp{F_{1}}{f} \myId \myApp{F_{2}}{f}\) by induction on
  \(f \myElemOf \myApp{\myOGraph}{A \myComma x}\). Both \(F_{1}\) and
  \(F_{2}\) send the root in \(\myApp{\myOGraph}{A \myComma x}\) to
  the root \(\myApp{\myOGraph}{A' \myComma x'}\) by
  \cref{prop-oset-tree}. Let \(h \myElemOf g \myMorphism f\) be an
  edge in \(\myApp{\myOGraph}{A \myComma x}\) and suppose that
  \(\myApp{F_{1}}{f} \myId \myApp{F_{2}}{f}\). If
  \(f \myElemOf A \mySourceSlice x\), then
  \(\myApp{F_{1}}{g} \myId \myApp{F_{2}}{g}\) by
  \cref{prop-oset-graph-slice-source} and by the induction hypothesis
  for the domain of \(f\). If \(f \myElemOf A \myTwoSlice x\), then
  \(\myApp{F_{1}}{g} \myId \myApp{F_{2}}{g}\) by
  \cref{prop-oset-graph-slice-2step}.
\end{proof}

\begin{myLemma}
  \label{prop-mor-oset-kslice-id}
  Let \(F_{1} \myComma F_{2} \myElemOf A \myMorphism A'\) be morphisms
  of opetopic sets, let \(k \myGe 0\), let \(x \myElemOf A\), and
  let \(x' \myElemOf A'\) such that
  \(\myApp{F_{1}}{x} \myId \myApp{F_{2}}{x} \myId x'\). Then the
  induced maps
  \(F_{1} \myOneSlice x \myComma F_{2} \myKSlice{k} x \myElemOf A
  \myKSlice{k} x \myMorphism A' \myKSlice{k} x'\) are identical.
\end{myLemma}
\begin{proof}
  By case analysis on \(k\). The case when \(k \myId 0\) is trivial
  since \(A \myKSlice{0} x\) is the singleton
  \(\{\myIdMorphism_{x}\}\). The case when \(k \myId 1\) is
  \cref{prop-mor-oset-1slice-id}. The case when \(k \myGe 2\) follows
  from \cref{prop-oset-kslice-normal-form,prop-mor-oset-1slice-id}.
\end{proof}

\begin{myLemma}
  \label{prop-slice-id-on-obj-to-id}
  Let \(F_{1} \myComma F_{2} \myElemOf C \myMorphism C'\) be functors
  between gaunt categories. Suppose that, for every \(x \myElemOf C\)
  and \(x' \myElemOf C\) such that
  \(\myApp{F_{1}}{x} \myId \myApp{F_{2}}{x} \myId x'\), the induced
  maps
  \(\myApp{\myObj}{F_{1} \mySlice x} \myComma \myApp{\myObj}{F_{2}
    \mySlice x} \myElemOf \myApp{\myObj}{C \mySlice x} \myMorphism
  \myApp{\myObj}{C' \mySlice x'}\) are identical. Then, for every
  \(x \myElemOf C\) and \(x' \myElemOf C'\) such that
  \(\myApp{F_{1}}{x} \myId \myApp{F_{2}}{x} \myId x'\), the induced
  functors
  \(F_{1} \mySlice x \myComma F_{2} \mySlice x \myElemOf C \mySlice x
  \myMorphism C' \mySlice x'\) are identical.
\end{myLemma}
\begin{proof}
  The object parts of \(F_{1} \mySlice x\) and \(F_{2} \mySlice x\)
  are identical by assumption. The arrow parts of them are also
  identical by \cref{prop-slice-category-arrow}.
\end{proof}

\begin{myProposition}
  \label{prop-mor-oset-slice-id}
  Let \(F_{1} \myComma F_{2} \myElemOf A \myMorphism A'\) be morphisms
  of opetopic sets, let \(x \myElemOf A\), and let
  \(x' \myElemOf A'\) such that
  \(\myApp{F_{1}}{x} \myId \myApp{F_{2}}{x} \myId x'\). Then the
  induced functors
  \(F_{1} \mySlice x \myComma F_{2} \mySlice x \myElemOf A \mySlice x
  \myMorphism A' \mySlice x'\) are identical.
\end{myProposition}
\begin{proof}
  By \cref{prop-mor-oset-kslice-id}, \(F_{1} \mySlice x\) and
  \(F_{2} \mySlice x\) agree on objects. Then apply
  \cref{prop-slice-id-on-obj-to-id}.
\end{proof}

\subsection{Local equivalence}
\label{sec:local-equivalence}

We show that any morphism of opetopic sets induces an equivalence
between slices (\cref{prop-mor-oset-slice-equiv}).

\begin{myLemma}
  \label{prop-mor-oset-1slice-equiv}
  Let \(F \myElemOf A \myMorphism A'\) be a morphism of opetopic
  sets and let \(x \myElemOf A\). Then the induced map
  \(F \myOneSlice x \myElemOf A \myOneSlice x \myMorphism A'
  \myOneSlice \myApp{F}{x}\) is an equivalence.
\end{myLemma}
\begin{proof}
  Let \(x' \myDefEq \myApp{F}{x}\). We prove that
  \(\myApp{\myOGraph}{F \myComma x} \myElemOf \myApp{\myOGraph}{A
    \myComma x} \myMorphism \myApp{\myOGraph}{A' \myComma x'}\) is an
  equivalence on vertices by induction on
  \(\myApp{\myDegree}{x}\). This implies that
  \(F \mySourceSlice x \myElemOf A \mySourceSlice x \myMorphism A'
  \mySourceSlice x'\) is an equivalence. Then, since the fiber of
  \(F\) over \(\myApp{\myTargetOf}{x'}\) is the singleton
  \(\{\myApp{\myTargetOf}{x}\}\), we see that \(F \myOneSlice x\) is
  an equivalence. The case when \(\myApp{\myDegree}{x} \myId 0\) is
  trivial since \(\myApp{\myOGraph}{A \myComma x}\) is empty. Suppose
  that \(\myApp{\myDegree}{x} \myGe 1\). We prove that the fiber of
  \(F\) over \(f'\) is contractible by induction on
  \(f' \myElemOf \myApp{\myOGraph}{A' \myComma x'}\). The fiber of
  \(F\) over the root in \(\myApp{\myOGraph}{A' \myComma x}\) consists
  of only the root in \(\myApp{\myOGraph}{A' \myComma x}\) by
  \cref{prop-oset-tree}. Let \(h' \myElemOf g' \myMorphism f'\) be an
  edge in \(\myApp{\myOGraph}{A' \myComma x'}\) and suppose that the
  fiber of \(F\) over \(f'\) is contractible with center \(f\). If
  \(f \myElemOf A \mySourceSlice x\), then the fiber of \(F\) over
  \(g'\) is contractible by \cref{prop-oset-graph-slice-source} and by
  the induction hypothesis for the domain of \(f\). If
  \(f \myElemOf A \myTwoSlice x\), then the fiber of \(F\) over \(g'\)
  is contractible by \cref{prop-oset-graph-slice-2step}.
\end{proof}

\begin{myLemma}
  \label{prop-mor-oset-kslice-equiv}
  Let \(F \myElemOf A \myMorphism A'\) be a morphism of opetopic sets,
  let \(x \myElemOf A\), and let \(k \myGe 0\). Then the induced map
  \(F \myKSlice{k} x \myElemOf A \myKSlice{k} x \myMorphism A'
  \myKSlice{k} \myApp{F}{x}\) is an equivalence.
\end{myLemma}
\begin{proof}
  By case analysis on \(k\). The case when \(k \myId 0\) is trivial
  since \(A \myKSlice{0} x\) is the singleton
  \(\{\myIdMorphism_{x}\}\). The case when \(k \myId 1\) is
  \cref{prop-mor-oset-1slice-id}. The case when \(k \myGe 2\) follows
  from \cref{prop-oset-kslice-normal-form,prop-mor-oset-1slice-equiv}.
\end{proof}

\begin{myLemma}
  \label{prop-slice-equiv-on-obj-to-equiv}
  Let \(F \myElemOf C \myMorphism C'\) be a functor. Suppose that the
  map
  \(\myApp{\myObj}{F \mySlice x} \myElemOf \myApp{\myObj}{C \mySlice
    x} \myMorphism \myApp{\myObj}{C' \mySlice \myApp{F}{x}}\) is an
  equivalence for every \(x \myElemOf C\). Then the functor
  \(F \mySlice x \myElemOf C \mySlice x \myMorphism C' \mySlice
  \myApp{F}{x}\) is an equivalence for every \(x \myElemOf C\).
\end{myLemma}
\begin{proof}
  The object part of \(F \mySlice x\) is an equivalence by
  assumption. The arrow part of it is also an equivalence by
  \cref{prop-slice-category-arrow}.
\end{proof}

\begin{myProposition}
  \label{prop-mor-oset-slice-equiv}
  Let \(F \myElemOf A \myMorphism A'\) be a morphism of opetopic sets
  and let \(x \myElemOf A\). Then the induced functor
  \(F \mySlice x \myElemOf A \mySlice x \myMorphism A' \mySlice
  \myApp{F}{x}\) is an equivalence.
\end{myProposition}
\begin{proof}
  By \cref{prop-mor-oset-kslice-equiv},
  \(\myApp{\myObj}{F \mySlice x}\) is an equivalence. Then apply
  \cref{prop-slice-equiv-on-obj-to-equiv}.
\end{proof}

\subsection{Opetopic sets over \(A\) as presheaves}
\label{sec:opetopic-sets-as}

We show that each slice \(\myOSet \mySlice A\) is a presheaf category
(\cref{prop-oset-slice-presheaf}).

\begin{myConstruction}
  Let \(A\) be a preopetopic set and let
  \(p \myElemOf B \myMorphism A\) be a morphism of
  \(\myFinOrd\)-direct categories. We extend \(B\) to a preopetopic
  set by
  \(\myApp{\mySource}{B} \myDefEq
  \myApp{p^{-1}}{\myApp{\mySource}{A}}\) and
  \(\myApp{\myTarget}{B} \myDefEq
  \myApp{p^{-1}}{\myApp{\myTarget}{A}}\). As a special case, every
  slice \(A \mySlice x\) is extended to a preopetopic set. This
  preopetopic set structure on \(B\) is characterized as the unique
  one that makes \(p\) a morphism of preopetopic sets. In other words,
  the forgetful functor
  \(\myPreOSet \mySlice A \myMorphism \myDirCat{\myFinOrd} \mySlice
  A\) is an equivalence. When \(A\) is an opetopic set, we regard
  \(\myOSet \mySlice A\) as a full subcategory of
  \(\myDirCat{\myFinOrd} \mySlice A\) via the equivalence
  \(\myPreOSet \mySlice A \myEquiv \myDirCat{\myFinOrd} \mySlice A\).
\end{myConstruction}

\begin{myDefinition}
  We say a property \(P\) on preopetopic sets is \myDefine{local} if a
  preopetopic set \(A\) satisfies \(P\) if and only if every slice
  \(A \mySlice x\) satisfies \(P\).
\end{myDefinition}

\begin{myLemma}
  \label{prop-oset-axioms-local}
  \Crefrange{ax-oset-first}{ax-oset-last} are local properties on
  preopetopic sets.
\end{myLemma}
\begin{proof}
  Straightforward.
\end{proof}

\begin{myDefinition}
  Let \(p \myElemOf B \myMorphism A\) be a functor. We say \(p\) is a
  \myDefine{right fibration} if the induced functor
  \(p \mySlice y \myElemOf B \mySlice y \myMorphism A \mySlice
  \myApp{p}{y}\) is an equivalence for every \(y \myElemOf B\). We say
  \(p\) is a \myDefine{discrete fibration} if it is a right fibration
  and every fiber \(B_{x}\) is a set. For a small category
  \(A\), the small discrete fibrations over \(A\) and the functors over
  \(A\) between them form a category \(\myApp{\myDFib}{A}\). Every
  discrete fibration is conservative. Thus, when \(A\) is an
  \(\myFinOrd\)-direct category, \(\myApp{\myDFib}{A}\) is regarded as
  a full subcategory of \(\myDirCat{\myFinOrd} \mySlice A\).
\end{myDefinition}

It is a standard fact that \(\myApp{\myDFib}{A}\) is equivalent to the
category of set-valued presheaves on \(A\).

\begin{myProposition}
  \label{prop-oset-slice-presheaf}
  Let \(A\) be an opetopic set. Then
  \(\myOSet \mySlice A \myId \myApp{\myDFib}{A}\) in the poset of
  full subcategories of \(\myDirCat{\myFinOrd} \mySlice A\).
\end{myProposition}
\begin{proof}
  Let \(P \myElemOf B \myMorphism A\) be a morphism of
  \(\myFinOrd\)-direct categories. If \(P\) is a right fibration, then
  \(B\) is an opetopic set by \cref{prop-oset-axioms-local} since
  \(B \mySlice y \myEquiv A \mySlice \myApp{P}{y}\) for every
  \(y \myElemOf B\). If \(B\) is an opetopic set, then \(P\) is a
  right fibration by \cref{prop-mor-oset-slice-equiv}.
\end{proof}

\section{Opetopes}
\label{sec:opetopes}

Opetopes are defined as special opetopic sets.

\begin{myDefinition}
  An \myDefine{opetope} is an opetopic set in which a terminal object
  exists. We refer to the terminal object in an opetope \(A\) as
  \(\myOpeTerminal_{A}\). Let \(\myOpetope \mySub \myOSet\) denote the
  full subcategory spanned by the opetopes.
\end{myDefinition}

An interesting phenomenon is that \(\myOpetope\) is extended to a
small opetopic set (\cref{prop-opetope-is-oset,prop-opetope-small}),
and thus \(\myOpetope\) is regarded as an object in \(\myOSet\).

\begin{myConstruction}
  We extend \(\myOpetope\) to a preopetopic set as follows. The degree
  functor is \(A \myMapsTo
  \myApp{\myDegree}{\myOpeTerminal_{A}}\). Clearly it is a functor to
  \(\myFinOrd\). It reflects equivalences by
  \cref{prop-mor-oset-slice-equiv}. It is then gaunt by
  \cref{prop-mor-oset-slice-id}. Thus, \(\myOpetope\) is an
  \(\myFinOrd\)-direct category. We say a \(1\)-step morphism
  \(F \myElemOf A \myOneMorphism A'\) of opetopes is a source/target
  arrow if the arrow
  \(\myApp{F}{\myOpeTerminal_{A}} \myOneMorphism \myOpeTerminal_{A'}\)
  is a source/target arrow in \(A'\).
\end{myConstruction}

\begin{myLemma}
  \label{prop-oset-canonical-equiv}
  Let \(A\) be an opetopic set. The morphism of preopetopic sets
  \(A \myMorphism \myOpetope \mySlice A\) that sends \(x \myElemOf A\)
  to the forgetful functor
  \(x_{\myBang} \myElemOf A \mySlice x \myMorphism A\) is an
  equivalence.
\end{myLemma}
\begin{proof}
  The inverse is given by
  \((F \myElemOf A' \myMorphism A) \myMapsTo
  \myApp{F}{\myOpeTerminal_{A'}}\). For \(x \myElemOf A\), we have
  \(\myApp{x_{\myBang}}{\myIdMorphism_{x}} \myId x\). For
  \(F \myElemOf A' \myMorphism A\), we have
  \(A' \myEquiv A \mySlice \myApp{F}{\myOpeTerminal_{A'}}\) by
  \cref{prop-mor-oset-slice-equiv}.
\end{proof}

\begin{myProposition}
  \label{prop-opetope-is-oset}
  \(\myOpetope\) is an opetopic set.
\end{myProposition}
\begin{proof}
  By \cref{prop-oset-canonical-equiv,prop-oset-axioms-local}.
\end{proof}

\begin{myProposition}
  \label{prop-opetope-small}
  \(\myOpetope\) is small
\end{myProposition}
\begin{proof}
  This is because every opetope is finite by
  \cref{prop-oset-slice-finite}.
\end{proof}

We show that \(\myOpetope \myElemOf \myOSet\) is moreover the terminal
object (\cref{prop-opetope-terminal}). As a corollary, \(\myOSet\) is
a presheaf category (\cref{prop-oset-presheaf}).

\begin{myLemma}
  \label{prop-natural-cone-to-terminal}
  Let \(C\) be a category and let \(x \myElemOf C\) be an
  object. Suppose that we have a natural transformation
  \(t \myElemOf \myIdMorphism_{C} \myMMorphism x\) such that
  \(t_{x} \myId \myIdMorphism_{x}\). Then \(x\) is a terminal object.
\end{myLemma}
\begin{proof}
  Let \(x' \myElemOf C\) be an object. We show that
  \(\myApp{\myArr_{C}}{x' \myComma x}\) is contractible. We have the
  arrow \(t_{x'} \myElemOf x' \myMorphism x\). Let
  \(f \myElemOf x' \myMorphism x\) be an arrow. By the naturality of
  \(t\), we have \(t_{x'} \myId t_{x} \myComp f\). Since
  \(t_{x} \myId \myIdMorphism_{x}\), we have \(t_{x'} \myId f\).
\end{proof}

\begin{myProposition}
  \label{prop-opetope-terminal}
  \(\myOpetope \myElemOf \myOSet\) is the terminal object.
\end{myProposition}
\begin{proof}
  For \(A \myElemOf \myOSet\), we define a morphism of opetopic sets
  \(t_{A} \myElemOf A \myMorphism \myOpetope\) by
  \(\myApp{t_{A}}{x} \myDefEq A \mySlice x\). This is natural in \(A\)
  by \cref{prop-mor-oset-slice-equiv}. The component at
  \(\myOpetope \myElemOf \myOSet\) is the morphism
  \((X \myMapsTo \myOpetope \mySlice X) \myElemOf \myOpetope
  \myMorphism \myOpetope\), which is equivalent to the identity on
  \(\myOpetope\) by \cref{prop-oset-canonical-equiv}. Then apply
  \cref{prop-natural-cone-to-terminal}.
\end{proof}

\begin{myTheorem}
  \label{prop-oset-presheaf}
  \(\myOSet \myEquiv \myApp{\myDFib}{\myOpetope}\).
\end{myTheorem}
\begin{proof}
  By \cref{prop-oset-slice-presheaf,prop-opetope-terminal}.
\end{proof}

Let us determine the opetopes of low degrees.

\begin{myNotation}
  Let \(A\) be an \(\myFinOrd\)-direct category and let
  \(n \myElemOf \myFinOrd\). The fiber of
  \(\myDegree \myElemOf A \myMorphism \myFinOrd\) over \(n\) is
  denoted by \(A_{n}\).
\end{myNotation}

\begin{myProposition}
  \label{prop-opetope-degree0}
  \(\myOpetope_{0} \myEquiv \myTerminal\).
\end{myProposition}
\begin{proof}
  The singleton \(\{0\}\) with \(\myApp{\myDegree}{0} \myDefEq 0\) is
  the only opetope of degree \(0\).
\end{proof}

\begin{myProposition}
  \label{prop-opetope-degree1}
  \(\myOpetope_{1} \myEquiv \myTerminal\). Moreover, for the unique
  opetope \(A\) of degree \(1\), there exist a unique source
  morphism into \(A\) and a unique target morphism into \(A\).
\end{myProposition}
\begin{proof}
  Let \(A\) be an opetope of degree \(1\). By \cref{ax-oset-source},
  there exists a unique source arrow
  \(x \mySourceMorphism \myOpeTerminal_{A}\). Since
  \(\myOpeTerminal_{A}\) is the terminal object and since the sets of
  source and target arrows are disjoint, we see that
  \(x \myNotId \myApp{\myTargetOf}{\myOpeTerminal_{A}}\). Hence, \(A\)
  must look like
  \begin{equation*}
    x \mySourceMorphism \myOpeTerminal_{A} \myTargetMorphismAlt
    \myApp{\myDomTargetOf}{\myOpeTerminal_{A}}.
  \end{equation*}
  It is straightforward to check that this preopetopic set is indeed
  an opetopic set.
\end{proof}

\section{The category of opetopic sets}
\label{sec:categ-opet-sets}

We study the category \(\myOSet\) of opetopic sets in more detail. We
first give a way to detect equivalences in \(\myOSet\). Of course,
equivalences in \(\myOSet \myEquiv \myApp{\myDFib}{\myOpetope}\)
(\cref{prop-oset-presheaf}) are detected fiberwise, but since we do
not know yet much about opetopes, this is not so helpful. A more
useful sufficient condition is degreewise equivalence
(\cref{prop-oset-fiber-conservative}).

\begin{myLemma}
  \label{prop-oset-dircat-conservative}
  The forgetful functor \(\myOSet \myMorphism \myGauntCat\) is
  conservative.
\end{myLemma}
\begin{proof}
  By definition, a morphism \(F \myElemOf A \myMorphism B\) of
  opetopic sets is an equivalence if and only if its underlying
  functor is an equivalence and it reflects source and target
  arrows. The second condition automatically holds since the sets of
  source and target arrows are complement to each other.
\end{proof}

\begin{myLemma}
  \label{prop-obj-equiv-slice-equiv-to-equiv}
  Let \(F \myElemOf C \myMorphism D\) be a functor and suppose that
  \(\myApp{\myObj}{F} \myElemOf \myApp{\myObj}{C} \myMorphism
  \myApp{\myObj}{D}\) is an equivalence and that
  \(F \mySlice x \myElemOf C \mySlice x \myMorphism D \mySlice
  \myApp{F}{x}\) is an equivalence for every \(x \myElemOf C\). Then
  \(F\) is an equivalence.
\end{myLemma}
\begin{proof}
  By assumption, \(F\) is an equivalence on objects. It is also fully
  faithful by \cref{prop-category-arrow-slice-fiber}.
\end{proof}

\begin{myProposition}
  \label{prop-oset-fiber-conservative}
  The functors
  \((A \myMapsTo A_{n}) \myElemOf \myOSet \myMorphism \mySet\) for all
  \(n \myElemOf \myFinOrd\) are jointly conservative.
\end{myProposition}
\begin{proof}
  Let \(F \myElemOf A \myMorphism B\) be a morphism of opetopic sets
  and suppose that \(F_{n} \myElemOf A_{n} \myMorphism B_{n}\) is an
  equivalence for every \(n \myElemOf \myFinOrd\). To see that \(F\)
  is an equivalence, by \cref{prop-oset-dircat-conservative}, it
  suffices to see that the underlying functor of \(F\) is an
  equivalence, but this follows from
  \cref{prop-mor-oset-slice-equiv,prop-obj-equiv-slice-equiv-to-equiv}
  since \(F\) is an equivalence on objects by assumption.
\end{proof}

We give some tools to compute colimits in \(\myOSet\). Fiberwise
computation of colimits in
\(\myOSet \myEquiv \myApp{\myDFib}{\myOpetope}\) is not helpful, and
degreewise computation (\cref{prop-oset-colimit-fiberwise}) is what we
want.

\begin{myConstruction}
  Let \(A\) be an opetopic set and let \(n \myElemOf \myFinOrd\). We
  define an opetopic set \(A_{\myLt n}\) to be the category of
  elements for the proposition-valued presheaf
  \(x \myMapsTo (\myApp{\myDegree}{x} \myLt n)\) on \(A\).
\end{myConstruction}

\begin{myProposition}
  \label{prop-oset-restrict-colimit}
  Let \(n \myElemOf \myFinOrd\). Then the functor
  \((A \myMapsTo A_{\myLt n}) \myElemOf \myOSet \myMorphism \myOSet\)
  preserves small colimits and pullbacks.
\end{myProposition}
\begin{proof}
  By construction and by \cref{prop-opetope-terminal}, the functor
  \(A \myMapsTo A_{\myLt n}\) factors as the pullback functor
  \(\myOSet \myEquiv \myOSet \mySlice \myOpetope \myMorphism \myOSet
  \mySlice \myOpetope_{\myLt n}\) followed by the forgetful functor
  \(\myOSet \mySlice \myOpetope_{\myLt n} \myMorphism \myOSet\). These
  functors preserve small colimits and pullbacks.
\end{proof}

\begin{myProposition}
  \label{prop-oset-fiber-colimit}
  Let \(n \myElemOf \myFinOrd\). Then the functor
  \((A \myMapsTo A_{n}) \myElemOf \myOSet \myMorphism \mySet\)
  preserves small colimits and pullbacks.
\end{myProposition}
\begin{proof}
  The functor \(A \myMapsTo A_{n}\) factors as the pullback functor
  \(\myOSet \myEquiv \myApp{\myDFib}{\myOpetope} \myMorphism \mySet
  \mySlice \myOpetope_{n}\) followed by the forgetful functor
  \(\mySet \mySlice \myOpetope_{n} \myMorphism \mySet\). These
  functors preserve small colimits and pullbacks.
\end{proof}

\begin{myProposition}
  \label{prop-oset-colimit-fiberwise}
  Let \(A \myElemOf I \myMorphism \myOSet\) be a diagram. Then a
  cocone
  \((\myApp{f}{i} \myElemOf \myApp{A}{i} \myMorphism B)_{i \myElemOf
    I}\) under \(A\) is a colimit cocone if and only if
  \((\myApp{f}{i}_{n} \myElemOf \myApp{A}{i}_{n} \myMorphism B_{n})_{i
    \myElemOf I}\) is a colimit cocone for every
  \(n \myElemOf \myFinOrd\).
\end{myProposition}
\begin{proof}
  By \cref{prop-oset-fiber-conservative,prop-oset-fiber-colimit}.
\end{proof}

\begin{myProposition}
  \label{prop-n-oset-colimit-fiberwise}
  Let \(n \myGe 0\) and let
  \(A \myElemOf I \myMorphism \myOSet \mySlice \myOpetope_{\myLt n +
    1}\). Then a cocone
  \((\myApp{f}{i} \myElemOf \myApp{A}{i} \myMorphism B)_{i \myElemOf
    I}\) under \(A\) is a colimit cocone if and only if
  \((\myApp{f}{i}_{n} \myElemOf \myApp{A}{i}_{n} \myMorphism B_{n})_{i
    \myElemOf I}\) and
  \((\myApp{f}{i}_{\myLt n} \myElemOf \myApp{A}{i}_{\myLt n}
  \myMorphism B_{\myLt n})_{i \myElemOf I}\) are colimit cocones.
\end{myProposition}
\begin{proof}
  Note that \(X_{m} \myEquiv \myInitial\) for any
  \(X \myElemOf \myOSet \mySlice \myOpetope_{\myLt n + 1}\) and
  \(m \myGe n + 1\). Thus, the claim follows from
  \cref{prop-oset-restrict-colimit,prop-oset-colimit-fiberwise}.
\end{proof}

\section{Boundaries and pasting diagrams}
\label{sec:bound-past-diagr}

We introduce boundaries and pasting diagrams. We show that an opetope
is completely determined by its pasting diagram of source objects
(\cref{prop-opetope-equiv-pd}).

\begin{myDefinition}
  Let \(n \myElemOf \myFinOrd\). An \myDefine{\(n\)-opetopic set} is
  an opetopic set \(A\) whose degree functor factors through
  \(n \mySub \myFinOrd\). This is equivalent to that the morphism
  \(A \myMorphism \myOpetope\) factors through
  \(\myOpetope_{\myLt n}\).
\end{myDefinition}

\begin{myDefinition}
  Let \(n \myElemOf \myFinOrd\). An \myDefine{\(n\)-preboundary} is an
  \(n\)-opetopic set \(A\) equipped with a subset
  \(\myApp{\myBdSource}{A} \mySub A_{n - 1}\) with complement
  \(\myApp{\myBdTarget}{A}\). Objects in \(\myApp{\myBdSource}{A}\)
  are called \myDefine{source objects}. Objects in
  \(\myApp{\myBdTarget}{A}\) are called \myDefine{target objects}.
\end{myDefinition}

\begin{myDefinition}
  Let \(n \myElemOf \myFinOrd\). An \myDefine{\(n\)-prepasting
    diagram} is an \((n + 1)\)-opetopic set \(A\) equipped with two
  families of sets
  \(\myPDLeaf_{A} \myComma \myPDRoot_{A} \myElemOf A_{n - 1}
  \myMorphism \mySet\). We say an object \(x \myElemOf A_{n - 1}\) is
  a \myDefine{leaf object} if \(\myApp{\myPDLeaf_{A}}{x}\) is
  inhabited. We say an object \(x \myElemOf A_{n - 1}\) is a
  \myDefine{root object} if \(\myApp{\myPDRoot_{A}}{x}\) is
  inhabited. Let
  \(\myApp{\myPDLeaf}{A} \myDefEq (x \myElemOf A_{n - 1}) \myBinProd
  \myApp{\myPDLeaf_{A}}{x}\) and
  \(\myApp{\myPDRoot}{A} \myDefEq (x \myElemOf A_{n - 1}) \myBinProd
  \myApp{\myPDRoot_{A}}{x}\). When extending an \((n + 1)\)-opetopic
  set \(A\) to an \(n\)-prepasting diagram, we specify either
  \(\myPDLeaf_{A}\) and \(\myPDRoot_{A}\) or
  \(\myApp{\myPDLeaf}{A} \myMorphism A_{n - 1}\) and
  \(\myApp{\myPDRoot}{A} \myMorphism A_{n - 1}\).
\end{myDefinition}

\begin{myConstruction}
  Let \(n \myGe 1\) and let \(A\) be an \(n\)-preboundary. We define
  an \(n\)-opetopic set \(\myApp{\mySourceHorn}{A}\) called the
  \myDefine{source horn of \(A\)} to be the category of elements for
  the following proposition-valued presheaf on \(A\).
  \begin{equation*}
    x \myMapsTo \left\{
      \begin{array}{ll}
        \myTextProp{\(x\) is a source object}
        & \text{if \(\myApp{\myDegree}{x} \myId n - 1\)} \\
        \myTerminal
        & \text{if \(\myApp{\myDegree}{x} \myLt n - 1\)}
      \end{array}
    \right.
  \end{equation*}
  We extend \(\myApp{\mySourceHorn}{A}\) to an \((n - 1)\)-prepasting
  diagram by
  \(\myApp{\myPDLeaf_{\myApp{\mySourceHorn}{A}}}{x} \myDefEq (y
  \myElemOf \myApp{\myBdTarget}{A}) \myBinProd (x \mySourceMorphism
  y)\) and
  \(\myApp{\myPDRoot_{\myApp{\mySourceHorn}{A}}}{x} \myDefEq (y
  \myElemOf \myApp{\myBdTarget}{A}) \myBinProd (x \myTargetMorphism
  y)\).
\end{myConstruction}

\begin{myConstruction}
  Let \(n \myGe 0\) and let \(A\) be an \(n\)-prepasting diagram. We
  define an \(n\)-opetopic set \(\myApp{\myBoundary}{A}\) called the
  \myDefine{boundary of \(A\)} to be the category of elements for the
  following set-valued presheaf on \(A\).
  \begin{equation*}
    x \myMapsTo \left\{
      \begin{array}{ll}
        \myInitial
        & \text{if \(\myApp{\myDegree}{x} \myId n\)} \\
        \myApp{\myPDLeaf_{A}}{x} \myBinCoprod \myApp{\myPDRoot_{A}}{x}
        & \text{if \(\myApp{\myDegree}{x} \myId n - 1\)} \\
        \myTerminal
        & \text{if \(\myApp{\myDegree}{x} \myLt n - 1\)}
      \end{array}
    \right.
  \end{equation*}
  We extend \(\myApp{\myBoundary}{A}\) to an \(n\)-preboundary by
  \(\myApp{\myBdSource}{\myApp{\myBoundary}{A}} \myDefEq
  \myApp{\myPDLeaf}{A}\) and
  \(\myApp{\myBdTarget}{\myApp{\myBoundary}{A}} \myDefEq
  \myApp{\myPDRoot}{A}\). When \(n \myGe 1\),
  \(\myApp{\mySourceHorn}{\myApp{\myBoundary}{A}}\) is abbreviated to
  \(\myApp{\mySourceHorn}{A}\).
\end{myConstruction}

\begin{myDefinition}
  Let \(n \myGe 0\). We mutually define \(n\)-boundaries and
  \(n\)-pasting diagrams. An \myDefine{\(n\)-boundary} is an
  \(n\)-preboundary \(A\) satisfying the following axioms.
  \begin{myWithCustomLabel}{axiom}
    \begin{enumerate}[label=\textbf{Bd\arabic*.},ref=Bd\arabic*]
    \item \label{ax-bd-target} When \(n \myGe 1\), there exists a
      unique target object in \(A\).
    \item \label{ax-bd-horn} When \(n \myGe 1\),
      \(\myApp{\mySourceHorn}{A}\) is an \((n - 1)\)-pasting diagram.
    \end{enumerate}
  \end{myWithCustomLabel}
  An \myDefine{\(n\)-pasting diagram} is an \(n\)-prepasting diagram
  \(A\) satisfying the following axioms.
  \begin{myWithCustomLabel}{axiom}
    \begin{enumerate}[label=\textbf{PD\arabic*.},ref=PD\arabic*]
    \item \label{ax-pd-first} \label{ax-pd-finite} \(A_{n}\) is
      finite.
    \item \label{ax-pd-leaf} For every \(x \myElemOf A_{n - 1}\),
      the type \(\myApp{\myPDLeaf_{A}}{x}\) is a proposition, and it
      holds if and only if there is no target arrow from \(x\).
    \item \label{ax-pd-root} For every \(x \myElemOf A_{n - 1}\),
      the type \(\myApp{\myPDRoot_{A}}{x}\) is a proposition, and it
      holds if and only if there is no source arrow from \(x\).
    \item \label{ax-pd-degree0} When \(n \myId 0\), \(A_{0}\) is
      contractible.
    \item \label{ax-pd-target} For every object
      \(x \myElemOf A_{n - 1}\), there is at most one target arrow from
      \(x\).
    \item \label{ax-pd-source} For every object \(x \myElemOf A_{n -
        1}\), there is at most one source arrow from \(x\).
    \item \label{ax-pd-path} When \(n \myGe 1\), there exists an
      object \(r \myElemOf A_{n - 1}\) such that, for every object
      \(x \myElemOf A_{n - 1}\), there exists a zigzag
      \begin{equation}
        \label{eq-pd-zigzag}
        x \myId x_{0} \myXMorphism{f_{0}} y_{0} \myXMorphismAlt{g_{0}}
        x_{1} \myXMorphism{f_{1}} \myDots \myXMorphism{f_{m - 1}} y_{m -
          1} \myXMorphismAlt{g_{m - 1}} x_{m} \myId r,
      \end{equation}
      where \(x_{i} \myElemOf A_{n - 1}\), \(y_{i} \myElemOf A_{n}\),
      \(f_{i}\)'s are source arrows, and \(g_{i}\)'s are target arrows.
      \label{ax-pd-last-but-bd}
    \item \label{ax-pd-boundary} \(\myApp{\myBoundary}{A}\) is an
      \(n\)-boundary.
    \end{enumerate}
  \end{myWithCustomLabel}
  Let \(\myNBd{n}\) denote the category of small \(n\)-boundaries
  whose morphisms are those morphisms of opetopic sets preserving
  source and target objects. Let \(\myNPD{n}\) denote the category of
  small \(n\)-pasting diagrams whose morphisms are those morphisms of
  opetopic sets preserving leaf and root objects. By definition,
  \(\mySourceHorn\) is a functor
  \(\myNBd{n} \myMorphism \myNPD{n - 1}\) when \(n \myGe 1\), and
  \(\myBoundary\) is a functor \(\myNPD{n} \myMorphism \myNBd{n}\).
\end{myDefinition}

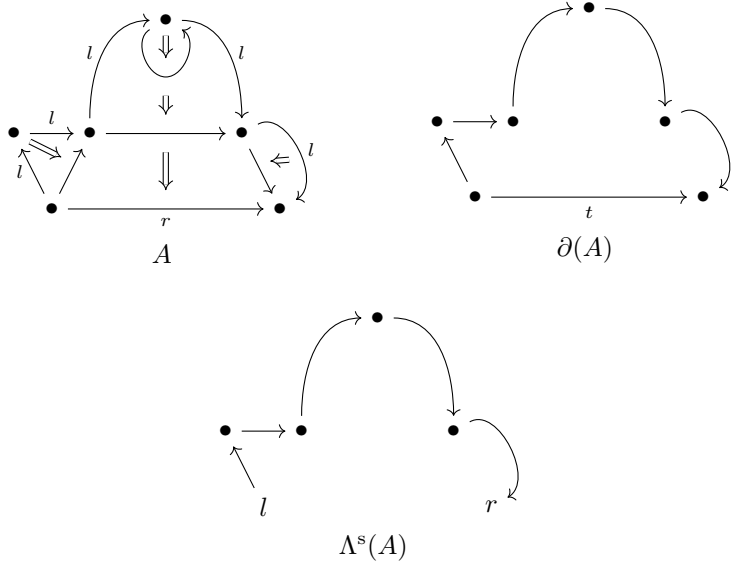
\begin{figure}
  \begin{mathpar}
    \begin{tikzpicture}
      \node (x0) at (0, 0) {\(\bullet\)};
      \node (x1) at (-0.5, 1) {\(\bullet\)};
      \node (x2) at (0.5, 1) {\(\bullet\)};
      \node (x3) at (1.5, 2.5) {\(\bullet\)};
      \node (x4) at (2.5, 1) {\(\bullet\)};
      \node (x5) at (3, 0) {\(\bullet\)};
      \draw[->] (x0) to node[left] {\scriptsize \(l\)} (x1);
      \draw[->] (x0) to node[above,name=y0] {} (x2);
      \draw[->] (x0) to node[above,name=y1] {}
      node[below] {\scriptsize \(r\)} (x5);
      \draw[->] (x1) to node[above] {\scriptsize \(l\)} (x2);
      \draw[->,out=90,in=180] (x2) to node[left] {\scriptsize \(l\)} (x3);
      \draw[->] (x2) to node[above,name=y2] {} node[below,name=y3] {} (x4);
      \draw[->] (x3)
      to [out=210,in=180] (1.5, 1.75)
      to node[above,name=y4] {} node[below,name=y5] {} (1.5, 1.75)
      to [out=0,in=330] (x3);
      \draw[->,out=0,in=90] (x3) to node[right] {\scriptsize \(l\)} (x4);
      \draw[->] (x4) to node[above,name=y6] {} (x5);
      \draw[->,out=30,in=30] (x4) to node[below,name=y7] {}
      node[right] {\scriptsize \(l\)} (x5);
      \draw[double equal sign distance,-Implies] (x1) to (y0);
      \draw[double equal sign distance,-Implies] (x3) to (y4);
      \draw[double equal sign distance,-Implies] (y5) to (y2);
      \draw[double equal sign distance,-Implies] (y7) to (y6);
      \draw[double equal sign distance,-Implies] (y3) to (y1);
      \node[below] (title) at (current bounding box.south) {\(A\)};
    \end{tikzpicture}
    \and
    \begin{tikzpicture}
      \node (x0) at (0, 0) {\(\bullet\)};
      \node (x1) at (-0.5, 1) {\(\bullet\)};
      \node (x2) at (0.5, 1) {\(\bullet\)};
      \node (x3) at (1.5, 2.5) {\(\bullet\)};
      \node (x4) at (2.5, 1) {\(\bullet\)};
      \node (x5) at (3, 0) {\(\bullet\)};
      \draw[->] (x0) to (x1);
      \draw[->] (x0) to node[above,name=y1] {}
      node[below] {\scriptsize \(t\)} (x5);
      \draw[->] (x1) to (x2);
      \draw[->,out=90,in=180] (x2) to (x3);
      \draw[->,out=0,in=90] (x3) to (x4);
      \draw[->,out=30,in=30] (x4) to node[below,name=y7] {} (x5);
      \node[below] (title) at (current bounding box.south)
      {\(\myApp{\myBoundary}{A}\)};
    \end{tikzpicture}
    \and
    \begin{tikzpicture}
      \node (x0) at (0, 0) {\(l\)};
      \node (x1) at (-0.5, 1) {\(\bullet\)};
      \node (x2) at (0.5, 1) {\(\bullet\)};
      \node (x3) at (1.5, 2.5) {\(\bullet\)};
      \node (x4) at (2.5, 1) {\(\bullet\)};
      \node (x5) at (3, 0) {\(r\)};
      \draw[->] (x0) to (x1);
      \draw[->] (x1) to (x2);
      \draw[->,out=90,in=180] (x2) to (x3);
      \draw[->,out=0,in=90] (x3) to (x4);
      \draw[->,out=30,in=30] (x4) to node[below,name=y7] {} (x5);
      \node[below] (title) at (current bounding box.south)
      {\(\myApp{\mySourceHorn}{A}\)};
    \end{tikzpicture}
  \end{mathpar}
  \caption{Examples of pasting diagrams and boundaries. \(A\) is a
    \(2\)-pasting diagram, \(\myApp{\myBoundary}{A}\) is its boundary,
    and \(\myApp{\mySourceHorn}{A}\) is its source horn. Leaves in
    pasting diagrams are marked with ``\(l\)''. The roots in pasting
    diagrams are marked with ``\(r\)''. The targets in boundaries are
    marked with ``\(t\)''.}
  \label{fig-example-pd-bd}
\end{figure}

We illustrate some examples of pasting diagrams and boundaries in
\cref{fig-example-pd-bd}. \Crefrange{ax-pd-first}{ax-pd-last-but-bd}
express that a pasting diagram forms a tree; see
\cref{fig-pd-tree}. \Cref{ax-pd-boundary} is necessary to exclude, for
example, the following diagram.
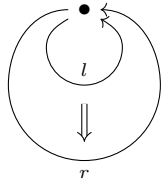
\begin{mathpar}
  \begin{tikzpicture}
    \node (x0) at (0, 0) {\(\bullet\)};
    \draw[->] (x0)
    to[out=-150,in=90] (-0.5, -0.5)
    to[out=-90,in=180] (0, -1)
    to node[above] {\scriptsize \(l\)} node[below,name=y1] {} (0, -1)
    to[out=0,in=-90] (0.5, -0.5)
    to[out=90,in=-30] (x0);
    \draw[->] (x0)
    to[out=180,in=90] (-1, -1)
    to[out=-90,in=180] (0, -2)
    to node[above,name=y2] {} node[below] {\scriptsize \(r\)} (0, -2)
    to[out=0,in=-90] (1, -1)
    to[out=90,in=0] (x0);
    \draw[double equal sign distance,-Implies] (y1) to (y2);
  \end{tikzpicture}
\end{mathpar}
It forms a tree with root \(r\) and leaf \(l\), but it is not
considered as a pasting diagram due to the hole surrounded by \(l\).
There are \(n\)-pasting diagrams \(A\) such that \(A_{n}\) is empty,
in which case the root object in \(A\) is also a leaf object. For
example, the following is a \(2\)-pasting diagram.
\begin{mathpar}
  \begin{tikzpicture}
    \node (x0) at (0, 0) {\(\bullet\)};
    \node (x1) at (1, 0) {\(\bullet\)};
    \draw[->] (x0) to (x1);
  \end{tikzpicture}
\end{mathpar}
Its boundary is the following, where the target object is marked with
``\(t\)''.
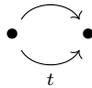
\begin{mathpar}
  \begin{tikzpicture}
    \node (x0) at (0, 0) {\(\bullet\)};
    \node (x1) at (1, 0) {\(\bullet\)};
    \draw[->,out=60,in=120] (x0) to (x1);
    \draw[->,out=-60,in=-120] (x0) to node[below] {\scriptsize \(t\)} (x1);
  \end{tikzpicture}
\end{mathpar}

We construct the boundary of an opetope
(\cref{prop-opetope-boundary}).

\begin{myConstruction}
  Let \(A\) be an opetope of degree \(n \myGe 0\). We define
  \(\myApp{\myBoundary}{A} \myDefEq A_{\myLt n}\). We extend it to an
  \(n\)-preboundary where the source/target objects are those
  \(x \myElemOf A_{n - 1}\) such that the unique arrow
  \(x \myOneMorphism \myOpeTerminal_{A}\) is a source/target arrow. We
  call \(\myApp{\myBoundary}{A}\) the \myDefine{boundary of
    \(A\)}. When \(n \myGe 1\),
  \(\myApp{\mySourceHorn}{\myApp{\myBoundary}{A}}\) is abbreviated to
  \(\myApp{\mySourceHorn}{A}\) and called the \myDefine{source horn of
    \(A\)}.
\end{myConstruction}

\begin{myLemma}
  \label{prop-opetope-horn-boundary}
  Let \(A\) be an opetope of degree \(n \myGe 1\). Then
  \(\myApp{\myBoundary}{\myApp{\mySourceHorn}{A}} \myId
  \myApp{\myBoundary}{\myApp{\myDomTargetOf_{\myOpetope}}{A}}\) in the
  type of \((n - 1)\)-preboundaries.
\end{myLemma}
\begin{proof}
  \(\myApp{\myBoundary}{\myApp{\mySourceHorn}{A}}\) is the category of
  elements for the following set-valued presheaf on \(A\).
  \begin{equation*}
    x \myMapsTo \left\{
      \begin{array}{ll}
        \myInitial & \text{if \(\myApp{\myDegree}{x} \myGe n - 1\)} \\
        (x \mySourceMorphism
        \myApp{\myDomTargetOf}{\myOpeTerminal_{A}}) \myBinCoprod (x
        \myTargetMorphism \myApp{\myDomTargetOf}{\myOpeTerminal_{A}})
                   & \text{if \(\myApp{\myDegree}{x} \myId n - 2\)} \\
        \myTerminal & \text{if \(\myApp{\myDegree}{x} \myLt n - 2\)}
      \end{array}
    \right.
  \end{equation*}
  \(\myApp{\myBoundary}{\myApp{\myDomTargetOf_{\myOpetope}}{A}}\) is
  the category of elements for the following set-valued presheaf on
  \(A\).
  \begin{equation*}
    x \myMapsTo \left\{
      \begin{array}{ll}
        \myInitial & \text{if \(\myApp{\myDegree}{x} \myGe n - 1\)} \\
        (x \mySourceMorphism
        \myApp{\myDomTargetOf}{\myOpeTerminal_{A}}) \myBinCoprod (x
        \myTargetMorphism \myApp{\myDomTargetOf}{\myOpeTerminal_{A}})
                   & \text{if \(\myApp{\myDegree}{x} \myId n - 2\)} \\
        \myApp{\myArr_{A}}{x \myComma
        \myApp{\myDomTargetOf}{\myOpeTerminal_{A}}}
                   & \text{if \(\myApp{\myDegree}{x} \myLt n - 2\)}
      \end{array}
    \right.
  \end{equation*}
  They are equivalent because
  \(\myApp{\myArr_{A}}{x \myComma
    \myApp{\myDomTargetOf}{\myOpeTerminal_{A}}} \myEquiv
  \myApp{\myKNF{k - 1}_{A}}{x \myComma
    \myApp{\myDomTargetOf}{\myOpeTerminal_{A}}} \myEquiv
  \myApp{\myKNF{k}_{A}}{x \myComma \myOpeTerminal_{A}} \myEquiv
  \myApp{\myArr_{A}}{x \myComma \myOpeTerminal_{A}}\)
  (\cref{prop-oset-normalization}) is contractible when
  \(\myApp{\myDegree}{x} \myId n - k\) for \(k \myGe 3\).
\end{proof}

\begin{myProposition}
  \label{prop-opetope-boundary}
  \(\myApp{\myBoundary}{A}\) is an \(n\)-boundary for every opetope
  \(A\) of degree \(n \myGe 0\).
\end{myProposition}
\begin{proof}
  We proceed by induction on \(n\). \Cref{ax-bd-target} is by
  \cref{ax-oset-target}. Suppose that \(n \myGe 1\) and we verity the
  \((n - 1)\)-pasting diagram axioms for
  \(\myApp{\mySourceHorn}{A}\). \Cref{ax-pd-finite} is by
  \cref{ax-oset-finite}. For \(x \myElemOf A_{n - 2}\), the type
  \(\myApp{\myPDLeaf_{\myApp{\mySourceHorn}{A}}}{x}\) is the type of
  factorizations of \(x \myTwoMorphism \myOpeTerminal_{A}\) into a
  source arrow followed by a target arrow. Thus, \cref{ax-pd-leaf}
  follows from \cref{ax-oset-hetero}. Similarly \cref{ax-pd-root}
  follows from \cref{ax-oset-homo}. \Cref{ax-pd-degree0} is by
  \cref{ax-oset-source}. \Cref{ax-pd-target} follows from
  \cref{ax-oset-hetero}. \Cref{ax-pd-source} follows from
  \cref{ax-oset-homo}. \Cref{ax-pd-path} is by
  \cref{ax-oset-path}. \Cref{ax-pd-boundary} is by
  \cref{prop-opetope-horn-boundary} and induction hypothesis.
\end{proof}

We then show that
\(\myBoundary \myElemOf \myOpetope_{n} \myMorphism \myNBd{n}\) is an
equivalence (\cref{prop-boundary-equiv}). We prepare basic lemmas for
boundaries and pasting diagrams.

\begin{myConstruction}
  Let \(n \myGe 0\) and let \(A\) be an \(n\)-prepasting diagram. We
  define a \(0\)-graph \(\myApp{\myPDGraph}{A}\) as follows. The set
  of vertices in \(\myApp{\myPDGraph}{A}\) is
  \(A_{n} \myBinCoprod A_{n - 1}\). There is no edge between vertices
  from \(A_{n}\). There is no edge between vertices from
  \(A_{n - 1}\). An edge from \(x \myElemOf A_{n}\) to
  \(y \myElemOf A_{n - 1}\) is a target arrow
  \(f \myElemOf y \myTargetMorphism x\). An edge from
  \(y \myElemOf A_{n - 1}\) to \(x \myElemOf A_{n}\) is a source arrow
  \(f \myElemOf y \mySourceMorphism x\).
\end{myConstruction}

\begin{myLemma}
  \label{prop-pd-tree}
  Let \(n \myGe 0\) and let \(A\) be an \(n\)-pasting diagram.
  \begin{enumerate}
  \item \(\myApp{\myPDGraph}{A}\) is a tree.
  \item When \(n \myId 0\), \(\myApp{\myPDGraph}{A}\) is a singleton
    set with no edge.
  \item When \(n \myGe 1\), the root in \(\myApp{\myPDGraph}{A}\) is
    the unique root object in \(A\), which exists by
    \cref{ax-bd-target} for \(\myApp{\myBoundary}{A}\).
  \end{enumerate}
\end{myLemma}
\begin{proof}
  When \(n \myId 0\), \(\myApp{\myPDGraph}{A}\) is a singleton set
  with no edge by \cref{ax-pd-degree0} and by definition. In
  particular, it is a tree. Suppose that \(n \myGe 1\). Take
  \(r \myElemOf A_{n - 1}\) as in \cref{ax-pd-path}. Then Zigzag
  \labelcref{eq-pd-zigzag} is a path from \(x\) to \(r\) in
  \(\myApp{\myPDGraph}{A}\). Let \(t\) be the unique root object in
  \(A\). There is a path \(p\) from \(t\) to \(r\). By
  \cref{ax-pd-root}, the length of \(p\) must be \(0\), and thus
  \(t \myId r\). For every \(x \myElemOf A_{n - 1}\), there is a path
  from \(x\) to \(r\). Such a path is unique by
  \cref{ax-pd-source,ax-oset-target,ax-pd-root}. For every
  \(x \myElemOf A_{n}\), we have a unique edge \(x \myMorphism y\) by
  \cref{ax-oset-target} and then a unique path from \(y\) to \(r\).
\end{proof}

\begin{myLemma}
  \label{prop-pd-root}
  Let \(n \myGe 1\), let \(A\) be an \(n\)-pasting diagram, and let
  \(x \myElemOf A_{n - 1}\). Then either \(x\) is a root object or
  there exists a unique source arrow from \(x\).
\end{myLemma}
\begin{proof}
  Let \(d\) be the length of the unique path in
  \(\myApp{\myPDGraph}{A}\) from \(x\) to the root
  (\cref{prop-pd-tree}). If \(d \myId 0\), then there is no source
  arrow from \(x\), and thus \(x\) is a root object by
  \cref{ax-pd-root}. If \(d \myGe 1\), then there is a source arrow
  from \(x\), and such a source arrow is unique by
  \cref{ax-pd-source}.
\end{proof}

\begin{myLemma}
  \label{prop-pd-finite}
  Let \(n \myGe 0\) and let \(A\) be an \(n\)-pasting diagram. Then
  \(A_{n - 1}\) is finite.
\end{myLemma}
\begin{proof}
  The case when \(n \myId 0\) is trivial since \(A_{-1}\) is
  empty. Suppose that \(n \myGe 1\). We show that
  \(A_{n - 1} \myEquiv \myTerminal \myBinCoprod ((x \myElemOf A_{n})
  \myBinProd (A \mySourceSlice x))\), which is finite by
  \cref{ax-pd-finite,ax-oset-finite}. Let \(x \myElemOf A_{n -
    1}\). We proceed by case analysis on \(x\) by
  \cref{prop-pd-root}. If \(x\) is a root object, we map \(x\) to
  \(\myTerminalEl \myElemOf \myTerminal\). If there is a unique source
  arrow \(f \myElemOf x \mySourceMorphism y\), we map \(x\) to
  \((y \myComma f) \myElemOf (x' \myElemOf A_{n}) \myBinProd (A
  \mySourceSlice x')\). This gives a one-to-one correspondence between
  \(A_{n - 1}\) and
  \(\myTerminal \myBinCoprod ((x \myElemOf A_{n}) \myBinProd (A
  \mySourceSlice x))\).
\end{proof}

\begin{myLemma}
  \label{prop-pd-leaf}
  Let \(n \myGe 1\), let \(A\) be an \(n\)-pasting diagram, and let
  \(x \myElemOf A_{n - 1}\). Then either \(x\) is a leaf object or
  there exists a unique target arrow from \(x\).
\end{myLemma}
\begin{proof}
  By \cref{ax-pd-finite,prop-pd-finite}, the proposition ``there is a
  target arrow from \(x\)'' is decidable: check if
  \(\myApp{\myDomTargetOf}{y} \myId x\) for all \(y \myElemOf
  A_{n}\). Thus, by \cref{ax-pd-leaf}, either \(x\) is a leaf object
  or there is a target arrow from \(x\). In the latter case, such a
  target arrow is unique by \cref{ax-pd-target}.
\end{proof}

\begin{myLemma}
  \label{prop-oset-kstep-factor}
  Let \(A\) be a preopetopic set satisfying
  \crefrange{ax-oset-first}{ax-oset-core-last}. Then \(A\) satisfies
  \cref{ax-oset-kstep-unique,ax-oset-kstep-factor} if and only if, for
  every \(k \myGe 3\) and \(x \myComma y \myElemOf A\) such that
  \(\myApp{\myDegree}{y} + k \myId \myApp{\myDegree}{x}\), the
  postcomposition map
  \(\myApp{\myTargetOf}{x}_{\myBang} \myElemOf \myApp{\myArr_{A}}{y
    \myComma \myApp{\myDomTargetOf}{x}} \myMorphism
  \myApp{\myArr_{A}}{y \myComma x}\) is an equivalence.
\end{myLemma}
\begin{proof}
  The ``only if'' direction follows from
  \cref{prop-oset-normalization} as
  \(\myApp{\myArr_{A}}{y \myComma \myApp{\myDomTargetOf}{x}} \myEquiv
  \myApp{\myKNF{k - 1}_{A}}{y \myComma \myApp{\myDomTargetOf}{x}}
  \myEquiv \myApp{\myKNF{k}_{A}}{y \myComma x} \myEquiv
  \myApp{\myArr_{A}}{y \myComma x}\). We show the ``if''
  direction. \Cref{ax-oset-kstep-unique} is immediate. By assumption
  and \cref{ax-oset-homo}, we have
  \(\myApp{\myKNF{k}_{A}}{y \myComma x} \myEquiv \myApp{\myArr_{A}}{y
    \myComma x}\) for every \(k \myGe 2\) and
  \(x \myComma y \myElemOf A\) such that
  \(\myApp{\myDegree}{y} + k \myId
  \myApp{\myDegree}{x}\). \Cref{ax-oset-kstep-factor} thus follows.
\end{proof}

\begin{myConstruction}
  Let \(n \myGe 1\) and let \(A\) be an \(n\)-boundary. We refer to
  the unique target object in \(A\), which exists by
  \cref{ax-bd-target}, as \(\myApp{\myBdTargetOf}{A}\).
\end{myConstruction}

\begin{myLemma}
  \label{prop-bd-arr-to-target}
  Let \(k \myGe 2\), let \(n \myGe 1\), let \(A\) be an
  \(n\)-boundary, and let \(x \myElemOf A\) be an object of degree
  \(n - 1 - k\). Then there exists a unique \(k\)-step arrow
  \(x \myKMorphism{k} \myApp{\myBdTargetOf}{A}\).
\end{myLemma}
\begin{proof}
  By induction on \(k \myGe 2\). Suppose that \(k \myId 2\). Then
  \(x\) is an object in
  \(\myApp{\mySourceHorn}{\myApp{\mySourceHorn}{A}}\) of degree
  \(n - 3\). By \cref{prop-pd-root}, either there is a unique target
  arrow
  \(f \myElemOf x \myTargetMorphism
  \myApp{\myBdTargetOf}{\myApp{\myBoundary}{\myApp{\mySourceHorn}{A}}}\)
  or there is a unique source arrow
  \(f \myElemOf x \mySourceMorphism y\) to a source object \(y\) in
  \(\myApp{\myBoundary}{\myApp{\mySourceHorn}{A}}\). In the former
  case, by construction,
  \(\myApp{\myBdTargetOf}{\myApp{\myBoundary}{\myApp{\mySourceHorn}{A}}}
  \myId \myApp{\myDomTargetOf}{\myApp{\myBdTargetOf}{A}}\), and thus
  we have the unique target arrow
  \(g \myDefEq \myApp{\myTargetOf}{\myApp{\myBdTargetOf}{A}} \myElemOf
  \myApp{\myBdTargetOf}{\myApp{\myBoundary}{\myApp{\mySourceHorn}{A}}}
  \myTargetMorphism \myApp{\myBdTargetOf}{A}\). In the latter case,
  \(y\) is a leaf object in \(\myApp{\mySourceHorn}{A}\) by the
  definition of
  \(\myApp{\myBoundary}{\myApp{\mySourceHorn}{A}}\). Then we have a
  unique source arrow
  \(g \myElemOf y \mySourceMorphism \myApp{\myBdTargetOf}{A}\). In
  both cases, \((g \myComma f)\) is the unique homogeneous pair of
  \(1\)-step arrows
  \(f \myElemOf y \myOneMorphism \myApp{\myBdTargetOf}{A}\) and
  \(g \myElemOf x \myOneMorphism y\). Then the composite
  \(g \myComp f\) is the unique arrow
  \(x \myTwoMorphism \myApp{\myBdTargetOf}{A}\) by
  \cref{ax-oset-homo}. Suppose that \(k \myGe 3\). As we have seen,
  \(\myApp{\myBdTargetOf}{\myApp{\myBoundary}{\myApp{\mySourceHorn}{A}}}
  \myId \myApp{\myDomTargetOf}{\myApp{\myBdTargetOf}{A}}\). By
  induction hypothesis, we have a unique arrow
  \(f \myElemOf x \myKMorphism{k - 1}
  \myApp{\myBdTargetOf}{\myApp{\myBoundary}{\myApp{\mySourceHorn}{A}}}\). Then
  the composite
  \(\myApp{\myTargetOf}{\myApp{\myBdTargetOf}{A}} \myComp f\) is the
  unique arrow \(x \myKMorphism{k} \myApp{\myBdTargetOf}{A}\) by
  \cref{prop-oset-kstep-factor}.
\end{proof}

We now construct an inverse of
\(\myBoundary \myElemOf \myOpetope_{n} \myMorphism \myNBd{n}\).

\begin{myConstruction}
  Let \(n \myGe 0\) and let \(A\) be an \(n\)-preboundary. We
  construct a category \(\myApp{\myFill}{A}\) from \(A\) by freely
  adjoining a terminal object \(\myOpeTerminal\). We extend it to an
  \(\myFinOrd\)-direct category by extending \(\myDegree_{A}\) by
  \(\myApp{\myDegree_{\myApp{\myFill}{A}}}{\myOpeTerminal} \myDefEq
  n\). We further extend it to a preopetopic set where an arrow is a
  source/target arrow if eigher it is a source/target arrow in \(A\)
  or it is the arrow \(x \myOneMorphism \myOpeTerminal\) from a
  source/target object \(x\) in \(A\).
\end{myConstruction}

\begin{myLemma}
  \label{prop-bd-fill-opetope}
  Let \(n \myGe 0\) and let \(A\) be an \(n\)-boundary. Then
  \(\myApp{\myFill}{A}\) is an opetope of degree \(n\).
\end{myLemma}
\begin{proof}
  Since \(\myApp{\myFill}{A} \mySlice x \myEquiv A \mySlice x\) for
  \(x \myElemOf A\) by construction, it suffices to verify
  \crefrange{ax-oset-first}{ax-oset-last} for
  \(x \myDefEq \myOpeTerminal\). Since
  \(\myApp{\myFill}{A} \myOneSlice \myOpeTerminal \myEquiv A_{n - 1}\)
  by construction, \cref{ax-oset-finite} is trivial when \(n \myId 0\)
  and follows from \cref{ax-bd-target} and \cref{ax-pd-finite} for
  \(\myApp{\mySourceHorn}{A}\) when \(n \myGe
  1\). \Cref{ax-oset-target} is by
  \cref{ax-bd-target}. \Cref{ax-oset-source} is by
  \cref{ax-pd-degree0} for
  \(\myApp{\mySourceHorn}{A}\). \Cref{ax-oset-homo} follows from
  \cref{prop-pd-root} for
  \(\myApp{\mySourceHorn}{A}\). \Cref{ax-oset-hetero} follows from
  \cref{prop-pd-leaf} for
  \(\myApp{\mySourceHorn}{A}\). \Cref{ax-oset-path} is by
  \cref{ax-pd-path} for
  \(\myApp{\mySourceHorn}{A}\). \Cref{ax-oset-kstep-unique,ax-oset-kstep-factor}
  follow from \cref{prop-oset-kstep-factor,prop-bd-arr-to-target}.
\end{proof}

\begin{myProposition}
  \label{prop-boundary-equiv}
  Let \(n \myGe 0\). Then the functor
  \(\myBoundary \myElemOf \myOpetope_{n} \myMorphism \myNBd{n}\) is an
  equivalence.
\end{myProposition}
\begin{proof}
  The inverse is given by \(\myFill\) (\cref{prop-bd-fill-opetope}).
\end{proof}

\begin{myCorollary}
  \label{prop-boundary-discrete}
  Let \(n \myGe 0\). Then \(\myNBd{n}\) is a discrete category over a
  set.
\end{myCorollary}
\begin{proof}
  By \cref{prop-boundary-equiv}.
\end{proof}

We show that
\(\mySourceHorn \myElemOf \myNBd{n + 1} \myMorphism \myNPD{n}\) is an
equivalence (\cref{prop-bd-equiv-pd}).

\begin{myLemma}
  \label{prop-bd-horn-boundary}
  Let \(n \myGe 0\). Then the diagram
  \begin{equation}
    \label{eq-bd-horn-boundary}
    \begin{tikzcd}
      \myNBd{n + 1}
      \arrow[r, "\mySourceHorn"]
      \arrow[d, "\myBdTargetOf"'] &
      \myNPD{n}
      \arrow[d, "\myBoundary"] \\
      \myOpetope_{n}
      \arrow[r, "\myBoundary"'] &
      \myNBd{n}
    \end{tikzcd}
  \end{equation}
  is a pullback of categories, where we identify
  \(\myApp{\myBdTargetOf}{A} \myElemOf A_{n}\) and the associated
  opetope \(A \mySlice \myApp{\myBdTargetOf}{A}\) for
  \(A \myElemOf \myNBd{n + 1}\).
\end{myLemma}
\begin{proof}
  We first note that Diagram \labelcref{eq-bd-horn-boundary} commutes,
  that is,
  \(\myApp{\myBoundary}{\myApp{\mySourceHorn}{A}} \myId
  \myApp{\myBoundary}{\myApp{\myBdTargetOf}{A}}\) for all
  \(A \myElemOf \myNBd{n + 1}\). By \cref{prop-boundary-equiv}, it
  suffices to show the case when \(A\) is of the form
  \(\myApp{\myBoundary}{A'}\) for an opetope \(A'\) of degree
  \(n + 1\), but this is just \cref{prop-opetope-horn-boundary}. To
  see that Diagram \labelcref{eq-bd-horn-boundary} is a pullback, let
  \(A \myElemOf \myOpetope_{n}\), \(B \myElemOf \myNPD{n}\), and
  \(C \myElemOf \myNBd{n}\) and suppose
  \(\myApp{\myBoundary}{A} \myId \myApp{\myBoundary}{B} \myId C\). Let
  \(\myApp{D}{A \myComma B}\) be the following pushout in \(\myOSet\).
  \begin{equation*}
    \begin{tikzcd}
      C
      \arrow[r]
      \arrow[d]
      \arrow[dr, myPOMark] &
      B
      \arrow[d] \\
      A
      \arrow[r] &
      \myApp{D}{A \myComma B}
    \end{tikzcd}
  \end{equation*}
  By \cref{prop-oset-fiber-colimit}, we have
  \(\myTerminal \myBinCoprod B_{n} \myEquiv \myApp{D}{A \myComma
    B}_{n}\). We extend \(\myApp{D}{A \myComma B}\) to an
  \((n + 1)\)-preboundary by
  \(\myApp{\myBdSource}{\myApp{D}{A \myComma B}} \myDefEq B_{n}\) and
  \(\myApp{\myBdTarget}{\myApp{D}{A \myComma B}} \myDefEq
  \myTerminal\). By \cref{prop-oset-restrict-colimit}, we have
  \(B_{\myLt n} \myEquiv \myApp{D}{A \myComma B}_{\myLt n}\). It then
  follows that
  \(B \myEquiv \myApp{\mySourceHorn}{\myApp{D}{A \myComma B}}\) by
  \cref{prop-oset-fiber-conservative}, and thus
  \(\myApp{D}{A \myComma B}\) is an \((n + 1)\)-boundary. By
  \cref{prop-mor-oset-slice-equiv},
  \(A \myEquiv \myApp{\myBdTargetOf}{\myApp{D}{A \myComma B}}\). We
  thus have a section \(D\) of
  \(\myNBd{n + 1} \myMorphism \myOSet_{n} \myBinProd_{\myNBd{n}}
  \myNPD{n}\). To see that \(D\) is moreover an inverse, let
  \(X \myElemOf \myNBd{n + 1}\). We have a canonical morphism
  \(\myApp{D}{\myApp{\myBdTargetOf}{X} \myComma
    \myApp{\mySourceHorn}{X}} \myMorphism X\) of
  \((n + 1)\)-boundaries, which is an equivalence by
  \cref{prop-boundary-discrete}.
\end{proof}

\begin{myProposition}
  \label{prop-bd-equiv-pd}
  Let \(n \myGe 1\). Then the functor
  \(\mySourceHorn \myElemOf \myNBd{n} \myMorphism \myNPD{n - 1}\) is
  an equivalence.
\end{myProposition}
\begin{proof}
  By \cref{prop-boundary-equiv,prop-bd-horn-boundary}.
\end{proof}

\begin{myCorollary}
  \label{prop-opetope-equiv-pd}
  Let \(n \myGe 1\). Then the map
  \(\mySourceHorn \myElemOf \myOpetope_{n} \myMorphism \myNPD{n - 1}\)
  is an equivalence.
\end{myCorollary}
\begin{proof}
  By \cref{prop-boundary-equiv,prop-bd-equiv-pd}.
\end{proof}

\begin{myCorollary}
  \label{prop-pd-discrete}
  Let \(n \myGe 0\). Then \(\myNPD{n}\) is a discrete category over a
  set.
\end{myCorollary}
\begin{proof}
  By \cref{prop-opetope-equiv-pd}.
\end{proof}

\section{Substitution and grafting}
\label{sec:subst-graft}

We introduce two operators on pasting diagrams, substitution and
grafting. Substitution is the operator on \(n\)-pasting diagrams
replacing \(n\)-cells in an \(n\)-pasting diagram by \(n\)-pasting
diagrams of the same boundaries. For example, let \(A\),
\(\myApp{B}{x_{0}}\), \(\myApp{B}{x_{1}}\), and \(\myApp{B}{x_{2}}\)
be the following \(2\)-pasting diagrams.
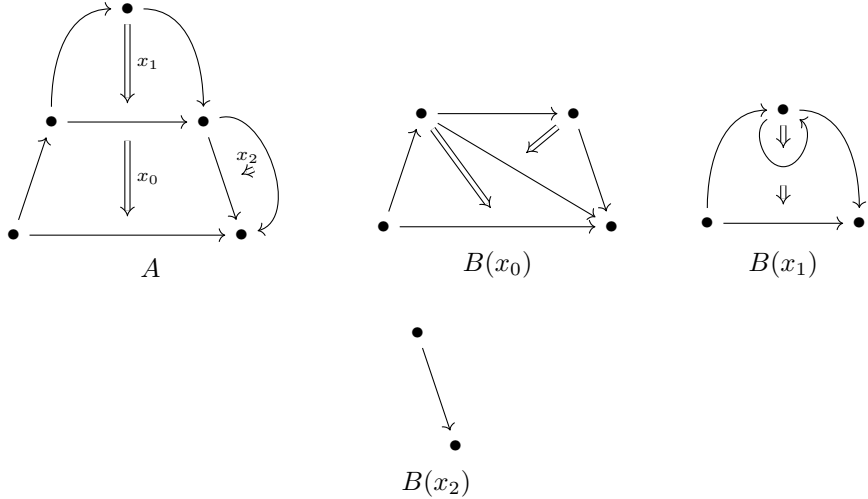
\begin{mathpar}
  \begin{tikzpicture}
    \node (x0) at (0, 0) {\(\bullet\)};
    \node (x1) at (0.5, 1.5) {\(\bullet\)};
    \node (x2) at (2.5, 1.5) {\(\bullet\)};
    \node (x3) at (3, 0) {\(\bullet\)};
    \node (x4) at (1.5, 3) {\(\bullet\)};
    \draw[->] (x0) to (x1);
    \draw[->] (x0) to node[above,name=y0] {} (x3);
    \draw[->,out=90,in=180] (x1) to (x4);
    \draw[->] (x1) to node[below,name=y1] {}
    node[above,name=y2] {} (x2);
    \draw[->] (x2) to node[right,name=y6] {} (x3);
    \draw[->,out=15,in=15] (x2) to node[left,name=y7] {} (x3);
    \draw[->,out=0,in=90] (x4) to (x2);
    \draw[double equal sign distance,-Implies] (y1) to
    node[right] {\scriptsize \(x_{0}\)} (y0);
    \draw[double equal sign distance,-Implies] (x4) to
    node[right] {\scriptsize \(x_{1}\)} (y2);
    \draw[double equal sign distance,-Implies] (y7) to
    node[above] {\scriptsize \(x_{2}\)} (y6);
    \node[below] (title) at (current bounding box.south) {\(A\)};
  \end{tikzpicture}
  \and
  \begin{tikzpicture}
    \node (x0) at (0, 0) {\(\bullet\)};
    \node (x1) at (0.5, 1.5) {\(\bullet\)};
    \node (x2) at (2.5, 1.5) {\(\bullet\)};
    \node (x3) at (3, 0) {\(\bullet\)};
    \draw[->] (x0) to (x1);
    \draw[->] (x0) to node[above,name=y0] {} (x3);
    \draw[->] (x1) to node[below,name=y1] {}
    node[above,name=y2] {} (x2);
    \draw[->] (x2) to (x3);
    \draw[->] (x1) to node[above,name=y3] {} (x3);
    \draw[double equal sign distance,-Implies] (x1) to (y0);
    \draw[double equal sign distance,-Implies] (x2) to (y3);
    \node[below] (title) at (current bounding box.south)
    {\(\myApp{B}{x_{0}}\)};
  \end{tikzpicture}
  \and
  \begin{tikzpicture}
    \node (x1) at (0.5, 1.5) {\(\bullet\)};
    \node (x2) at (2.5, 1.5) {\(\bullet\)};
    \node (x4) at (1.5, 3) {\(\bullet\)};
    \draw[->,out=90,in=180] (x1) to (x4);
    \draw[->] (x1) to node[below,name=y1] {}
    node[above,name=y2] {} (x2);
    \draw[->,out=0,in=90] (x4) to (x2);
    \draw[->,out=210,in=180] (x4)
    to (1.5, 2.25)
    to node[below,name=y4] {} node[above,name=y5] {} (1.5, 2.25)
    to[out=0,in=-30] (x4);
    \draw[double equal sign distance,-Implies] (y4) to (y2);
    \draw[double equal sign distance,-Implies] (x4) to (y5);
    \node[below] (title) at (current bounding box.south)
    {\(\myApp{B}{x_{1}}\)};
  \end{tikzpicture}
  \and
  \begin{tikzpicture}
    \node (x2) at (2.5, 1.5) {\(\bullet\)};
    \node (x3) at (3, 0) {\(\bullet\)};
    \draw[->] (x2) to node[right,name=y6] {} (x3);
    \node[below] (title) at (current bounding box.south)
    {\(\myApp{B}{x_{2}}\)};
  \end{tikzpicture}
\end{mathpar}
Then the result of substitution of \(B\) in \(A\) is the following
\(2\)-pasting diagram.
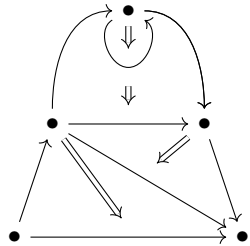
\begin{mathpar}
  \begin{tikzpicture}
    \node (x0) at (0, 0) {\(\bullet\)};
    \node (x1) at (0.5, 1.5) {\(\bullet\)};
    \node (x2) at (2.5, 1.5) {\(\bullet\)};
    \node (x3) at (3, 0) {\(\bullet\)};
    \node (x4) at (1.5, 3) {\(\bullet\)};
    \draw[->] (x0) to (x1);
    \draw[->] (x0) to node[above,name=y0] {} (x3);
    \draw[->,out=90,in=180] (x1) to (x4);
    \draw[->] (x1) to node[below,name=y1] {}
    node[above,name=y2] {} (x2);
    \draw[->] (x2) to (x3);
    \draw[->,out=0,in=90] (x4) to (x2);
    \draw[->] (x1) to node[above,name=y3] {} (x3);
    \draw[double equal sign distance,-Implies] (x1) to (y0);
    \draw[double equal sign distance,-Implies] (x2) to (y3);
    \draw[->,out=0,in=90] (x4) to (x2);
    \draw[->,out=210,in=180] (x4)
    to (1.5, 2.25)
    to node[below,name=y4] {} node[above,name=y5] {} (1.5, 2.25)
    to[out=0,in=-30] (x4);
    \draw[double equal sign distance,-Implies] (y4) to (y2);
    \draw[double equal sign distance,-Implies] (x4) to (y5);
  \end{tikzpicture}
\end{mathpar}
Grafting is the operator on \(n\)-pasting diagrams attaching an
\(n\)-pasting diagram to each leaf of an \(n\)-pasting diagram. For
example, let \(A\), \(\myApp{B}{y_{0}}\), \(\myApp{B}{y_{1}}\), and
\(\myApp{B}{y_{2}}\) be the following \(2\)-pasting diagrams.
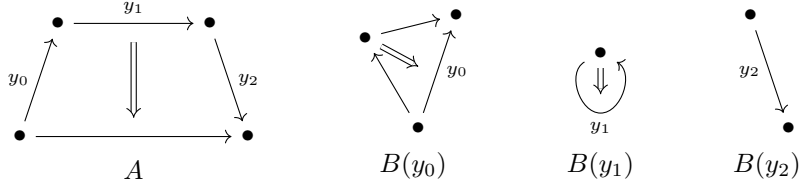
\begin{mathpar}
  \begin{tikzpicture}
    \node (x0) at (0, 0) {\(\bullet\)};
    \node (x1) at (0.5, 1.5) {\(\bullet\)};
    \node (x2) at (2.5, 1.5) {\(\bullet\)};
    \node (x3) at (3, 0) {\(\bullet\)};
    \draw[->] (x0) to node[left] {\scriptsize \(y_{0}\)} (x1);
    \draw[->] (x0) to node[above,name=y0] {} (x3);
    \draw[->] (x1) to node[below,name=y1] {}
    node[above] {\scriptsize \(y_{1}\)} (x2);
    \draw[->] (x2) to node[right] {\scriptsize \(y_{2}\)} (x3);
    \draw[double equal sign distance,-Implies] (y1) to (y0);
    \node[below] (title) at (current bounding box.south) {\(A\)};
  \end{tikzpicture}
  \and
  \begin{tikzpicture}
    \node (x0) at (0, 0) {\(\bullet\)};
    \node (x1) at (0.5, 1.5) {\(\bullet\)};
    \node (x4) at (-0.7, 1.2) {\(\bullet\)};
    \draw[->] (x0) to node[left,name=y2] {}
    node[right] {\scriptsize \(y_{0}\)} (x1);
    \draw[->] (x0) to (x4);
    \draw[->] (x4) to (x1);
    \draw[double equal sign distance,-Implies] (x4) to (y2);
    \node[below] (title) at (current bounding box.south)
    {\(\myApp{B}{y_{0}}\)};
  \end{tikzpicture}
  \and
  \begin{tikzpicture}
    \node (x1) at (0.5, 1.5) {\(\bullet\)};
    \node (x2) at (0.5, 1.5) {\(\phantom{\bullet}\)};
    \draw[->] (x1)
    to[out=-150,in=180] (0.5, 0.7)
    to node[below,name=y3] {\scriptsize \(y_{1}\)}
    node[above,name=y4] {} (0.5, 0.7)
    to[out=0,in=-30] (x2);
    \draw[double equal sign distance,-Implies] (x1) to (y4);
    \node[below] at (current bounding box.south)
    {\(\myApp{B}{y_{1}}\)};
  \end{tikzpicture}
  \and
  \begin{tikzpicture}
    \node (x2) at (2.5, 1.5) {\(\bullet\)};
    \node (x3) at (3, 0) {\(\bullet\)};
    \draw[->] (x2) to node[right,name=y5] {}
    node[left] {\scriptsize \(y_{2}\)} (x3);
    \node[below] at (current bounding box.south)
    {\(\myApp{B}{y_{2}}\)};
  \end{tikzpicture}
\end{mathpar}
Then the result of grafting of \(B\) to \(A\) is the following
\(2\)-pasting diagram.
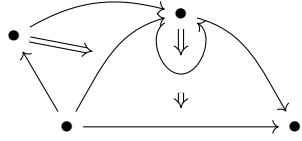
\begin{mathpar}
  \begin{tikzpicture}
    \node (x0) at (0, 0) {\(\bullet\)};
    \node (x1) at (1.5, 1.5) {\(\bullet\)};
    \node (x2) at (1.5, 1.5) {\(\phantom{\bullet}\)};
    \node (x3) at (3, 0) {\(\bullet\)};
    \node (x4) at (-0.7, 1.2) {\(\bullet\)};
    \draw[->,out=60,in=-165] (x0) to node[left,name=y2] {} (x1);
    \draw[->] (x0) to node[above,name=y0] {} (x3);
    \draw[->,out=-15,in=120] (x2) to node[right,name=y5] {} (x3);
    \draw[->] (x0) to (x4);
    \draw[->,out=30,in=165] (x4) to (x1);
    \draw[->] (x1)
    to[out=-150,in=180] (1.5, 0.7)
    to node[below,name=y3] {} node[above,name=y4] {} (1.5, 0.7)
    to[out=0,in=-30] (x2);
    \draw[double equal sign distance,-Implies] (x1) to (y4);
    \draw[double equal sign distance,-Implies] (x4) to (y2);
    \draw[double equal sign distance,-Implies] (y3) to (y0);
  \end{tikzpicture}
\end{mathpar}

We begin with the formal definition of substitution.

\begin{myNotation}
  Let \(f \myElemOf B \myMorphism A\) be a map between types. We write
  the fiber of \(f\) over \(x \myElemOf A\) as
  \(\myFiber{B}{f}{x}\).
\end{myNotation}

\begin{myConstruction}
  Let \(n \myGe 0\), let \(A\) be an \(n\)-pasting diagram, and let
  \(B \myElemOf (x \myElemOf A_{n}) \myMorphism
  \myFiber{\myNPD{n}}{\myBoundary}{\myApp{\myBoundary}{A \mySlice
      x}}\). We define an opetopic set
  \(\myApp{\mySubst}{A \myComma B}\) called the \myDefine{substitution
    of \(B\) in \(A\)} by the following pushout in \(\myOSet\).
  \begin{equation*}
    \begin{tikzcd}
      \myCoprod_{x \myElemOf A_{n}} \myApp{\myBoundary}{A \mySlice x}
      \arrow[r]
      \arrow[d]
      \arrow[dr, myPOMark] &
      A_{\myLt n}
      \arrow[d] \\
      \myCoprod_{x \myElemOf A_{n}} \myApp{B}{x}
      \arrow[r] &
      \myApp{\mySubst}{A \myComma B}
    \end{tikzcd}
  \end{equation*}
  We extend \(\myApp{\mySubst}{A \myComma B}\) to an \(n\)-prepasting
  diagram by
  \(\myApp{\myPDLeaf}{\myApp{\mySubst}{A \myComma B}} \myDefEq
  \myApp{\myPDLeaf}{A}\) and
  \(\myApp{\myPDRoot}{\myApp{\mySubst}{A \myComma B}} \myDefEq
  \myApp{\myPDRoot}{A}\).
\end{myConstruction}

\begin{myLemma}
  \label{prop-subst-degree-n}
  Let \(n \myGe 0\), let \(A\) be an \(n\)-pasting diagram, and let
  \(B \myElemOf (x \myElemOf A_{n}) \myMorphism
  \myFiber{\myNPD{n}}{\myBoundary}{\myApp{\myBoundary}{A \mySlice
      x}}\). Then
  \((x \myElemOf A_{n}) \myBinProd \myApp{B}{x}_{n} \myEquiv
  \myApp{\mySubst}{A \myComma B}_{n}\).
\end{myLemma}
\begin{proof}
  By \cref{prop-oset-fiber-colimit}.
\end{proof}

\begin{myLemma}
  \label{prop-subst-degree-n-1}
  Let \(n \myGe 1\), let \(A\) be an \(n\)-pasting diagram, and let
  \(B \myElemOf (x \myElemOf A_{n}) \myMorphism
  \myFiber{\myNPD{n}}{\myBoundary}{\myApp{\myBoundary}{A \mySlice
      x}}\).
  \begin{enumerate}
  \item
    \(\myApp{\myPDRoot}{A} \myBinCoprod ((x \myElemOf A_{n - 1})
    \myBinProd (\myApp{B}{x}_{n - 1} \mySetMinus
    \myApp{\myPDRoot}{\myApp{B}{x}})) \myEquiv \myApp{\mySubst}{A
      \myComma B}_{n - 1}\)
  \item
    \(\myApp{\myPDLeaf}{A} \myBinCoprod ((x \myElemOf A_{n - 1})
    \myBinProd (\myApp{B}{x}_{n - 1} \mySetMinus
    \myApp{\myPDLeaf}{\myApp{B}{x}})) \myEquiv \myApp{\mySubst}{A
      \myComma B}_{n - 1}\)
  \end{enumerate}
\end{myLemma}
\begin{proof}
  We prove the first claim. The second one is similarly proved. Since
  \((x \myElemOf A_{n}) \myBinProd (A \myOneSlice x) \myEquiv ((x
  \myElemOf A_{n}) \myBinProd (A \mySourceSlice x)) \myBinCoprod ((x
  \myElemOf A_{n}) \myBinProd (A \myTargetSlice x))\),
  \((x \myElemOf A_{n}) \myBinProd \myApp{B}{x}_{n - 1} \myEquiv ((x
  \myElemOf A_{n}) \myBinProd (\myApp{B}{x}_{n - 1} \mySetMinus
  \myApp{\myPDRoot}{\myApp{B}{x}})) \myBinCoprod ((x \myElemOf A_{n})
  \myBinProd \myApp{\myPDRoot}{\myApp{B}{x}})\), and
  \((x \myElemOf A_{n}) \myBinProd (A \myTargetSlice x) \myEquiv (x
  \myElemOf A_{n}) \myBinProd \myApp{\myPDRoot}{\myApp{B}{x}}\) as
  \(\myApp{\myBoundary}{A \mySlice x} \myId
  \myApp{\myBoundary}{\myApp{B}{x}}\), we have the following pushout.
  \begin{equation*}
    \begin{tikzcd}
      (x \myElemOf A_{n}) \myBinProd (A \mySourceSlice x)
      \arrow[r]
      \arrow[d]
      \arrow[dr, myPOMark] &
      (x \myElemOf A_{n}) \myBinProd (A \myOneSlice x)
      \arrow[d] \\
      (x \myElemOf A_{n}) \myBinProd (\myApp{B}{x}_{n - 1} \mySetMinus
      \myApp{\myPDRoot}{\myApp{B}{x}})
      \arrow[r] &
      (x \myElemOf A_{n}) \myBinProd \myApp{B}{x}_{n - 1}
    \end{tikzcd}
  \end{equation*}
  By \cref{prop-pd-root},
  \((x \myElemOf A_{n}) \myBinProd (A \mySourceSlice x) \myEquiv A_{n
    - 1} \mySetMinus \myApp{\myPDRoot}{A}\). We then have the
  following pushout by \cref{prop-oset-fiber-colimit}.
  \begin{equation*}
    \begin{tikzcd}
      A_{n - 1} \mySetMinus \myApp{\myPDRoot}{A}
      \arrow[r, mySub]
      \arrow[d]
      \arrow[dr, myPOMark] &
      A_{n - 1}
      \arrow[d] \\
      (x \myElemOf A_{n}) \myBinProd (\myApp{B}{x}_{n - 1} \mySetMinus
      \myApp{\myPDRoot}{\myApp{B}{x}})
      \arrow[r] &
      \myApp{\mySubst}{A \myComma B}_{n - 1}
    \end{tikzcd}
  \end{equation*}
  Since
  \(\myApp{\myPDRoot}{A} \myBinCoprod (A_{n - 1} \mySetMinus
  \myApp{\myPDRoot}{A}) \myEquiv A\), we have
  \(\myApp{\myPDRoot}{A} \myBinCoprod ((x \myElemOf A_{n}) \myBinProd
  (\myApp{B}{x}_{n - 1} \mySetMinus
  \myApp{\myPDRoot}{\myApp{B}{x}}))\).
\end{proof}

\begin{myLemma}
  \label{prop-subst-boundary}
  Let \(n \myGe 1\), let \(A\) be an \(n\)-pasting diagram, and let
  \(B \myElemOf (x \myElemOf A_{n}) \myMorphism
  \myFiber{\myNPD{n}}{\myBoundary}{\myApp{\myBoundary}{A \mySlice
      x}}\). Then
  \(\myApp{\myBoundary}{\myApp{\mySubst}{A \myComma B}} \myId
  \myApp{\myBoundary}{A}\) in the type of \(n\)-preboundaries.
\end{myLemma}
\begin{proof}
  By \cref{prop-oset-restrict-colimit},
  \(A_{\myLt n - 1} \myEquiv \myApp{\mySubst}{A \myComma B}_{\myLt n -
    1}\). Thus, the claim is true by construction.
\end{proof}

\begin{myProposition}
  \label{prop-pd-subst}
  Let \(n \myGe 0\), let \(A\) be an \(n\)-pasting diagram, and let
  \(B \myElemOf (x \myElemOf A_{n}) \myMorphism
  \myFiber{\myNPD{n}}{\myBoundary}{\myApp{\myBoundary}{A \mySlice
      x}}\). Then \(\myApp{\mySubst}{A \myComma B}\) is an
  \(n\)-pasting diagram.
\end{myProposition}
\begin{proof}
  When \(n \myId 0\), \(A_{0}\) is contractible by
  \cref{ax-pd-degree0} with center \(\myTerminalEl\), and then
  \(\myApp{B}{\myTerminalEl} \myEquiv \myApp{\mySubst}{A \myComma
    B}\). Suppose that \(n \myGe 1\). \Cref{ax-pd-finite} is by
  \cref{prop-subst-degree-n}. \Cref{ax-pd-leaf,ax-pd-root} follow from
  \cref{prop-subst-degree-n-1}. \Cref{ax-pd-degree0} is vacuously
  true. \Cref{ax-pd-target,ax-pd-source,ax-pd-path} follow from
  \cref{prop-subst-degree-n-1}. \Cref{ax-pd-boundary} is by
  \cref{prop-subst-boundary}.
\end{proof}

The substitution operator is associative in the following sense.

\begin{myProposition}
  \label{prop-subst-assoc}
  Let \(n \myGe n\), let \(A\) be an \(n\)-pasting diagram, let
  \(B \myElemOf (x \myElemOf A_{n}) \myMorphism
  \myFiber{\myNPD{n}}{\myBoundary}{\myApp{\myBoundary}{A \mySlice
      x}}\), and let
  \(C \myElemOf (x \myElemOf A_{n}) \myMorphism (y \myElemOf
  \myApp{B}{x}_{n}) \myMorphism
  \myFiber{\myNPD{n}}{\myBoundary}{\myApp{\myBoundary}{\myApp{B}{x}
      \mySlice y}}\). Then
  \begin{equation*}
    \myApp{\mySubst}{A \myComma (x \myMapsTo
      \myApp{\mySubst}{\myApp{B}{x} \myComma \myApp{C}{x}})} \myEquiv
    \myApp{\mySubst}{\myApp{\mySubst}{A \myComma B} \myComma ((x
      \myComma y) \myMapsTo \myApp{C}{x \myComma y})}.
  \end{equation*}
\end{myProposition}
\begin{proof}
  We first note that both sides of the stated equivalence are
  well-typed. The left side is well-typed by
  \cref{prop-subst-boundary}. For the right side, use
  \cref{prop-subst-degree-n} and the equivalence
  \(\myApp{B}{x} \mySlice y \myEquiv \myApp{\mySubst}{A \myComma B}
  \mySlice (x \myComma y)\) for every \(x \myElemOf A_{n}\) and
  \(y \myElemOf \myApp{B}{x}_{n}\) by
  \cref{prop-mor-oset-slice-equiv}. Let \(X\) be the following
  pushout
  \begin{equation*}
    \begin{tikzcd}[column sep = 0.5em]
      & \myCoprod_{x \myElemOf A_{n}} \myApp{\myBoundary}{A \mySlice
        x}
      \arrow[r]
      \arrow[d]
      \arrow[dr, myPOMark] &
      A_{\myLt n}
      \arrow[d] \\
      \myCoprod_{x \myElemOf A_{n}} \myCoprod_{y \myElemOf
        \myApp{B}{x}_{n}} \myApp{\myBoundary}{\myApp{B}{x} \mySlice y}
      \arrow[r]
      \arrow[d]
      \arrow[dr, myPOMark] &
      \myCoprod_{x \myElemOf A_{n}} \myApp{B}{x}_{\myLt n}
      \arrow[r]
      \arrow[d]
      \arrow[dr, myPOMark] &
      \myApp{\mySubst}{A \myComma B}_{\myLt n}
      \arrow[d] \\
      \myCoprod_{x \myElemOf A_{n}} \myCoprod_{y \myElemOf
        \myApp{B}{x}_{n}} \myApp{C}{x \myComma y}
      \arrow[r] &
      \myCoprod_{x \myElemOf A_{n}} \myApp{\mySubst}{\myApp{B}{x}
        \myComma \myApp{C}{x}}
      \arrow[r] &
      X,
    \end{tikzcd}
  \end{equation*}
  where the upper right square is a pushout by
  \cref{prop-oset-restrict-colimit}. The composite of the upper right
  and lower right pushouts exhibits \(X\) as
  \(\myApp{\mySubst}{A \myComma (x \myMapsTo
    \myApp{\mySubst}{\myApp{B}{x} \myComma \myApp{C}{x}})}\). The
  composite of the lower left and lower right pushouts exhibits \(X\)
  as
  \(\myApp{\mySubst}{\myApp{\mySubst}{A \myComma B} \myComma ((x
    \myComma y) \myMapsTo \myApp{C}{x \myComma y})}\).
\end{proof}

Any opetope of degree \(n\) can be turned into an \(n\)-pasting
diagram, which plays the role of the unit for substitution
(\cref{prop-shift-unit-a,prop-shift-unit-b}).

\begin{myConstruction}
  Let \(A\) be an opetope of degree \(n \myGe 0\). We define an
  \(n\)-prepasting diagram \(\myApp{\myShift}{A}\) as follows. The
  underlying opetopic set of \(\myApp{\myShift}{A}\) is \(A\). An
  object \(x \myElemOf A_{n - 1}\) is a leaf/root object if the arrow
  \(x \myOneMorphism \myOpeTerminal_{A}\) is a source/target arrow.
\end{myConstruction}

\begin{myLemma}
  \label{prop-shift-boundary}
  Let \(A\) be an opetope of degree \(n \myGe 0\). Then
  \(\myApp{\myBoundary}{\myApp{\myShift}{A}} \myEquiv
  \myApp{\myBoundary}{A}\).
\end{myLemma}
\begin{proof}
  By construction.
\end{proof}

\begin{myProposition}
  \label{prop-shift-pd}
  Let \(A\) be an opetope of degree \(n \myGe 0\). Then
  \(\myApp{\myShift}{A}\) is an \(n\)-pasting diagram.
\end{myProposition}
\begin{proof}
  Straightforward. For \cref{ax-pd-boundary}, use
  \cref{prop-shift-boundary}.
\end{proof}

\begin{myProposition}
  \label{prop-shift-unit-a}
  Let \(A\) be an opetope of degree \(n \myGe 0\) and let
  \(B \myElemOf
  \myFiber{\myNPD{n}}{\myBoundary}{\myApp{\myBoundary}{A}}\). Since
  \(\myApp{\myShift}{A}_{n} \myId \{\myOpeTerminal_{A}\}\), we may
  regard \(B\) as a map
  \((x \myElemOf \myApp{\myShift}{A}_{n}) \myMorphism
  \myFiber{\myNPD{n}}{\myBoundary}{\myApp{\myBoundary}{A \mySlice
      x}}\). Then
  \(\myApp{\mySubst}{\myApp{\myShift}{A} \myComma B} \myEquiv B\).
\end{myProposition}
\begin{proof}
  By construction,
  \(\myApp{\mySubst}{\myApp{\myShift}{A} \myComma B}\) is the pushout
  of the equivalence \(\myApp{\myBoundary}{A} \myEquiv A_{\myLt n}\)
  along
  \(\myApp{\myBoundary}{A} \myId \myApp{\myBoundary}{B} \myMorphism
  B\) and thus equivalent to \(B\).
\end{proof}

\begin{myProposition}
  \label{prop-shift-unit-b}
  Let \(n \myGe 0\) and let \(A\) be an \(n\)-pasting diagram. Then
  \(\myApp{\mySubst}{A \myComma (x \myMapsTo \myApp{\myShift}{A
      \mySlice x})} \myEquiv A\).
\end{myProposition}
\begin{proof}
  We first note that
  \(\myApp{\mySubst}{A \myComma (x \myMapsTo \myApp{\myShift}{A
      \mySlice x})}\) is well-typed by \cref{prop-shift-boundary}. We
  have the following commutative square.
  \begin{equation}
    \label{eq-shift-unit-b}
    \begin{tikzcd}
      \myCoprod_{x \myElemOf A_{n}} \myApp{\myBoundary}{A \mySlice x}
      \arrow[r]
      \arrow[d] &
      A_{\myLt n}
      \arrow[d] \\
      \myCoprod_{x \myElemOf A_{n}} A \mySlice x
      \arrow[r] &
      A
    \end{tikzcd}
  \end{equation}
  It suffices to show that Square \labelcref{eq-shift-unit-b} is a
  pushout. By \cref{prop-oset-fiber-colimit}, the fiber of Square
  \labelcref{eq-shift-unit-b} is
  \begin{equation*}
    \begin{tikzcd}
      \myInitial
      \arrow[r, myId]
      \arrow[d] &
      \myInitial
      \arrow[d] \\
      (x \myElemOf A_{n}) \myBinProd \{x\}
      \arrow[r, "\myEquiv"'] &
      A_{n},
    \end{tikzcd}
  \end{equation*}
  which is a pushout. By \cref{prop-oset-restrict-colimit}, the
  restriction of Square \labelcref{eq-shift-unit-b} to \(\myLt n\) is
  \begin{equation*}
    \begin{tikzcd}
      \myCoprod_{x \myElemOf A_{n}} (A \mySlice x)_{\myLt n}
      \arrow[r]
      \arrow[d, myId] &
      A_{\myLt n}
      \arrow[d, myId] \\
      \myCoprod_{x \myElemOf A_{n}} (A \mySlice x)_{\myLt n}
      \arrow[r] &
      A_{\myLt n},
    \end{tikzcd}
  \end{equation*}
  which is a pushout. Therefore, Square \labelcref{eq-shift-unit-b} is
  a pushout by \cref{prop-n-oset-colimit-fiberwise}.
\end{proof}

We then define grafting.

\begin{myConstruction}
  Let \(n \myGe 1\) and let \(A\) be an \(n\)-pasting diagram. We
  refer to the unique root object in \(A\) as \(\myPDRootObj_{A}\) and
  define an opetope \(\myApp{\myPDTargetOf}{A}\) to be
  \(A \mySlice \myPDRootObj_{A}\).
\end{myConstruction}

\begin{myConstruction}
  Let \(n \myGe 0\), let \(A\) be an \((n + 1)\)-pasting diagram, and
  let
  \(B \myElemOf (x \myElemOf \myApp{\myPDLeaf}{A}) \myMorphism
  \myFiber{\myNPD{n + 1}}{\myPDTargetOf}{A \mySlice x}\). We define an
  opetopic set \(\myApp{\myGraft}{A \myComma B}\) called the
  \myDefine{grafting of \(B\) onto \(A\)} by the following pushout in
  \(\myOSet\).
  \begin{equation*}
    \begin{tikzcd}
      \myCoprod_{x \myElemOf \myApp{\myPDLeaf}{A}} A \mySlice x
      \arrow[r]
      \arrow[d]
      \arrow[dr, myPOMark] &
      A
      \arrow[d] \\
      \myCoprod_{x \myElemOf \myApp{\myPDLeaf}{A}} \myApp{B}{x}
      \arrow[r] &
      \myApp{\myGraft}{A \myComma B}
    \end{tikzcd}
  \end{equation*}
  We extend \(\myApp{\myGraft}{A \myComma B}\) to an
  \((n + 1)\)-prepasting diagram by
  \(\myApp{\myPDLeaf}{\myApp{\myGraft}{A \myComma B}} \myDefEq (x
  \myElemOf \myApp{\myPDLeaf}{A}) \myBinProd
  \myApp{\myPDLeaf}{\myApp{B}{x}}\) and
  \(\myApp{\myPDRoot}{\myApp{\myGraft}{A \myComma B}} \myDefEq
  \myApp{\myPDRoot}{A}\).
\end{myConstruction}

\begin{myLemma}
  \label{prop-graft-degree-np1}
  Let \(n \myGe 0\), let \(A\) be an \((n + 1)\)-pasting diagram, and
  let
  \(B \myElemOf (x \myElemOf \myApp{\myPDLeaf}{A}) \myMorphism
  \myFiber{\myNPD{n + 1}}{\myPDTargetOf}{A \mySlice x}\). Then
  \(A_{n + 1} \myBinCoprod ((x \myElemOf \myApp{\myPDLeaf}{A})
  \myBinProd \myApp{B}{x}_{n + 1}) \myEquiv \myApp{\myGraft}{A
    \myComma B}_{n + 1}\).
\end{myLemma}
\begin{proof}
  By \cref{prop-oset-fiber-colimit}.
\end{proof}

\begin{myLemma}
  \label{prop-graft-degree-n}
  Let \(n \myGe 0\), let \(A\) be an \((n + 1)\)-pasting diagram, and
  let
  \(B \myElemOf (x \myElemOf \myApp{\myPDLeaf}{A}) \myMorphism
  \myFiber{\myNPD{n + 1}}{\myPDTargetOf}{A \mySlice x}\).
  \begin{enumerate}
  \item
    \((A_{n} \mySetMinus \myApp{\myPDLeaf}{A}) \myBinCoprod ((x
    \myElemOf \myApp{\myPDLeaf}{A}) \myBinProd \myApp{B}{x}_{n})
    \myEquiv \myApp{\myGraft}{A \myComma B}_{n}\)
  \item
    \(A_{n} \myBinCoprod ((x \myElemOf \myApp{\myPDLeaf}{A})
    \myBinProd (\myApp{B}{x}_{n} \mySetMinus
    \myApp{\myPDRoot}{\myApp{B}{x}})) \myEquiv \myApp{\myGraft}{A
      \myComma B}_{n}\)
  \end{enumerate}
\end{myLemma}
\begin{proof}
  By \cref{prop-oset-fiber-colimit}, we have the following pushout
  \begin{equation*}
    \begin{tikzcd}
      \myApp{\myPDLeaf}{A}
      \arrow[r, mySub]
      \arrow[d, mySub]
      \arrow[dr, myPOMark] &
      A_{n}
      \arrow[d] \\
      (x \myElemOf \myApp{\myPDLeaf}{A}) \myBinProd \myApp{B}{x}_{n}
      \arrow[r] &
      \myApp{\myGraft}{A \myComma B}_{n},
    \end{tikzcd}
  \end{equation*}
  where the left map sends \(x \myElemOf \myApp{\myPDLeaf}{A}\) to
  \((x \myComma \myPDRootObj_{\myApp{B}{x}})\). The first claim
  directly follows from this pushout. Since
  \(\myApp{\myPDLeaf}{A} \myBinCoprod ((x \myElemOf
  \myApp{\myPDLeaf}{A}) \myBinProd (\myApp{B}{x}_{n} \mySetMinus
  \{\myPDRootObj_{\myApp{B}{x}}\})) \myEquiv (x \myElemOf
  \myApp{\myPDLeaf}{A}) \myBinProd \myApp{B}{x}_{n}\), the second
  claim follows.
\end{proof}

\begin{myLemma}
  \label{prop-graft-horn}
  Let \(n \myGe 0\), let \(A\) be an \((n + 1)\)-pasting diagram, and
  let
  \(B \myElemOf (x \myElemOf \myApp{\myPDLeaf}{A}) \myMorphism
  \myFiber{\myNPD{n + 1}}{\myPDTargetOf}{A \mySlice x}\). Then
  \(\myApp{\mySourceHorn}{\myApp{\myGraft}{A \myComma B}} \myId
  \myApp{\mySubst}{\myApp{\mySourceHorn}{A} \myComma (x \myMapsTo
    \myApp{\mySourceHorn}{\myApp{B}{x}})}\) in the type of
  \(n\)-prepasting diagrams.
\end{myLemma}
\begin{proof}
  We first note that
  \(\myApp{\mySourceHorn}{A}_{n} \myId \myApp{\myPDLeaf}{A}\) and
  \(\myApp{\myBoundary}{\myApp{\mySourceHorn}{\myApp{B}{x}}} \myId
  \myApp{\myBoundary}{\myApp{\myPDTargetOf}{\myApp{B}{x}}} \myId
  \myApp{\myBoundary}{A \mySlice x} \myId
  \myApp{\myBoundary}{\myApp{\mySourceHorn}{A} \mySlice x}\) for all
  \(x \myElemOf \myApp{\mySourceHorn}{A}_{n}\), and thus
  \(\myApp{\mySubst}{\myApp{\mySourceHorn}{A} \myComma (x \myMapsTo
    \myApp{\mySourceHorn}{\myApp{B}{x}})}\) is well-typed. We have the
  following commutative square in \(\myOSet\).
  \begin{equation}
    \label{eq-graft-horn}
    \begin{tikzcd}
      \myCoprod_{x \myElemOf \myApp{\mySourceHorn}{A}_{n}}
      \myApp{\myBoundary}{\myApp{\mySourceHorn}{A} \mySlice x}
      \arrow[r]
      \arrow[d] &
      \myApp{\mySourceHorn}{A}_{\myLt n}
      \arrow[d] \\
      \myCoprod_{x \myElemOf \myApp{\mySourceHorn}{A}_{n}}
      \myApp{\mySourceHorn}{\myApp{B}{x}}
      \arrow[r] &
      \myApp{\mySourceHorn}{\myApp{\myGraft}{A \myComma B}}
    \end{tikzcd}
  \end{equation}
  It suffices to show that Square \labelcref{eq-graft-horn} is a
  pushout. By \cref{prop-oset-fiber-colimit}, the fiber of Square
  \labelcref{eq-graft-horn} over \(n\) is
  \begin{equation*}
    \begin{tikzcd}
      \myInitial
      \arrow[r, myId]
      \arrow[d] &
      \myInitial
      \arrow[d] \\
      (x \myElemOf \myApp{\myPDLeaf}{A}) \myBinProd
      \myApp{\myPDLeaf}{\myApp{B}{x}}
      \arrow[r, myId] &
      (x \myElemOf \myApp{\myPDLeaf}{A}) \myBinProd
      \myApp{\myPDLeaf}{\myApp{B}{x}},
    \end{tikzcd}
  \end{equation*}
  which is a pushout. By \cref{prop-oset-restrict-colimit}, the
  restriction of Square \labelcref{eq-graft-horn} to \(\myLt n\) is
  \begin{equation*}
    \begin{tikzcd}
      \myCoprod_{x \myElemOf \myApp{\myPDLeaf}{A}} (A \mySlice
      x)_{\myLt n}
      \arrow[r]
      \arrow[d] &
      A_{\myLt n}
      \arrow[d] \\
      \myCoprod_{x \myElemOf \myApp{\myPDLeaf}{A}} \myApp{B}{x}_{\myLt
        n}
      \arrow[r] &
      \myApp{\myGraft}{A \myComma B}_{\myLt n},
    \end{tikzcd}
  \end{equation*}
  which is a pushout by the definition of
  \(\myApp{\myGraft}{A \myComma B}\) and
  \cref{prop-oset-restrict-colimit}. Therefore, Square
  \labelcref{eq-graft-horn} is a pushout by
  \cref{prop-n-oset-colimit-fiberwise}.
\end{proof}

\begin{myProposition}
  \label{prop-graft-pd}
  Let \(n \myGe 0\), let \(A\) be an \((n + 1)\)-pasting diagram, and
  let
  \(B \myElemOf (x \myElemOf \myApp{\myPDLeaf}{A}) \myMorphism
  \myFiber{\myNPD{n + 1}}{\myPDTargetOf}{A \mySlice x}\). Then
  \(\myApp{\myGraft}{A \myComma B}\) is an \((n + 1)\)-pasting
  diagram.
\end{myProposition}
\begin{proof}
  \Cref{ax-pd-finite} is by
  \cref{prop-graft-degree-np1}. \Cref{ax-pd-leaf,ax-pd-root} follow
  from \cref{prop-graft-degree-n}. \Cref{ax-pd-degree0} is vacuously
  true. \Cref{ax-pd-target,ax-pd-source,ax-pd-path} follow from
  \cref{prop-graft-degree-n}. \Cref{ax-pd-boundary} is by
  \cref{prop-graft-horn}.
\end{proof}

The grafting operator is associative in the following sense.

\begin{myProposition}
  \label{prop-graft-assoc}
  Let \(n \myGe 0\), let \(A\) be an \((n + 1)\)-pasting diagram, let
  \(B \myElemOf (x \myElemOf \myApp{\myPDLeaf}{A}) \myMorphism
  \myFiber{\myNPD{n + 1}}{\myPDTargetOf}{A \mySlice x}\), and let
  \(C \myElemOf (x \myElemOf \myApp{\myPDLeaf}{A}) \myMorphism (y
  \myElemOf \myApp{\myPDLeaf}{\myApp{B}{x}}) \myMorphism
  \myFiber{\myNPD{n + 1}}{\myPDTargetOf}{\myApp{B}{x} \mySlice
    y}\). Then
  \begin{equation*}
    \myApp{\myGraft}{A \myComma (x \myMapsTo
      \myApp{\myGraft}{\myApp{B}{x} \myComma \myApp{C}{x}})} \myEquiv
    \myApp{\myGraft}{\myApp{\myGraft}{A \myComma B} \myComma ((x
      \myComma y) \myMapsTo \myApp{C}{x \myComma y})}.
  \end{equation*}
\end{myProposition}
\begin{proof}
  We first note that the both sides of the stated equivalence are
  well-typed. The left side is well-typed by
  \cref{prop-graft-horn}. The right side is well-typed since
  \(\myApp{B}{x} \mySlice y \myEquiv \myApp{\myGraft}{A \myComma B}
  \mySlice (x \myComma y)\) for all
  \(x \myElemOf \myApp{\myPDLeaf}{A}\) and
  \(y \myElemOf \myApp{\myPDLeaf}{\myApp{B}{x}}\) by
  \cref{prop-mor-oset-slice-equiv}. Let \(X\) be the following
  pushout.
  \begin{equation*}
    \begin{tikzcd}[column sep = 0.5em]
      & \myCoprod_{x \myElemOf \myApp{\myPDLeaf}{A}} A \mySlice x
      \arrow[r]
      \arrow[d]
      \arrow[dr, myPOMark] &
      A
      \arrow[d] \\
      \myCoprod_{x \myElemOf \myApp{\myPDLeaf}{A}} \myCoprod_{y
        \myElemOf \myApp{\myPDLeaf}{\myApp{B}{x}}} \myApp{B}{x}
      \mySlice y
      \arrow[r]
      \arrow[d]
      \arrow[dr, myPOMark] &
      \myCoprod_{x \myElemOf \myApp{\myPDLeaf}{A}} \myApp{B}{x}
      \arrow[r]
      \arrow[d]
      \arrow[dr, myPOMark] &
      \myApp{\myGraft}{A \myComma B}
      \arrow[d] \\
      \myCoprod_{x \myElemOf \myApp{\myPDLeaf}{A}} \myCoprod_{y
        \myElemOf \myApp{\myPDLeaf}{\myApp{B}{x}}} \myApp{C}{x
        \myComma y}
      \arrow[r] &
      \myCoprod_{x \myElemOf \myApp{\myPDLeaf}{A}}
      \myApp{\myGraft}{\myApp{B}{x} \myComma \myApp{C}{x}}
      \arrow[r] &
      X
    \end{tikzcd}
  \end{equation*}
  The composite of the upper right and lower right pushouts exhibits
  \(X\) as
  \(\myApp{\myGraft}{A \myComma (x \myMapsTo
    \myApp{\myGraft}{\myApp{B}{x} \myComma \myApp{C}{x}})}\). The
  composite of the lower left and lower right pushouts exhibits \(X\)
  as
  \(\myApp{\myGraft}{\myApp{\myGraft}{A \myComma B} \myComma ((x
    \myComma y) \myMapsTo \myApp{C}{x \myComma y})}\).
\end{proof}

Any opetope of degree \(n\) can be turned into an \((n + 1)\)-pasting
diagram, which plays the role of the unit for grafting
(\cref{prop-degen-unit-a,prop-degen-unit-b}).

\begin{myConstruction}
  Let \(A\) be an opetope of degree \(n \myGe 0\). We define an
  \((n + 1)\)-prepasting diagram \(\myApp{\myDegen}{A}\) called the
  \myDefine{degenerate pasting diagram on \(A\)} as follows. The
  underlying opetopic set of \(\myApp{\myDegen}{A}\) is \(A\). The
  terminal object \(\myOpeTerminal_{A} \myElemOf A_{n}\) is both a
  leaf object and a root object.
\end{myConstruction}

\begin{myLemma}
  \label{prop-degen-horn}
  Let \(A\) be an opetope of degree \(n \myGe 0\). Then
  \(\myApp{\mySourceHorn}{\myApp{\myDegen}{A}} \myEquiv A\) and
  \(\myApp{\myPDTargetOf}{\myApp{\myDegen}{A}} \myEquiv A\).
\end{myLemma}
\begin{proof}
  By construction.
\end{proof}

\begin{myProposition}
  \label{prop-degen-pd}
  Let \(A\) be an opetope of degree \(n \myGe 0\). Then
  \(\myApp{\myDegen}{A}\) is an \((n + 1)\)-pasting diagram.
\end{myProposition}
\begin{proof}
  Straightforward. For \cref{ax-pd-boundary}, use
  \cref{prop-degen-horn}.
\end{proof}

\begin{myProposition}
  \label{prop-degen-unit-a}
  Let \(A\) be an opetope of degree \(n \myGe 0\) and let
  \(B \myElemOf \myFiber{\myNPD{n + 1}}{\myPDTargetOf}{A}\). Since
  \(\myApp{\myPDLeaf}{\myApp{\myDegen}{A}} \myId
  \{\myOpeTerminal_{A}\}\), we may regard \(B\) as a map
  \((x \myElemOf \myApp{\myPDLeaf}{\myApp{\myDegen}{A}}) \myMorphism
  \myFiber{\myNPD{n + 1}}{\myPDTargetOf}{A \mySlice x}\). Then
  \(\myApp{\myGraft}{\myApp{\myDegen}{A} \myComma B} \myEquiv B\).
\end{myProposition}
\begin{proof}
  By construction, \(\myApp{\myGraft}{A \myComma B}\) is the pushout
  of the equivalence \(A \mySlice \myOpeTerminal_{A} \myEquiv A\)
  along \(A \mySlice \myOpeTerminal_{A} \myMorphism B\) and thus
  equivalent to \(B\).
\end{proof}

\begin{myProposition}
  \label{prop-degen-unit-b}
  Let \(n \myGe 0\) and let \(A\) be an \((n + 1)\)-pasting
  diagram. Then
  \(\myApp{\myGraft}{A \myComma (x \myMapsTo \myApp{\myDegen}{A
      \mySlice x})} \myEquiv A\).
\end{myProposition}
\begin{proof}
  We first note that
  \(\myApp{\myGraft}{A \myComma (x \myMapsTo \myApp{\myDegen}{A
      \mySlice x})}\) is well-typed by \cref{prop-degen-horn}. By
  construction,
  \(\myApp{\myGraft}{A \myComma (x \myMapsTo \myApp{\myDegen}{A
      \mySlice x})}\) is the pushout of the equivalence
  \(\myCoprod_{x \myElemOf \myApp{\myPDLeaf}{A}} A \mySlice x \myEquiv
  \myCoprod_{x \myElemOf \myApp{\myPDLeaf}{A}} \myApp{\myDegen}{A
    \mySlice x}\) along
  \(\myCoprod_{x \myElemOf \myApp{\myPDLeaf}{A}} A \mySlice x
  \myMorphism A\) and thus equivalent to \(A\).
\end{proof}

\section{Equivalence with existing definitions}
\label{sec:baez-dolan-constr}

We show that our definition of opetopes is equivalent to the
polynomial monad definition given by
\myTextCite{kock2010polynomial}. We also see that the category of
opetopes is presented by the generators and relations described by
\myTextCite{ho-thanh2020equivalence-arxiv}.

\subsection{Equivalence with the polynomial monad definition}
\label{sec:equiv-with-polyn}

We first review the polynomial monad definition of opetopes given by
\myTextCite{kock2010polynomial}. Let \(I\) be a set. We define
\(\myObjFinSet_{I}\) to be the type of finite sets \(E\) equipped with
a map \(E \myMorphism I\). The category \(\myPolynomial_{I}\) of
\myDefine{(finitary) polynomials on \(I\)} is defined to be
\(\mySet \mySlice (\myObjFinSet_{I} \myBinProd I)\). Concretely, a
polynomial \(P\) on \(I\) consists of a set \(\myApp{\myPolyNode}{P}\)
and two maps
\begin{equation*}
  \myObjFinSet_{I} \myXMorphismAlt{\myPolyInput_{P}}
  \myApp{\myPolyNode}{P} \myXMorphism{\myPolyTarget_{P}} I.
\end{equation*}
We regard \(\myApp{\myPolyNode}{P}\) as an object in
\(\mySet \mySlice I\) with \(\myPolyTarget_{P}\). We refer to
\(\myCoprod_{b \myElemOf \myApp{\myPolyNode}{P}}
\myApp{\myPolyInput_{P}}{b} \myElemOf \mySet \mySlice I\) as
\(\myPolySource_{P} \myElemOf \myApp{\myPolyInput}{P} \myMorphism
I\). A polynomial on \(I\) is thus equivalently presented by three
maps
\begin{equation*}
  I \myXMorphismAlt{\myPolySource_{P}} \myApp{\myPolyInput}{P}
  \myXMorphism{\myPolyProj_{P}} \myApp{\myPolyNode}{P}
  \myXMorphism{\myPolyTarget_{P}} I,
\end{equation*}
which is a more standard definition of polynomials
\myCite{gambino2013polynomial}. For a map
\(f \myElemOf I \myMorphism J\), the postcomposition with \(f\)
induces a functor \(\myPolynomial_{I} \myMorphism
\myPolynomial_{J}\). In this way the map
\(I \myMapsTo \myPolynomial_{I}\) is functorial, and let
\(\myPolyColor \myElemOf \myPolynomial \myMorphism \mySet\) denote the
cocartesian fibration corresponding to it.

For a polynomial \(P\) on \(I\), we define a functor
\(\myPolyFun_{P} \myElemOf \mySet \mySlice I \myMorphism \mySet
\mySlice I\) by
\(\myApp{\myPolyFun_{P}}{A}_{i} \myDefEq (b \myElemOf
\myApp{\myPolyNode}{P}_{i}) \myBinProd \myApp{\myArr_{\mySet \mySlice
    I}}{\myApp{\myPolyInput_{P}}{b} \myComma A}\), which is identical
to
\((b \myElemOf \myApp{\myPolyNode}{P}_{i}) \myBinProd ((e \myElemOf
\myApp{\myPolyInput_{P}}{b}) \myMorphism
A_{\myApp{\myPolySource_{P}}{e}})\). Every morphism
\(P \myMorphism Q\) of polynomials on \(I\) induces a natural
transformation \(\myPolyFun_{P} \myMMorphism \myPolyFun_{Q}\) which is
cartesian in the sense that all the naturality squares are
pullbacks. Let \(\myApp{\myEndFunCart}{I}\) denote the category of
endofunctors on \(\mySet \mySlice I\) and cartesian natural
transformations between them. The mapping
\(P \myMapsTo \myPolyFun_{P}\) defines a fully faithful functor
\(\myPolynomial_{I} \myMorphism \myApp{\myEndFunCart}{\mySet \mySlice
  I}\) \myCite[Lemma 2.15]{gambino2013polynomial}. The monoidal
structure on \(\myApp{\myEndFunCart}{\mySet \mySlice I}\) given by
composition of endofunctors restricts to a monoidal structure on
\(\myPolynomial_{I}\). The category \(\myPolyMonad_{I}\) of
\myDefine{polynomial monads on \(I\)} is defined to be the category of
monoid objects in the monoidal category \(\myPolynomial_{I}\).

Let \(P\) be a polynomial monad on a set \(I\). By a
\myDefine{\(P\)-polynomial monad} we mean an object in
\(\myPolyMonad_{I} \mySlice P\). By definition,
\(\myPolynomial_{I} \mySlice P \myEquiv \mySet \mySlice
\myApp{\myPolyNode}{P}\). We thus regard
\(\myPolyMonad_{I} \mySlice P\) as a category over
\(\mySet \mySlice \myApp{\myPolyNode}{P}\). The \myDefine{Baez-Dolan
  construction} \(P^{\myBaezDolan}\) is the polynomial monad on
\(\myApp{\myPolyNode}{P}\) whose algebras are the \(P\)-polynomial
monads. The set \(\myNKJBMOpetope{n}\) of
\myDefine{Kock-Joyal-Batanin-Mascari ({\myKJBM}) opetopes of degree
  \(n\)} and the polynomial monad \(\myNBDPoly{n}\) on
\(\myNKJBMOpetope{n}\) are inductively defined by
\(\myNKJBMOpetope{0} \myDefEq \myTerminal\),
\(\myNBDPoly{0} \myDefEq (\myObjFinSet_{\myTerminal}
\myXMorphismAlt{\myTerminal} \myTerminal \myMorphism \myTerminal)\),
\(\myNKJBMOpetope{n + 1} \myDefEq
\myApp{\myPolyNode}{\myNBDPoly{n}}\), and
\(\myNBDPoly{n + 1} \myDefEq \myNBDPoly{n}^{\myBaezDolan}\).

Let us concretely describe the structure of a polynomial monad \(P\)
on \(I\). Because the identity polynomial on \(I\) is
\(\myObjFinSet_{I} \myXMorphismAlt{\myYoneda} I
\myXMorphism{\myIdMorphism} I\), where \(\myYoneda\) denotes the
Yoneda embedding, the unit of \(P\) is a map
\(\myMonadUnit_{P} \myElemOf (i \myElemOf I) \myMorphism
\myApp{\myPolyNode}{P}_{i}\) equipped with an equivalence
\(\myMonadUnitE_{P} \myElemOf \myApp{\myYoneda}{i} \myEquiv
\myApp{\myPolyInput_{P}}{\myApp{\myMonadUnit_{P}}{i}}\). Because the
composite \(P^{2}\) is defined by
\(\myApp{\myPolyNode}{P^{2}} \myDefEq
\myApp{\myPolyFun_{P}}{\myApp{\myPolyNode}{P}}\) and
\(\myApp{\myPolyInput_{P^{2}}}{b_{1} \myComma b_{2}} \myDefEq \myCoprod_{e
  \myElemOf \myApp{\myPolyInput_{P}}{b_{1}}}
\myApp{\myPolyInput_{P}}{\myApp{b_{2}}{\myApp{\myPolySource_{P}}{e}}}\),
the multiplication of \(P\) is a map
\(\myMonadMul_{P} \myElemOf
\myApp{\myPolyFun_{P}}{\myApp{\myPolyNode}{P}} \myMorphism
\myApp{\myPolyNode}{P}\) over \(I\) equipped with an equivalence
\(\myMonadMulE_{P} \myElemOf (\myCoprod_{e \myElemOf
  \myApp{\myPolyInput_{P}}{b_{1}}}
\myApp{\myPolyInput_{P}}{\myApp{b_{2}}{\myApp{\myPolySource_{P}}{e}}})
\myEquiv \myApp{\myPolyInput_{P}}{\myApp{\myMonadMul_{P}}{b_{1}
    \myComma b_{2}}}\). Let
\(A \myElemOf \mySet \mySlice \myApp{\myPolyNode}{P}\). The polynomial
\(Q\) over \(P\) corresponding to \(A\) is the composite
\(A \myMorphism \myApp{\myPolyNode}{P} \myXMorphism{(\myPolyInput_{P}
  \myComma \myPolyTarget_{P})} \myObjFinSet_{I} \myBinProd I\), and
\(\myApp{\myPolyFun_{Q}}{X}_{i} \myId ((b_{1} \myComma b_{2})
\myElemOf \myApp{\myPolyFun_{P}}{X}_{i}) \myBinProd A_{b_{1}}\). A
\(P\)-polynomial monad structure on \(A\) thus consists of a map
\(\myMonadUnit_{A} \myElemOf (i \myElemOf I) \myMorphism
A_{\myApp{\myMonadUnit_{P}}{i}}\) and a map
\(\myMonadMul_{A} \myElemOf \myImplicit{(b_{1} \myComma b_{2})
  \myElemOf \myApp{\myPolyFun_{P}}{\myApp{\myPolyNode}{P}}}
\myMorphism A_{b_{1}} \myMorphism ((e \myElemOf
\myApp{\myPolyInput_{P}}{b_{1}}) \myMorphism A_{\myApp{b_{2}}{e}})
\myMorphism A_{\myApp{\myMonadMul_{P}}{b_{1} \myComma b_{2}}}\)
satisfying suitable associativity and unit laws.

We also recall the notion of a \(P\)-tree. For a polynomial \(P\), we
define a graph \(\myApp{\myPolyGraph}{P}\) as follows. The set of
vertices in \(\myApp{\myPolyGraph}{P}\) is
\(\myApp{\myPolyColor}{P} \myBinCoprod \myApp{\myPolyNode}{P}\). There
is no edge between vertices from \(\myApp{\myPolyColor}{P}\). There is
no edge between vertices from \(\myApp{\myPolyNode}{P}\). An edge from
\(x \myElemOf \myApp{\myPolyColor}{P}\) to
\(y \myElemOf \myApp{\myPolyNode}{P}\) is an element
\(e \myElemOf \myApp{\myPolyInput_{P}}{y}_{x}\). An edge from
\(y \myElemOf \myApp{\myPolyNode}{P}\) to
\(x \myElemOf \myApp{\myPolyColor}{P}\) is an identification
\(\myApp{\myPolyTarget_{P}}{y} \myId x\). A \myDefine{polynomial tree}
is a polynomial \(P\) satisfying the following axioms.
\begin{myWithCustomLabel}{axiom}
  \begin{enumerate}[label=\textbf{PT\arabic*.},ref=PT\arabic*]
  \item \label{ax-polytree-finite} The sets
    \(\myApp{\myPolyColor}{P}\) and \(\myApp{\myPolyNode}{P}\) are
    finite.
  \item \label{ax-polytree-decidable} The maps \(\myPolyTarget_{P}\)
    and \(\myPolySource_{P}\) are injective.
  \item \label{ax-polytree-tree} The graph
    \(\myApp{\myPolyGraph}{P}\) is a tree.
  \end{enumerate}
\end{myWithCustomLabel}
Note that the image of any map between finite sets is decidable. For a
polynomial tree \(P\), let \(\myApp{\myPTrLeaf}{P}\) denote the
complement of the image of \(\myPolyTarget_{P}\) whose elements are
called \myDefine{leaves in \(P\)}. The complement of the image of
\(\myPolySource_{P}\) is the singleton consisting of the root in the
tree \(\myApp{\myPolyGraph}{P}\) by \cref{ax-polytree-tree}. We refer
to the root in \(\myApp{\myPolyGraph}{P}\) as \(\myPTrRoot_{P}\) and
called it the \myDefine{root in \(P\)}. For a polynomial \(P\), a
\myDefine{\(P\)-tree} is a polynomial tree \(T\) equipped with a
morphism \(\myPTrDeco_{T} \myElemOf T \myMorphism P\) in
\(\myPolynomial\). Let \(\myApp{\myPTree}{P}\) denote the category of
\(P\)-trees whose morphisms are those morphisms of polynomials over
\(P\) preserving roots and leaves.

The Baez-Dolan construction \(P^{\myBaezDolan}\) has an explicit
construction using \(P\)-trees. There is an equivalence
\(\myPTreeToBD \myElemOf \myApp{\myObj}{\myApp{\myPTree}{P}} \myEquiv
\myApp{\myPolyNode}{P^{\myBaezDolan}}\) characterized as follows. By
the definition of \(P^{\myBaezDolan}\), the object
\(\myApp{\myPolyNode}{P^{\myBaezDolan}} \myEquiv
\myApp{\myPolyFun_{P^{\myBaezDolan}}}{\myTerminal} \myElemOf \mySet
\mySlice \myApp{\myPolyNode}{P}\) is the free \(P\)-polynomial monad
over \(\myTerminal\). For a \(P\)-tree \(T\) and
\(i \myElemOf \myApp{\myPolyColor}{T}\), we define a polynomial
\(T \myStarSlice x\) as follows.
\(\myApp{\myPolyColor}{T \myStarSlice i}\) is the subset of
\(\myApp{\myPolyColor}{T}\) spanned by those \(i'\) such that there is
a (unique) path in \(\myApp{\myPolyGraph}{T}\) from \(i'\) to
\(i\). We define
\(\myApp{\myPolyNode}{T \myStarSlice i} \myDefEq
\myApp{\myPolyNode}{T} \myBinProd_{\myApp{\myPolyColor}{T}}
\myApp{\myPolyColor}{T \myStarSlice i}\). The map
\(\myPolySource_{T} \myElemOf \myApp{\myPolyInput_{T}}{b} \myMorphism
\myApp{\myPolyColor}{T}\) factors through
\(\myApp{\myPolyColor}{T \myStarSlice i}\) when
\(b \myElemOf \myApp{\myPolyNode}{T \myStarSlice i}\), and thus we can
define
\(\myApp{\myPolyInput_{T \myStarSlice i}}{b} \myDefEq
\myApp{\myPolyInput_{T}}{b}\). One can show that \(T \myStarSlice i\)
is a polynomial tree. By construction, we have a morphism
\(T \myStarSlice i \myMorphism T \myXMorphism{\myPTrDeco_{T}} P\) in
\(\myPolynomial\) by which we regard \(T \myStarSlice i\) as a
\(P\)-tree. Then \(\myApp{\myPTreeToBD}{T}\) for a \(P\)-tree \(T\) is
defined by induction on the size of \(\myApp{\myPolyColor}{T}\) as
follows.
\begin{itemize}
\item If the root \(\myPTrRoot_{T}\) is a leaf, then
  \(\myApp{\myPTreeToBD}{T} \myId
  \myApp{\myMonadUnit_{\myApp{\myPolyNode}{P^{\myBaezDolan}}}}{\myApp{\myPTrDeco_{T}}{\myPTrRoot_{T}}}\).
\item If there is a (unique)
  \(b \myElemOf \myApp{\myPolyNode}{T}_{\myPTrRoot_{T}}\), then
  \(\myApp{\myPTreeToBD}{T} \myId
  \myApp{\myMonadMul_{\myApp{\myPolyNode}{P^{\myBaezDolan}}}}{\myApp{\myMonadUnit_{P^{\myBaezDolan}}}{\myApp{\myPTrDeco_{T}}{b}}
    \myComma (e \myMapsTo \myApp{\myPTreeToBD}{T \myStarSlice
      \myApp{\myPolySource_{T}}{e}})}\).
\end{itemize}

We now prove the equivalence of opetopes in our sense and {\myKJBM}
opetopes, that is, \(\myOpetope_{n} \myEquiv \myNKJBMOpetope{n}\)
(\cref{prop-ope-equiv-kock}). We first make \(\myOpetope_{n}\)'s part
of polynomial monads.

\begin{myConstruction}
  Let \(A\) be an opetopic set. For \(n \myGe 0\), we define a
  polynomial \(\myApp{\myNOpePoly{n}}{A}\) to be
  \begin{equation*}
    \myObjFinSet_{A_{n}}
    \myXMorphismAlt{\myPolyInput_{\myApp{\myNOpePoly{n}}{A}}} A_{n +
      1} \myXMorphism{\myDomTargetOf_{A}} A_{n},
  \end{equation*}
  where
  \(\myApp{\myPolyInput_{\myApp{\myNOpePoly{n}}{A}}}{x} \myDefEq A
  \mySourceSlice x\), which is finite by \cref{ax-oset-finite}. A
  morphism \(F \myElemOf A \myMorphism B\) of opetopic sets induces a
  morphism of polynomials
  \(\myApp{\myNOpePoly{n}}{A} \myMorphism \myApp{\myNOpePoly{n}}{B}\)
  because
  \(A \mySourceSlice x \myEquiv B \mySourceSlice \myApp{F}{x}\) by
  \cref{prop-mor-oset-slice-equiv}. In particular,
  \(\myApp{\myNOpePoly{n}}{A} \myElemOf \myPolynomial \mySlice
  \myApp{\myNOpePoly{n}}{\myOpetope}\) by
  \cref{prop-opetope-terminal}.
\end{myConstruction}

\begin{myConstruction}
  Let \(n \myGe 0\). We define a polynomial \(\myNOpePolyAlt{n}\) to
  be
  \begin{equation*}
    \myObjFinSet_{\myNBd{n}}
    \myXMorphismAlt{\myPolyInput_{\myNOpePolyAlt{n}}} \myNPD{n}
    \myXMorphism{\myBoundary} \myNBd{n},
  \end{equation*}
  where
  \(\myApp{\myPolyInput_{\myNOpePolyAlt{n}}}{A} \myDefEq A_{n}\),
  which is finite by \cref{ax-pd-finite}, with
  \((x \myMapsTo \myApp{\myBoundary}{A \mySlice x}) \myElemOf A_{n}
  \myMorphism \myNBd{n}\). We extend \(\myNOpePolyAlt{n}\) to a
  polynomial monad on \(\myNBd{n}\) as follows. We define
  \(\myApp{\myMonadUnit_{\myNOpePolyAlt{n}}}{A} \myDefEq
  \myApp{\myShift}{\myApp{\myFill}{A}}\) and
  \(\myApp{\myMonadMul_{\myNOpePolyAlt{n}}}{A \myComma B} \myDefEq
  \myApp{\mySubst}{A \myComma B}\). We have
  \(\myApp{\myMonadUnit_{\myNOpePolyAlt{n}}}{A}_{n} \myEquiv
  \{\myOpeTerminal\}\) and
  \(\myApp{\myBoundary}{\myApp{\myMonadUnit_{\myNOpePolyAlt{n}}}{A}
    \mySlice \myOpeTerminal} \myEquiv A\) by construction, and thus
  \(\myApp{\myYoneda}{A} \myEquiv
  \myApp{\myPolyInput_{\myNOpePolyAlt{n}}}{\myApp{\myMonadUnit_{\myNOpePolyAlt{n}}}{A}}\). By
  \cref{prop-subst-degree-n},
  \(\myCoprod_{x \myElemOf
    \myApp{\myPolyInput_{\myNOpePolyAlt{n}}}{A}}
  \myApp{\myPolyInput_{\myNOpePolyAlt{n}}}{\myApp{B}{x}} \myEquiv
  \myApp{\myPolyInput_{\myNOpePolyAlt{n}}}{\myApp{\myMonadMul_{\myNOpePolyAlt{n}}}{A
      \myComma B}}\). The associativity and unit laws follow from
  \cref{prop-subst-assoc,prop-shift-unit-a,prop-shift-unit-b}. By
  \cref{prop-opetope-equiv-pd,prop-boundary-equiv},
  \(\myNOpePolyAlt{n} \myEquiv \myApp{\myNOpePoly{n}}{\myOpetope}\),
  so the polynomial monad structure on \(\myNOpePolyAlt{n}\) is
  transported to \(\myApp{\myNOpePoly{n}}{\myOpetope}\). That is, the
  polynomial monad structure on \(\myApp{\myNOpePoly{n}}{\myOpetope}\)
  is determined by
  \(\myApp{\mySourceHorn}{\myApp{\myMonadUnit_{\myApp{\myNOpePoly{n}}{\myOpetope}}}{A}}
  \myEquiv \myApp{\myShift}{A}\) and
  \(\myApp{\mySourceHorn}{\myApp{\myMonadMul_{\myApp{\myNOpePoly{n}}{\myOpetope}}}{A
      \myComma B}} \myEquiv \myApp{\mySubst}{\myApp{\mySourceHorn}{A}}
  \myComma (x \myMapsTo \myApp{\mySourceHorn}{\myApp{B}{x}})\).
\end{myConstruction}

It suffices to construct an equivalence
\(\myApp{\myNOpePoly{n}}{\myOpetope} \myEquiv \myNBDPoly{n}\)
(\cref{prop-ope-poly-equiv-baez-dolan}). We proceed by induction on
\(n \myElemOf \myNat\). The base case is easy.

\begin{myLemma}
  \label{prop-ope-poly-degree0}
  \(\myApp{\myNOpePoly{0}}{\myOpetope} \myEquiv \myNBDPoly{0}\).
\end{myLemma}
\begin{proof}
  By \cref{prop-opetope-degree0,prop-opetope-degree1}.
\end{proof}

For the successor case, we show that
\(\myApp{\myNOpePoly{n + 1}}{\myOpetope} \myEquiv
\myApp{\myNOpePoly{n}}{\myOpetope}^{\myBaezDolan}\)
(\cref{prop-ope-poly-bd}). By the definition of the Baez-Dolan
construction, to get a morphism
\(\myApp{\myNOpePoly{n}}{\myOpetope}^{\myBaezDolan} \myMorphism
\myApp{\myNOpePoly{n + 1}}{\myOpetope}\) of polynomial monads on
\(\myOpetope_{n + 1}\), it suffices to construct a functor
\(\myApp{\myAlg}{\myApp{\myNOpePoly{n + 1}}{\myOpetope}} \myMorphism
\myPolyMonad_{\myOpetope_{n}} \mySlice
\myApp{\myNOpePoly{n}}{\myOpetope}\) over
\(\mySet \mySlice \myOpetope_{n + 1}\), where
\(\myApp{\myAlg}{\myApp{\myNOpePoly{n + 1}}{\myOpetope}}\) is the
category of \(\myApp{\myNOpePoly{n + 1}}{\myOpetope}\)-algebras.

\begin{myConstruction}
  Let \(n \myGe 0\) and let \(A\) be a
  \(\myApp{\myNOpePoly{n + 1}}{\myOpetope}\)-algebra. That is,
  \(A \myElemOf \mySet \mySlice \myOpetope_{n + 1}\) is equipped with
  an operator
  \(\myAlgOp_{A} \myElemOf (X \myElemOf \myOpetope_{n + 2}) (a
  \myElemOf (x \myElemOf X \mySourceSlice \myOpeTerminal_{X})
  \myMorphism A_{X \mySlice x}) \myMorphism
  A_{\myApp{\myDomTargetOf}{X}}\) compatible with the polynomial monad
  structure on \(\myApp{\myNOpePoly{n + 1}}{\myOpetope}\). By
  \cref{prop-opetope-equiv-pd}, \(\myAlgOp_{A}\) is also regarded as
  an operator
  \((X \myElemOf \myNPD{n + 1}) (a \myElemOf (x \myElemOf X_{n + 1}))
  \myMorphism A_{\myApp{\myFill}{\myApp{\myBoundary}{X}}}\). We equip
  \(A\) with a \(\myApp{\myNOpePoly{n}}{\myOpetope}\)-polynomial monad
  structure as follows. For \(X \myElemOf \myOpetope_{n}\), we have
  \(\myApp{\myDegen}{X}_{n + 1} \myEquiv \myInitial\) by construction,
  so let
  \(\myApp{\myMonadUnit_{A}}{X} \myDefEq
  \myApp{\myAlgOp_{A}}{\myApp{\myDegen}{X} \myComma \myBang}\), where
  \(\myBang\) is the unique map from \(\myInitial\). Since
  \(\myApp{\mySourceHorn}{\myApp{\myDegen}{X}} \myEquiv
  \myApp{\myShift}{X}\) by construction, we see that
  \(\myApp{\myMonadUnit_{A}}{X}\) lies over
  \(\myApp{\myMonadUnit_{\myApp{\myNOpePoly{n}}{\myOpetope}}}{X}\). For
  \(X \myElemOf \myOpetope_{n + 1}\),
  \(X' \myElemOf (x \myElemOf X \mySourceSlice \myOpeTerminal_{X})
  \myMorphism \myFiber{\myOpetope_{n + 1}}{\myDomTargetOf}{X \mySlice
    x}\), \(a \myElemOf A_{X}\), and
  \(a' \myElemOf (x \myElemOf X \mySourceSlice \myOpeTerminal_{X})
  \myMorphism A_{\myApp{X'}{x}}\), we have
  \((\{\myOpeTerminal_{X}\} \myBinCoprod (\myCoprod_{x \myElemOf X
    \mySourceSlice \myOpeTerminal_{X}} \{\myOpeTerminal_{X'}\}))
  \myEquiv \myApp{\myGraft}{\myApp{\myShift}{X} \myComma (x \myMapsTo
    \myApp{\myShift}{\myApp{X'}{x}})}_{n + 1}\) by
  \cref{prop-graft-degree-np1}, and thus \(a\) and \(a'\) defines a
  map
  \((a \myComma a') \myElemOf \myApp{\myGraft}{\myApp{\myShift}{X}
    \myComma (x \myMapsTo \myApp{\myShift}{\myApp{X'}{x}})}_{n + 1}
  \myMorphism A\). We then define
  \(\myApp{\myMonadMul_{A}}{a \myComma a'} \myDefEq
  \myApp{\myAlgOp_{A}}{\myApp{\myGraft}{\myApp{\myShift}{X} \myComma
      (x \myMapsTo \myApp{\myShift}{\myApp{X'}{x}})} \myComma (a
    \myComma a')}\). Since
  \(\myApp{\mySourceHorn}{\myApp{\myGraft}{\myApp{\myShift}{X}
      \myComma (x \myMapsTo \myApp{\myShift}{\myApp{X'}{x}})}}
  \myEquiv \myApp{\mySubst}{\myApp{\mySourceHorn}{X} \myComma (x
    \myMapsTo \myApp{\mySourceHorn}{\myApp{X'}{x}})}\) by
  \cref{prop-graft-horn}, we see that
  \(\myApp{\myMonadMul_{A}}{a \myComma a'}\) lies over
  \(\myApp{\myMonadMul_{\myApp{\myNOpePoly{n}}{\myOpetope}}}{X
    \myComma X'}\). The associativity and unit laws follow from
  \cref{prop-graft-assoc,prop-degen-unit-a,prop-degen-unit-b}. This
  construction extends to a functor
  \(\myApp{\myAlg}{\myApp{\myNOpePoly{n + 1}}{\myOpetope}} \myMorphism
  \myPolyMonad_{\myOpetope_{n}} \mySlice
  \myApp{\myNOpePoly{n}}{\myOpetope}\) over
  \(\mySet \mySlice \myOpetope_{n + 1}\), and let
  \(\myAlgToPM \myElemOf
  \myApp{\myNOpePoly{n}}{\myOpetope}^{\myBaezDolan} \myMorphism
  \myApp{\myNOpePoly{n + 1}}{\myOpetope}\) be the corresponding
  morphism of polynomial monads on \(\myOpetope_{n + 1}\).
\end{myConstruction}

To see that \(\myAlgToPM\) is an equivalence, it suffices to show that
the composite
\begin{equation}
  \label{eq-tree-to-pd}
  \myApp{\myObj}{\myApp{\myPTree}{\myApp{\myNOpePoly{n}}{\myOpetope}}}
  \myEquiv
  \myApp{\myPolyNode}{\myApp{\myNOpePoly{n}}{\myOpetope}^{\myBaezDolan}}
  \myXMorphism{\myAlgToPM} \myApp{\myPolyNode}{\myApp{\myNOpePoly{n +
        1}}{\myOpetope}} \myId \myOpetope_{n + 2} \myEquiv \myNPD{n +
    1}
\end{equation}
is an equivalence. We construct an equivalence
\(\myApp{\myObj}{\myApp{\myPTree}{\myApp{\myNOpePoly{n}}{\myOpetope}}}
\myEquiv \myNPD{n + 1}\) (\cref{prop-pd-equiv-ptree}) and see that it
coincides with \cref{eq-tree-to-pd}.

\begin{myLemma}
  \label{prop-pd-poly-tree}
  Let \(n \myGe 0\) and let \(A\) be an \((n + 1)\)-pasting
  diagram. Then \(\myApp{\myNOpePoly{n}}{A}\) is a
  \(\myApp{\myNOpePoly{n}}{\myOpetope}\)-tree.
\end{myLemma}
\begin{proof}
  \Cref{ax-polytree-finite} is by
  \cref{ax-pd-finite,prop-pd-finite,ax-oset-finite}. \Cref{ax-polytree-decidable}
  follows from \cref{prop-pd-leaf,prop-pd-root}. Since
  \(\myApp{\myPolyGraph}{\myApp{\myNOpePoly{n}}{A}} \myEquiv
  \myApp{\myPDGraph}{A}\) by construction, \cref{ax-polytree-tree}
  follows from \cref{prop-pd-tree}.
\end{proof}

\begin{myConstruction}
  Let \(n \myGe 0\) and let \(T\) be a
  \(\myApp{\myNOpePoly{n}}{\myOpetope}\)-tree. We construct a graph
  \(\myApp{\myPolyGraph}{T}'\) from \(\myApp{\myPolyGraph}{T}\) by
  reversing the directions of the edges
  \(y \myMorphism \myApp{\myPolyTarget_{T}}{y}\) for all
  \(y \myElemOf \myApp{\myPolyNode}{T}\). We define a diagram
  \(\myPTrDiagram_{T} \myElemOf \myApp{\myPolyGraph}{T}' \myMorphism
  \myOSet\) as follows. A vertex \(x\) in \(\myApp{\myPolyGraph}{T}'\)
  is either in \(\myApp{\myPolyColor}{T}\) or
  \(\myApp{\myPolyNode}{T}\). In both cases,
  \(\myApp{\myPTrDiagram_{T}}{x} \myDefEq \myApp{\myPTrDeco_{T}}{x}\)
  defines an opetope. For an edge of the form
  \(\myApp{\myPolySource}{e} \myMorphism y\) for
  \(e \myElemOf \myApp{\myPolyInput_{T}}{y}\), let \(\myPTrDeco_{T}\)
  send it to the source morphism
  \(\myApp{\myPTrDeco_{T}}{e} \myElemOf
  \myApp{\myPTrDeco_{T}}{\myApp{\myPolySource}{e}} \myEquiv
  \myApp{\myPolySource}{e'} \myMorphism \myApp{\myPTrDeco_{T}}{y}\),
  where \(e'\) is the element corresponding to \(e\) via the
  equivalence
  \(\myApp{\myPolyInput_{\myApp{\myNOpePoly{n}}{\myOpetope}}}{\myApp{\myPTrDeco_{T}}{y}}
  \myEquiv \myApp{\myPolyInput_{T}}{y}\). For an edge of the form
  \(\myApp{\myPolyTarget_{T}}{y} \myMorphism y\) for
  \(y \myElemOf \myApp{\myPolyNode}{T}\), let \(\myPTrDeco_{T}\) send
  it to the target morphism
  \(\myApp{\myPTrDeco_{T}}{\myApp{\myPolyTarget_{T}}{y}} \myEquiv
  \myApp{\myPolyTarget_{\myApp{\myNOpePoly{n}}{\myOpetope}}}{\myApp{\myPTrDeco_{T}}{y}}
  \myMorphism \myApp{\myPTrDeco_{T}}{y}\). Finally, we define
  \(\myApp{\myPTrPaste}{T}\) to be the colimit of
  \(\myPTrDiagram_{T}\). We further extend \(\myApp{\myPTrPaste}{T}\)
  to an \((n + 1)\)-prepasting diagram by
  \(\myApp{\myPDLeaf}{\myApp{\myPTrPaste}{T}} \myDefEq \myCoprod_{x
    \myElemOf \myApp{\myPTrLeaf}{T}}
  \myApp{\myPDLeaf}{\myApp{\myPTrDiagram_{T}}{x}}\) and
  \(\myApp{\myPDRoot}{\myApp{\myPTrPaste}{T}} \myDefEq
  \myApp{\myPDRoot}{\myApp{\myPTrDiagram_{T}}{\myPTrRoot_{T}}}\).
\end{myConstruction}

\begin{myLemma}
  \label{prop-ptree-paste-computation}
  \label{prop-ptree-paste-pd}
  Let \(n \myGe 0\) and let \(T\) be a
  \(\myApp{\myNOpePoly{n}}{\myOpetope}\)-tree. Then
  \(\myApp{\myPTrPaste}{T}\) is an \((n + 1)\)-pasting
  diagram. Moreover, the following hold.
  \begin{enumerate}
  \item Suppose that the root \(\myPTrRoot_{T}\) is a leaf. Then
    \(\myApp{\myDegen}{\myApp{\myPTrDeco_{T}}{\myPTrRoot_{T}}}
    \myEquiv \myApp{\myPTrPaste}{T}\).
  \item Suppose that there is a (unique)
    \(b \myElemOf \myApp{\myPolyNode}{T}_{\myPTrRoot_{T}}\). Then
    \(\myApp{\myGraft}{\myApp{\myShift}{\myApp{\myPTrDeco_{T}}{b}}
      \myComma (x \myMapsTo \myApp{\myPTrPaste}{T \myStarSlice
        \myApp{\myPolySource_{T}}{x}})} \myEquiv
    \myApp{\myPTrPaste}{T}\), where we identify
    \(\myApp{\myPDLeaf}{\myApp{\myShift}{\myApp{\myPTrDeco_{T}}{b}}}
    \myEquiv \myApp{\myPTrDeco_{T}}{b} \mySourceSlice
    \myOpeTerminal_{\myApp{\myPTrDeco_{T}}{b}} \myEquiv
    \myApp{\myPolyInput_{T}}{b}\).
  \end{enumerate}
\end{myLemma}
\begin{proof}
  We proceed by induction on the size of
  \(\myApp{\myPolyColor}{T}\). Suppose that \(\myPTrRoot_{T}\) is a
  leaf. Then \(\myApp{\myPolyGraph}{T}'\) is the singleton
  \(\{\myPTrRoot_{T}\}\) with no edge. Then
  \(\myApp{\myDegen}{\myApp{\myPTrDeco_{T}}{\myPTrRoot_{T}}} \myEquiv
  \myApp{\myPTrPaste}{T}\), and thus \(\myApp{\myPTrPaste}{T}\) is an
  \((n + 1)\)-pasting diagram. Suppose that there is a (unique)
  \(b \myElemOf \myApp{\myPolyNode}{T}_{\myPTrRoot_{T}}\). By
  induction hypothesis,
  \(\myApp{\myPTrPaste}{T \myStarSlice \myApp{\myPolySource_{T}}{x}}\)
  is an \((n + 1)\)-pasting diagram for every
  \(x \myElemOf \myApp{\myPolyInput_{T}}{b}\). Observe that
  \(\myApp{\myPolyGraph}{T}'\) is the following pushout in the
  category of graphs
  \begin{equation*}
    \begin{tikzcd}
      \myCoprod_{e \myElemOf \myApp{\myPolyInput_{T}}{b}}
      \{\myApp{\myPolySource_{T}}{e}\}
      \arrow[r]
      \arrow[d]
      \arrow[dr, myPOMark] &
      X
      \arrow[d]\\
      \myCoprod_{e \myElemOf \myApp{\myPolyInput_{T}}{b}}
      \myApp{\myPolyGraph}{T \myStarSlice
        \myApp{\myPolySource_{T}}{e}}'
      \arrow[r] &
      \myApp{\myPolyGraph}{T}',
    \end{tikzcd}
  \end{equation*}
  where \(X\) is the full subgraph spanned by \(\myPTrRoot_{T}\),
  \(b\), and \(\myApp{\myPolySource_{T}}{e}\) for all
  \(e \myElemOf \myApp{\myPolyInput_{T}}{b}\). Then
  \(\myApp{\myPTrPaste}{T}\) is the following pushout in \(\myOSet\)
  \begin{equation*}
    \begin{tikzcd}
      \myCoprod_{x \myElemOf
        \myApp{\myPDLeaf}{\myApp{\myShift}{\myApp{\myPTrDeco_{T}}{b}}}}
      \myApp{\myPTrDeco_{T}}{b} \mySlice x
      \arrow[r]
      \arrow[d]
      \arrow[dr, myPOMark] &
      \myApp{\myPTrDeco_{T}}{b}
      \arrow[d] \\
      \myCoprod_{x \myElemOf
        \myApp{\myPDLeaf}{\myApp{\myShift}{\myApp{\myPTrDeco_{T}}{b}}}}
      \myApp{\myPTrPaste}{T \myStarSlice \myApp{\myPolySource_{T}}{x}}
      \arrow[r] &
      \myApp{\myPTrPaste}{T},
    \end{tikzcd}
  \end{equation*}
  where the colimit of the restriction of \(\myPTrDiagram_{T}\) to
  \(X\) is \(\myApp{\myPTrDeco_{T}}{b}\) because \(b\) is the terminal
  object in \(X\). Therefore,
  \(\myApp{\myGraft}{\myApp{\myShift}{\myApp{\myPTrDeco_{T}}{b}}
    \myComma (x \myMapsTo \myApp{\myPTrPaste}{T \myStarSlice
      \myApp{\myPolySource_{T}}{x}})} \myEquiv
  \myApp{\myPTrPaste}{T}\) by the definition of \(\myGraft\), and thus
  \(\myApp{\myPTrPaste}{T}\) is an \((n + 1)\)-pasting diagram.
\end{proof}

\begin{myLemma}
  \label{prop-ptree-groupoid}
  Let \(P\) be a polynomial. Then all the morphisms in
  \(\myApp{\myPTree}{P}\) are equivalences.
\end{myLemma}
\begin{proof}
  Let \(h \myElemOf T \myMorphism T'\) be a morphism of
  \(P\)-trees. Since \(h\) preserves leaves, we see that
  \(\myApp{\myPolyNode}{T} \myEquiv \myApp{\myPolyColor}{T}
  \myBinProd_{\myApp{\myPolyColor}{T'}}
  \myApp{\myPolyNode}{T'}\). Thus, it suffices to show that
  \(h_{\myPolyColor} \myElemOf \myApp{\myPolyColor}{T} \myMorphism
  \myApp{\myPolyColor}{T'}\) is an equivalence. We show that the fiber
  of \(h_{\myPolyColor}\) over
  \(i' \myElemOf \myApp{\myPolyColor}{T'}\) is contractible by
  induction on the length of the path in \(\myApp{\myPolyGraph}{T'}\)
  from \(i'\) to the root. If \(i'\) is the root, then the root in
  \(T\) is the unique element of the fiber of \(h_{\myPolyColor}\)
  over \(i'\). Suppose that there is a (unique) pair
  \((b' \myComma e')\) of \(b' \myElemOf \myApp{\myPolyNode}{T'}\) and
  \(e' \myElemOf \myApp{\myPolyInput_{T'}}{b'}_{i'}\). Since \(h\)
  preserves roots, we see that
  \(\myApp{\myPolyInput}{T} \myEquiv \myApp{\myPolyColor}{T}
  \myBinProd_{\myApp{\myPolyColor}{T'}}
  \myApp{\myPolyInput}{T'}\). Thus, it suffices to show that the fiber
  of
  \(h_{\myPolyInput} \myElemOf \myApp{\myPolyInput}{T} \myMorphism
  \myApp{\myPolyInput}{T'}\) over \((b' \myComma e')\) is
  contractible. Since
  \(\myApp{\myPolyInput}{T} \myEquiv \myApp{\myPolyNode}{T}
  \myBinProd_{\myApp{\myPolyNode}{T'}} \myApp{\myPolyInput}{T'}\) and
  \(\myApp{\myPolyNode}{T} \myEquiv \myApp{\myPolyColor}{T}
  \myBinProd_{\myApp{\myPolyColor}{T'}} \myApp{\myPolyNode}{T'}\),
  this follows from the induction hypothesis for
  \(\myApp{\myPolyTarget_{T'}}{b'} \myElemOf
  \myApp{\myPolyColor}{T'}\).
\end{proof}

\begin{myLemma}
  \label{prop-pd-equiv-ptree}
  Let \(n \myGe 0\). Then the functors
  \(\myNOpePoly{n} \myElemOf \myNPD{n + 1} \myMorphism
  \myApp{\myPTree}{\myApp{\myNOpePoly{n}}{\myOpetope}}\) induced by
  \cref{prop-pd-poly-tree} and
  \(\myPTrPaste \myElemOf
  \myApp{\myPTree}{\myApp{\myNOpePoly{n}}{\myOpetope}} \myMorphism
  \myNPD{n + 1}\) induced by \cref{prop-ptree-paste-pd} are mutually
  inverses.
\end{myLemma}
\begin{proof}
  For a \(\myApp{\myNOpePoly{n}}{\myOpetope}\)-tree \(T\), we have a
  canonical morphism
  \(T \myMorphism \myApp{\myNOpePoly{n}}{\myApp{\myPTrPaste}{T}}\) of
  polynomials over \(\myApp{\myNOpePoly{n}}{\myOpetope}\) by
  construction. This morphism preserves roots and leaves and thus is
  an equivalence by \cref{prop-ptree-groupoid}. For an
  \((n + 1)\)-pasting diagram \(A\), we have a canonical morphism
  \(\myApp{\myPTrPaste}{\myApp{\myNOpePoly{n}}{A}} \myMorphism A\) of
  opetopic sets by construction. This morphism preserves root and leaf
  objects and thus is an equivalence by \cref{prop-pd-discrete}.
\end{proof}

\begin{myLemma}
  \label{prop-ope-poly-bd}
  Let \(n \myGe 0\). Then the morphism
  \(\myAlgToPM \myElemOf
  \myApp{\myNOpePoly{n}}{\myOpetope}^{\myBaezDolan} \myMorphism
  \myApp{\myNOpePoly{n + 1}}{\myOpetope}\) of polynomial monads on
  \(\myOpetope_{n +1}\) is an equivalence.
\end{myLemma}
\begin{proof}
  By \cref{prop-ptree-paste-computation} and by the definition of the
  equivalence
  \(\myApp{\myObj}{\myApp{\myPTree}{\myApp{\myNOpePoly{n}}{\myOpetope}}}
  \myEquiv
  \myApp{\myPolyNode}{\myApp{\myNOpePoly{n}}{\myOpetope}^{\myBaezDolan}}\),
  we see that \cref{eq-tree-to-pd} is equivalent to
  \(\myPTrPaste \myElemOf
  \myApp{\myObj}{\myApp{\myPTree}{\myApp{\myNOpePoly{n}}{\myOpetope}}}
  \myMorphism \myNPD{n + 1}\), which is an equivalence by
  \cref{prop-pd-equiv-ptree}.
\end{proof}

\begin{myTheorem}
  \label{prop-ope-poly-equiv-baez-dolan}
  \(\myApp{\myNOpePoly{n}}{\myOpetope} \myEquiv \myNBDPoly{n}\) for
  all \(n \myElemOf \myNat\).
\end{myTheorem}
\begin{proof}
  By \cref{prop-ope-poly-degree0,prop-ope-poly-bd}.
\end{proof}

\begin{myCorollary}
  \label{prop-ope-equiv-kock}
  \(\myOpetope_{n} \myEquiv \myNKJBMOpetope{n}\) for all
  \(n \myElemOf \myNat\).
\end{myCorollary}
\begin{proof}
  By \cref{prop-ope-poly-equiv-baez-dolan}.
\end{proof}

\subsection{Equivalence with Ho Thanh's category of opetopes}
\label{sec:equivalence-with-ho}

We compare the canonical presentation of our category \(\myOpetope\)
of opetopes
(\cref{cst-oset-canonical-presentation,prop-oset-canonical-presentation})
and the presentation given by \myTextCite[Definition
3.6]{ho-thanh2020equivalence-arxiv}. The set of objects
\(\myApp{\myObj}{\myOpetope}\) is
\((n \myElemOf \myFinOrd) \myBinProd \myOpetope_{n}\), which coincides
with \myCite[Definition 3.6 (1)]{ho-thanh2020equivalence-arxiv} by
\cref{prop-ope-equiv-kock}. The generating morphisms for
\(\myOpetope\) are the source and target morphisms, which coincides
with \myCite[Definition 3.6 (2)]{ho-thanh2020equivalence-arxiv}. The
relations for \(\myOpetope\) are all the equations
\(f_{1} \myComp g_{1} \myId f_{2} \myComp g_{2}\) that hold in
\(\myOpetope\) such that exactly one of the following holds.
\begin{enumerate}
\item \(f_{1}\), \(f_{2}\), and \(g_{2}\) are source arrows and
  \(g_{1}\) is a target arrow. This corresponds to Eq.\@ (Inner) in
  \myCite[Definition 3.6 (3)]{ho-thanh2020equivalence-arxiv}.
\item \(f_{1}\) is a source arrow and \(g_{1}\), \(f_{2}\), and
  \(g_{2}\) are target arrows. This corresponds to Eq.\@ (Glob1) in
  \myCite[Definition 3.6 (3)]{ho-thanh2020equivalence-arxiv}.
\item \(f_{1}\) is a target arrow and \(g_{1}\), \(f_{2}\), and
  \(g_{2}\) are source arrows. This corresponds to Eq.\@ (Glob2) in
  \myCite[Definition 3.6 (3)]{ho-thanh2020equivalence-arxiv}.
\item \(f_{1}\), \(f_{2}\), and \(g_{2}\) are target arrows and
  \(g_{1}\) is a source arrow. This corresponds to Eq.\@ (Degen) in
  \myCite[Definition 3.6 (3)]{ho-thanh2020equivalence-arxiv}.
\end{enumerate}
Therefore:

\begin{myTheorem}
  \label{prop-opetope-ho-thanh}
  The category of opetopes \(\myOpetope\) is equivalent to the one
  given by \myTextCite[Definition
  3.6]{ho-thanh2020equivalence-arxiv}. \qed
\end{myTheorem}

\section*{Acknowledgements}
\label{sec:acknowledgements}

The author would like to thank Soichiro Fujii for helpful
feedback. This work was supported by JST (JPMJMS2033).

\printbibliography

\end{document}